\documentclass[11pt]{amsart}

\usepackage[english]{babel}

\usepackage[a4paper,top=2cm,bottom=2cm,left=3cm,right=3cm,marginparwidth=1.75cm]{geometry}

\usepackage{amsmath}
\usepackage{graphicx}
\usepackage{comment}
\usepackage[colorlinks=true, allcolors=blue]{hyperref}

\usepackage{amssymb,amsmath,amsthm,amscd,mathtools,mathdots}
\usepackage{multirow}
\usepackage{enumitem}
\usepackage{caption}
\usepackage{subcaption}
\usepackage{tikz}
\usepackage{tikz-cd}

\newtheorem{thm}{Theorem}[section]
\newtheorem{lm}[thm]{Lemma}
\newtheorem{coro}[thm]{Corollary}
\newtheorem{prop}[thm]{Proposition}

\theoremstyle{definition}
\newtheorem{df}[thm]{Definition}
\newtheorem{ex}[thm]{Example}
\theoremstyle{plain}
\newtheorem{rmk}[thm]{Remark}

\newtheorem{claim}[thm]{Claim}

\newcommand{\R}{\mathbb{R}}
\newcommand{\C}{\mathbb{C}}
\newcommand{\Z}{\mathbb{Z}}
\newcommand{\A}{\mathbb{A}}

\newcommand{\N}{\mathbb{N}}
\newcommand{\Q}{\mathbb{Q}}
\newcommand{\GW}{\operatorname{GW}}
\newcommand{\W}{\operatorname{W}}
\newcommand{\Tr}{\operatorname{Tr}}

\newcommand{\K}{\mathbb{K}}

\newcommand{\lra}{\longrightarrow}

\newcommand{\val}{\operatorname{val}}

\newcommand{\mult}{\operatorname{mult}^{\A^1}}

\newcommand{\Card}{\operatorname{Card}}

\newcommand{\qinv}[1]{\left\langle#1\right\rangle}

\newcommand{\Tra}[2][E/k]{\operatorname{Tr}_{#1}\left(#2\right)}

\newcommand{\Puiseux}{\{\!\{t\}\!\}}

\newcommand{\Wel}{\operatorname{Wel}}

\newcommand{\Tor}{\operatorname{Tor}}
\newcommand{\etale}{M}

\makeatletter
\let\save@mathaccent\mathaccent
\newcommand*\if@single[3]{%
    \setbox0\hbox{${\mathaccent"0362{#1}}^H$}%
    \setbox2\hbox{${\mathaccent"0362{\kern0pt#1}}^H$}%
    \ifdim\ht0=\ht2 #3\else #2\fi
}
\newcommand*\rel@kern[1]{\kern#1\dimexpr\macc@kerna}
\newcommand*\widebar[1]{\@ifnextchar^{{\wide@bar{#1}{0}}}{\wide@bar{#1}{1}}}
\newcommand*\wide@bar[2]{\if@single{#1}{\wide@bar@{#1}{#2}{1}}{\wide@bar@{#1}{#2}{2}}}
\newcommand*\wide@bar@[3]{%
    \begingroup
    \def\mathaccent##1##2{%
        \let\mathaccent\save@mathaccent
        \if#32 \let\macc@nucleus\first@char \fi
        \setbox\z@\hbox{$\macc@style{\macc@nucleus}_{}$}%
        \setbox\tw@\hbox{$\macc@style{\macc@nucleus}{}_{}$}%
        \dimen@\wd\tw@
        \advance\dimen@-\wd\z@
        \divide\dimen@ 3
        \@tempdima\wd\tw@
        \advance\@tempdima-\scriptspace
        \divide\@tempdima 10
        \advance\dimen@-\@tempdima
        \ifdim\dimen@>\z@ \dimen@0pt\fi
        \rel@kern{0.6}\kern-\dimen@
        \if#31
        \overline{\rel@kern{-0.6}\kern\dimen@\macc@nucleus\rel@kern{0.4}\kern\dimen@}%
        \advance\dimen@0.4\dimexpr\macc@kerna
        \let\final@kern#2%
        \ifdim\dimen@<\z@ \let\final@kern1\fi
        \if\final@kern1 \kern-\dimen@\fi
        \else
        \overline{\rel@kern{-0.6}\kern\dimen@#1}%
        \fi
    }%
    \macc@depth\@ne
    \let\math@bgroup\@empty \let\math@egroup\macc@set@skewchar
    \mathsurround\z@ \frozen@everymath{\mathgroup\macc@group\relax}%
    \macc@set@skewchar\relax
    \let\mathaccentV\macc@nested@a
    \if#31
    \macc@nested@a\relax111{#1}%
    \else
    \def\gobble@till@marker##1\endmarker{}%
    \futurelet\first@char\gobble@till@marker#1\endmarker
    \ifcat\noexpand\first@char A\else
    \def\first@char{}%
    \fi
    \macc@nested@a\relax111{\first@char}%
    \fi
    \endgroup
}
\makeatother

\title{A quadratically enriched correspondence theorem}
\author{Jaramillo Puentes, Andr\'es \and Pauli, Sabrina}

\begin{document}
\maketitle

\begin{abstract}
We quadratically enrich Mikhalkin's correspondence theorem. That is, we prove a correspondence between algebraic curves on a toric surface counted with Levine's quadratic enrichment of the Welschinger sign, and tropical curves counted with a quadratic enrichment of Mikhalkin's multiplicity for tropical curves.
\end{abstract}
{\hypersetup{linkcolor=black}\tableofcontents}
\setcounter{tocdepth}{2}
\section{Introduction}
Let $k$ be a field, $\Delta\subset\Z^2$ be a convex lattice polygon, $\Tor_k(\Delta)$ be the associated toric surface and $\Lambda_k(\Delta)$ be the linear system on $\Tor_k(\Delta)$ generated by monomials $x^iy^j$ with~$(i,j)\in \Delta\cap\Z^2$. 
Let
\[N(\Delta,g,\mathcal{P})\]
be the count
of genus $g$ curves in $\Tor_k(\Delta)$ passing through a configuration of  
\[n\coloneqq \operatorname{dim}_k(\Lambda_k(\Delta))+g-\delta=\vert \Delta\cap\Z^2\vert-1-\delta+g=\vert \partial\Delta\cap\Z^2\vert-1+g\]
points $\mathcal{P}=\{p_1,\ldots,p_n\}$ in general position in $\Tor_k(\Delta)$. Here, $\delta$ is the number of nodes of a generic rational 
curve in $\Tor_k(\Delta)$. 
This number is independent of the choice of configuration of points when we consider all curves defined over an algebraic closure $\overline{k}$ of~$k$,
and we can drop the $\mathcal{P}$ in the notation
\[N(\Delta, g)=N(\Delta, g,\mathcal{P}).\]
The invariant $N(\Delta,g)$ counts all curves defined over an algebraic closure.
When one is interested in the naive count of curves defined over a field $k$ which is not algebraically closed, one loses invariance. That means, in this case the count depends on the chosen configuration $\mathcal{P}$ of the $n$ points. For example when $k=\R$, there can be $8$, $10$ or $12$ real rational plane cubic curves through a generic configuration of $8$ real points in $\mathbb{P}^2_\R$, depending on the chosen configuration (cf.~\cite[Proposition 4.7.3]{Degtyarev_2000}).
In~\cite{Welschinger}, Jean-Yves Welschinger introduced a way to restore invariance by defining a signed count of real curves. A real curve $C$ has two \emph{types} of real nodes:
\begin{itemize}
    \item the \emph{real split node}, locally defined by the equation $x^2-y^2=(x-y)(x+y)$;
    \item and the \emph{real solitary node}, locally defined by the equation $x^2+y^2$.
\end{itemize}
If $f$ is a defining equation for $C$, then a real node $z$ is split (respectively, solitary) if and only if 
$\det \operatorname{Hessian} f(z)$ is negative (respectively, positive).
The \emph{Welschinger sign} of a real curve $C\subset \Tor_\R(\Delta)$ is the sign
\begin{equation}
\label{eq:realWelschinger}
\operatorname{Wel}_{\mathbb{R}}(C)\coloneqq (-1)^{\#\text{solitary nodes}}=\prod_{\text{$z$ a real node}}\operatorname{sign}\left(-\det \operatorname{Hessian} f(z)\right).
\end{equation}
Given a configuration of~$n$ real points $\mathcal{P}=\{p_1,\ldots,p_n\}$ in $\Tor_\R(\Delta)$, we define the sum 
\[W(\Delta,g,\mathcal{P})\coloneqq\sum \Wel_\R(C)\]
over all real curves $C$ in $\Tor_\R(\Delta)$ passing through the configuration $\mathcal{P}$.
Welschinger proved that if $\Tor_\R(\Delta)$ is a toric del Pezzo surface and $g=0$, then the sum~$W(\Delta,0,\mathcal{P})$ is independent of the choice of configuration of points as long as it is generic. Hence, in this case we write
\[W(\Delta)\coloneqq W(\Delta,0,\mathcal{P}).\]
In order to generalize $\Wel_\R(C)$ to curves with ground field $k$, we look at the different \emph{types} of nodes of an algebraic curve with residue field $k$, as in the case $k=\R$. We assume $\operatorname{char}k\neq 2$. There are:
\begin{itemize}
    \item split nodes, locally defined by $x^2-ay^2$ for $a\in(k^\times)^2$, that is, the unit $a$ is a square in~$k$;
    \item and solitary nodes of \emph{type} $a$, where $a$ is not a square in $k^\times$, locally defined by the equation $x^2-ay^2$.
\end{itemize}

Note that there is only one type of split node: Since $a$ has a square root, we can assume that the local equation for the split node is $x^2-y^2$ after a coordinate change. However, there might be different types of solitary nodes corresponding to the different non-square classes of $k^\times/(k^\times)^2$. 
We define the \emph{type} of a node $z$ defined over $k$ to be 
$\operatorname{Type}(z)\coloneqq a\in k^\times/(k^\times)^2.$
If the curve $C$ is defined by a polynomial~$f$, then one can compute the type of a $k$-rational node $z$ by computing the image of~$-\det \operatorname{Hessian}f(z)$
in~$k^\times/(k^\times)^2$. That is,
\[\operatorname{Type}(z)= -\det \operatorname{Hessian}f(z)\in k^\times/(k^\times)^2.\]
The type records the quadratic field extension of $k$ where the branches of the node live, i.e., they are Galois conjugate and defined over $k(\sqrt{a})$.

The different types of nodes correspond to the classes in $k^\times/(k^{\times})^2$. These classes can be seen as the generators of the Grothendieck-Witt ring $\GW(k)$ of non-degenerate symmetric bilinear forms over the field $k$ (see e.g. \cite{PauliWickelgren} or \cite[$\S 2$]{JPP}).
Methods from $\A^1$-homotopy theory allow to answer questions in enumerative geometry over an arbitrary field $k$. These answers are invariants that belong to 
 $\GW(k)$. 
 The ring $\GW(k)$ has a presentation, whose set of generators is in bijection with the classes $k^\times/(k^{\times})^2$ and is given by the rank one forms $\langle a\rangle$ where $a\in k^\times$ (see Definition~\ref{df:GW}). We use the notation $h$ for the hyperbolic form $h=\langle1\rangle+\langle-1\rangle$. 
 Additionally, for a finite separable field extension~$L/k$, there is a \emph{trace map} $\Tr_{L/k}\colon \GW(L)\rightarrow\GW(k)$. 
This leads to the following generalization of the Welschinger sign defined by Marc Levine in \cite{LevineWelschinger}.
\begin{df}
The \emph{(Levine-Welschinger) quadratic weight} of a curve $C$ with residue field $\kappa(C)$ is 
    \[\Wel^{\A^1}_{\kappa(C)}(C)\coloneqq\qinv{\prod_{\text{nodes } z}N_{\kappa(z)/\kappa(C)}(-\det \operatorname{Hessian}f(z))}\in \GW(\kappa(C)),\]
where $\kappa(z)$ is the residue field of $z$ (which we assume to be separable over $\kappa(C)$) and the map
  $N_{\kappa(z)/\kappa(C)}\colon \kappa(z)\rightarrow \kappa(C)$ is the field norm.
\end{df}
Observe that for $k=\R$ we have that 
\begin{equation}
\label{eq:WelR}
\Wel_\R^{\A^1}(C)=\langle\Wel_\R(C)\rangle\in \GW(\R).\end{equation}

Let
\begin{equation}
\label{eq:defnNA1}
N^{\A^1}_k(\Delta,g, \mathcal{P})\coloneqq \sum \Tr_{\kappa(C)/k}\left(\Wel^{\A^1}_{\kappa(C)}(C)\right),\end{equation}
where the sum runs over all genus $g$ curves $C$ in $\Tor_k(\Delta)$ passing through a configuration $\mathcal{P}=\{p_1,\ldots,p_n\}$ of $n$ $k$-rational points in general position. Note that this sum may include curves $C$ that are not defined over the base field $k$, but they are defined over a finite field extension $\kappa(C)$ of $k$. Taking the trace $\Tr_{\kappa(C)/k}$ of its quadratic weight yields an element of $\GW(k)$. 
\begin{rmk}
In case $k=\R$, the only possibilities for the residue field $\kappa(C)$ in \eqref{eq:defnNA1} are~$\R$ or~$\C$. In case $\kappa(C)=\R$, we have that $\Tr_{\kappa(C)/\R}=\operatorname{id}$. 
When $\kappa(C)=\C$, the quadratic weight is $\Wel^{\A^1}_{\kappa(C)}(C)=\langle1\rangle$, and $\Tr_{\kappa(C)/\R}(\Wel_{\kappa(C)}^{\A^1}(C))=\langle1\rangle+\langle-1\rangle=h$ which has signature~$0$. It follows from \eqref{eq:WelR} that the signature of $N^{\A^1}_\R(\Delta,\mathcal{P})$ equals
\[\operatorname{sgn}(N^{\A^1}_\R(\Delta,g,\mathcal{P}))=\sum\operatorname{sgn}(\Wel_\R^{\A^1}(C)) =W(\Delta,g,\mathcal{P})\]
where the sum runs over the curves with residue field $\R$.
\end{rmk}

Suppose $k$ is not of characteristic $2$ or $3.$ A definition of invariant quadratically enriched counts of rational curves through $n$ $k$-points in general position on a del Pezzo surface was given in \cite[Example 3.9]{LevineWelschinger} for $k$ infinite. Work of Jesse Kass, Marc Levine, Jake Solomon and Kirsten Wickelgren \cite{KLSWCount,KLSWorientation} announced in \cite{PauliWickelgren} gives an invariant quadratically enriched count of rational curves passing through given Galois orbits of not-necessarily $k$-rational points on a del Pezzo surface of degree at least $4$ when the field $k$ is perfect.  The invariant is defined as the $\A^1$-degree of an evaluation map which is shown to be relatively oriented. The $\A^1$-degree of the evaluation map is shown to coincide with the count of rational curves $C$ passing through given points in general position weighted by $\Wel_{k(C)}^{\A^1}(C).$ In particular, $N_k^{\A^1}(\Delta,0,\mathcal{P})$ is independent of the choice of points $\mathcal{P}$ when $\Tor_k(\Delta)$ is a del Pezzo surface.
In this case, we denote this invariant by 
\[N^{\A^1}_k(\Delta)\coloneqq N^{\A^1}_k(\Delta,0, \mathcal{P}).\]

We aim to find a way to compute $N^{\A^1}_k(\Delta)$ for $\Delta$ defining a toric del Pezzo surface~$\Tor_k(\Delta)$.
Tropical geometry provides a way to compute $N(\Delta,g)$ and $W(\Delta)$. More precisely, Grigory Mikhalkin's celebrated correspondence theorem \cite{Mikhalkin} establishes a correspondence between algebraic curves and tropical curves counted with the following multiplicities: Let $A$ be a simple tropical curve with Newton polygon $\Delta$ and dual subdivision $S$. For a~$3$-valent vertex $v$ of $A$, let $\Delta_v$ be the dual triangle in $S$.
Then complex \emph{vertex multiplicity} of $v$ is defined as 
\[ m_v^\C\coloneqq \vert \Delta_v\vert,\]
where $\vert \Delta_v\vert$ denotes the double Euclidean area of $\Delta_v$.  The real \emph{vertex multiplicity} of $v$ is defined as 
\[m_v^\R\coloneqq \begin{cases} (-1)^{\operatorname{Int}(\Delta_v)} &\text{if all edges in $\Delta_v$ have odd lattice length},\\0 &\text{else}.\end{cases}\]
Here, $\operatorname{Int}(\Delta_v)$ is the number of interior lattice points of $\Delta_v$.
The complex and the real multiplicities of the tropical curve $A$ are the products of the complex and the real vertex multiplicities of $A$ over all $3$-valent vertices, respectively,
\[\operatorname{mult}_\C(A)\coloneqq \prod_{v\text{ a $3$-valent vertex of $A$}}\mkern-45mu m_v^\C\mkern45mu \text{ and } \operatorname{mult}_\R(A)\coloneqq \prod_{v\text{ a $3$-valent vertex of $A$}}\mkern-45mu m_v^\R.\]
Set
\[N^{\operatorname{trop}}(\Delta,g)\coloneqq \sum \operatorname{mult}_\C(A) \text{ and }W^{\operatorname{trop}}(\Delta,g)\coloneqq \sum \operatorname{mult}_\R(A),\]
where both sums run over all genus $g$ tropical curves $A$ through a $\Delta$-generic configuration of $n$ points in $\R^2$.
Milkhalkin's correspondence theorems states that 
\[N(\Delta,g)=N^{\operatorname{trop}}(\Delta,g)\text{ and }W(\Delta)=W^{\operatorname{trop}}(\Delta,0).\]
This allows to translate the question of counting algebraic curves to the question of counting tropical curves, which can be done with combinatorial methods.
By now there are several proofs of Mikhalkin's correspondence theorem, e.g. \cite{Mikhalkin,Shustin,NishinouSiebert,Tyomkin,ArguzBousseau,MandelRuddat}. 
The goal of this paper is to  prove a correspondence theorem between algebraic curves over the field of \emph{Puiseux series} counted with their quadratic weight and tropical curves counted with the correct quadratic multiplicity.
Recall that the field of Puiseux series 
\[\K\coloneqq k\Puiseux=\{a_0t^{q_0}+a_1t^{q_1}+\ldots\mid  q_i\in \Q, q_0<q_1<\ldots,\, q_i=\frac{q_i'}{N},\, {q_i'},N\in\Z\}\]
has a non-Archimedean valuation
\[\operatorname{val}\colon \K\rightarrow \R\cup\{-\infty\}\colon \,a_0t^{q_0}+a_1t^{q_1}+\ldots\mapsto -q_0.\]
Let $p_1,\ldots,p_n\in (\K^\times)^2$ be in general position and let $x_i\coloneqq \operatorname{val}(p_i)\in \R^2$ for $i=1,\ldots,n$. We find a correspondence between genus $g$ algebraic curves in $\Tor_\K(\Delta)$ through $p_1,\ldots,p_n$ and genus $g$ tropical curves with Newton polygon $\Delta$ through $x_1,\ldots,x_n$, where we count the algebraic curves with their quadratic weight $\Wel^{\A^1}_{\kappa(C)}(C)$ and the tropical curves with their \emph{quadratically enriched multiplicity} $\mult_k$
defined as follows.
\begin{df}
\label{df: enriched multiplicity}
Let $A$ be a tropical curve.
We define the \emph{quadratically enriched vertex multiplicity}, or \emph{quadratic multiplicity} for short, of a $3$-valent vertex $v$ of $A$ to be  
\[m_v^{\A^1}\coloneqq\begin{cases}
    \qinv{ (-1)^{\operatorname{Int}(\Delta_v)}\vert\sigma_v\vert\vert\sigma'_v\vert\vert\sigma''_v\vert}+\frac{\vert\Delta_v\vert-1}{2}h&\text{if $\Delta_v$ has only odd edges,}\\
    \frac{\vert\Delta_v\vert}{2}h&\text{ else.}
\end{cases}\]
Here, $\sigma_v$, $\sigma_v'$ and $\sigma_v''$ are the three edges of the triangle $\Delta_v$ dual to $v$ and $\operatorname{Int}(\Delta_v)$ is the number of interior lattice points of $\Delta_v$.
The \emph{quadratically enriched multiplicity}, or \emph{quadratic multiplicity} for short, of the tropical curve $A$ is the product of its vertex multiplicities over all its~$3$-valent vertices
    \[\mult_k(A)\coloneqq\prod_{v\text{ a 3-valent vertex of $A$ }}\mkern-45mu m_v^{\A^1}.\]
\end{df}
\begin{rmk}
\label{rmk:rankandsign}
    The rank of the quadratic form $\mult_k(A)$ equals the complex multiplicity~$\operatorname{mult}_\C(A)$. If $k\subset\R$, then  the signature of $\mult_k(A)$ is the real multiplicity $\operatorname{mult}_\R(A)$.
\end{rmk}

\begin{ex}
    \begin{figure}
    \begin{tabular}{ccc}
    \includegraphics[scale=0.1]{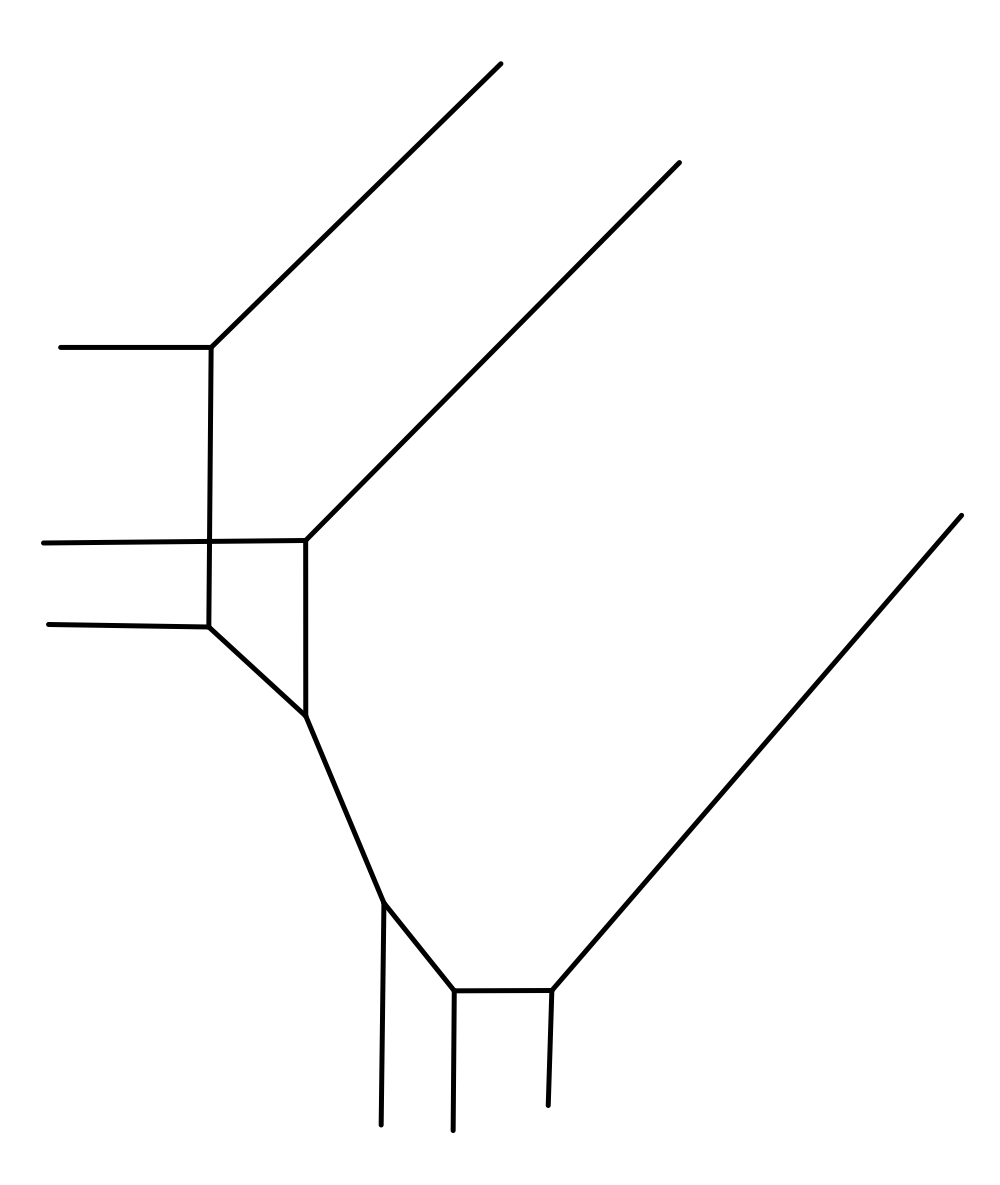}&
    \includegraphics[scale=0.1]{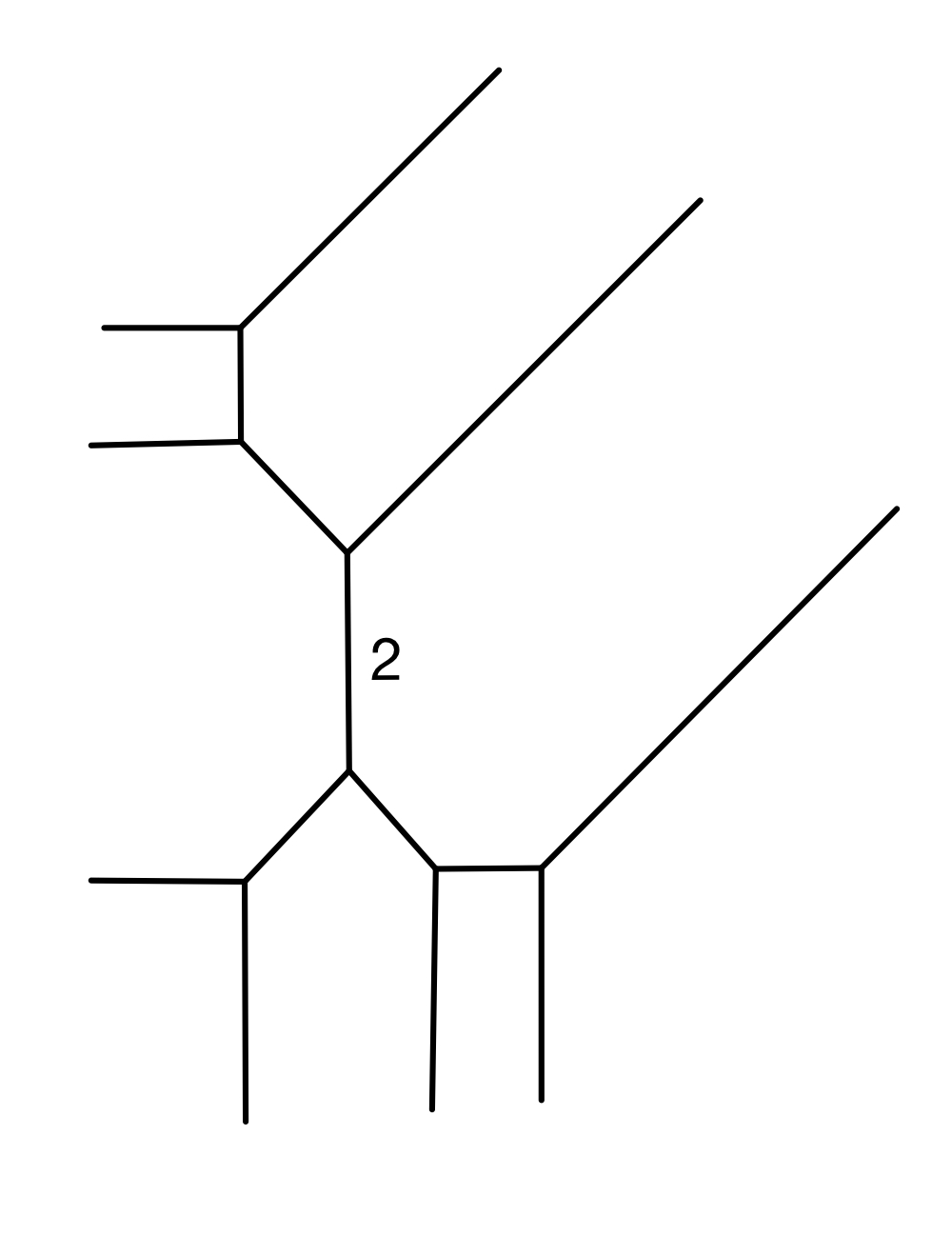}&
    \includegraphics[scale=0.1]{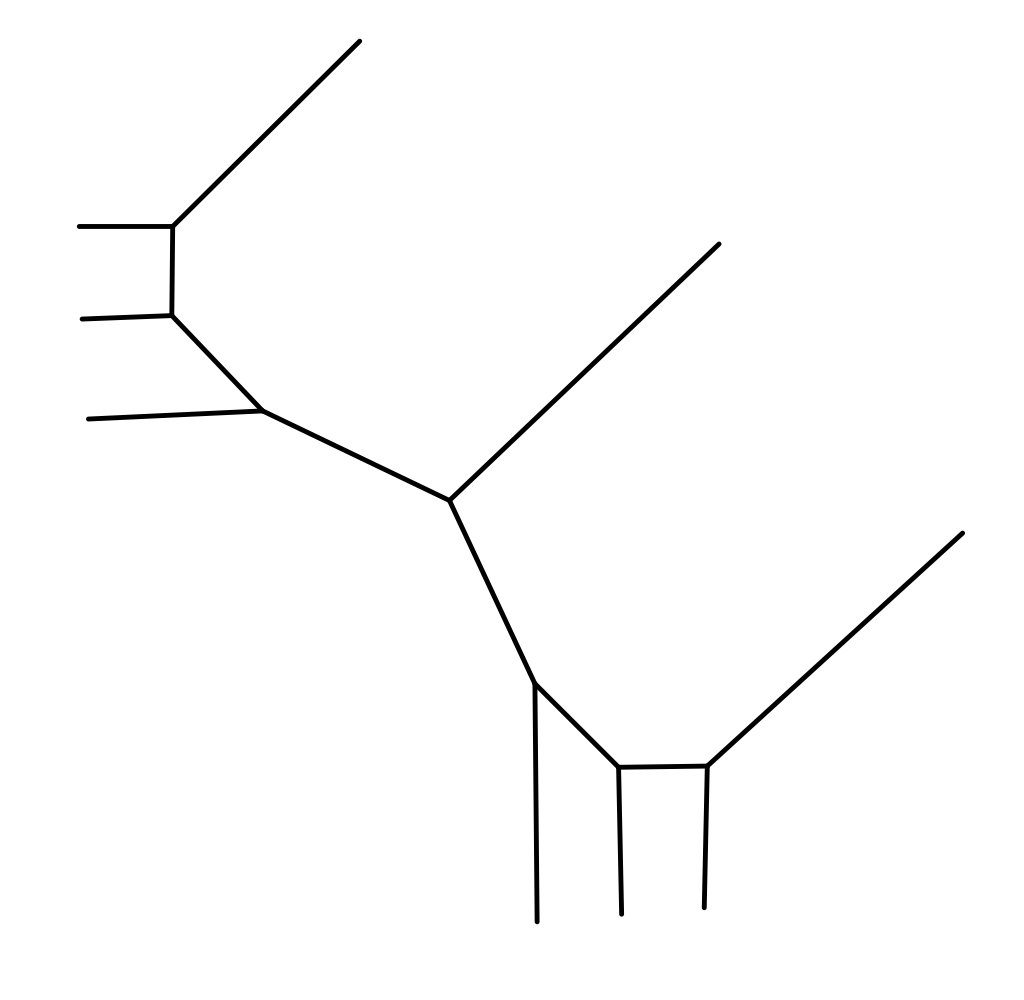}\\
    \includegraphics[scale=0.1]{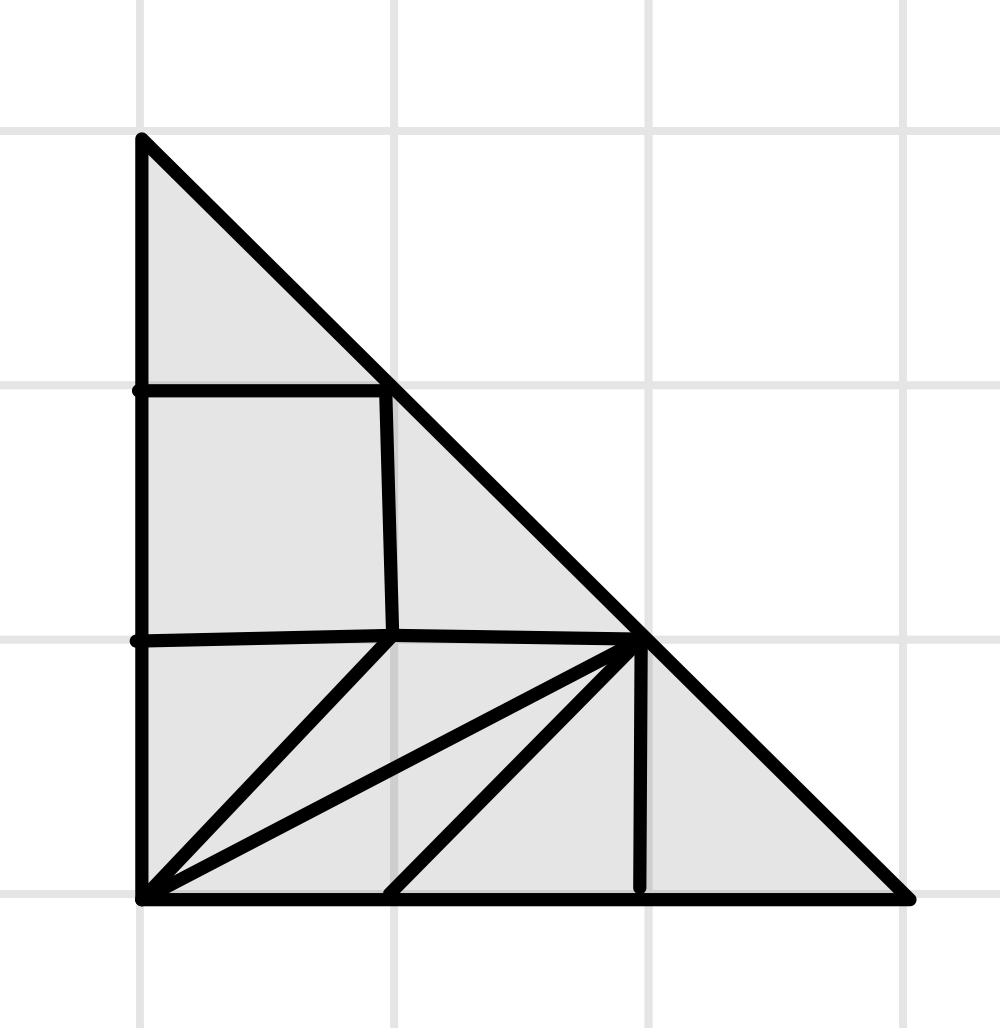}&
    \includegraphics[scale=0.1]{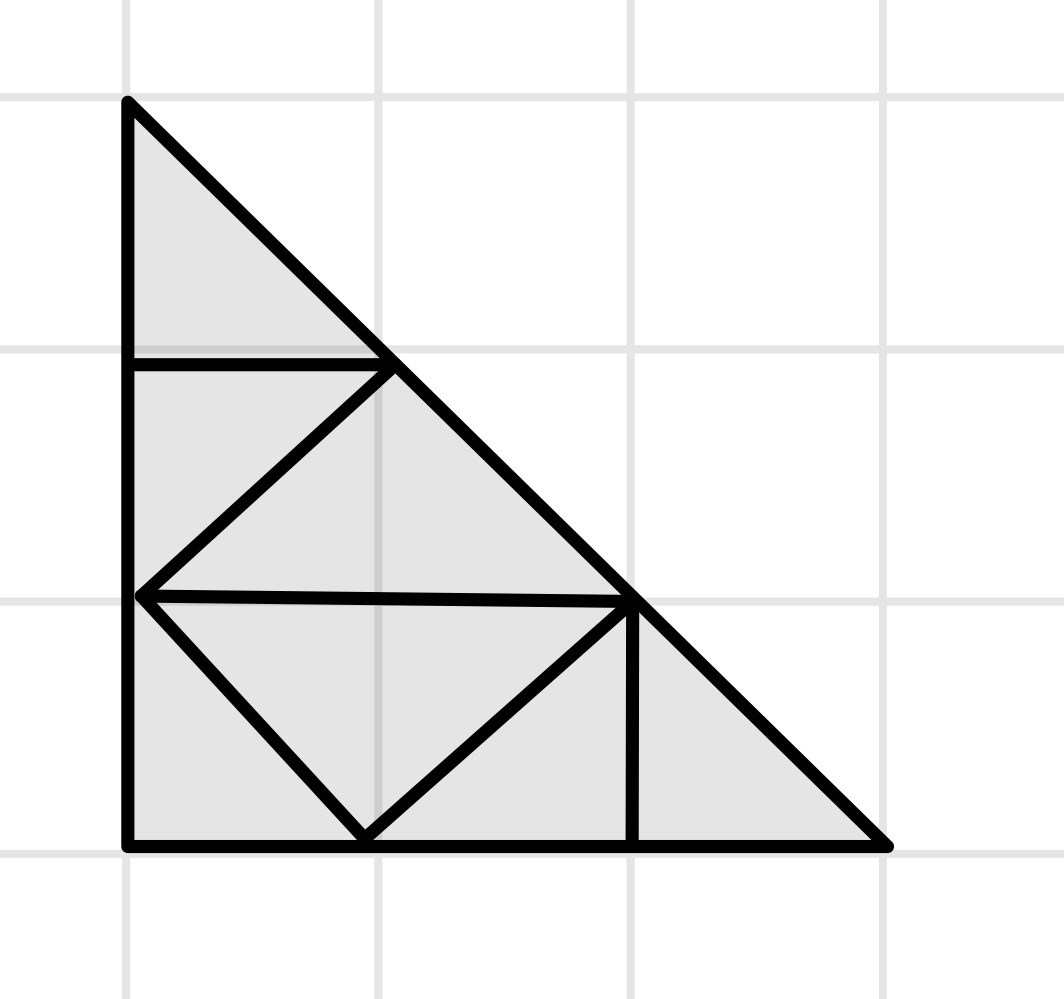}&
    \includegraphics[scale=0.1]{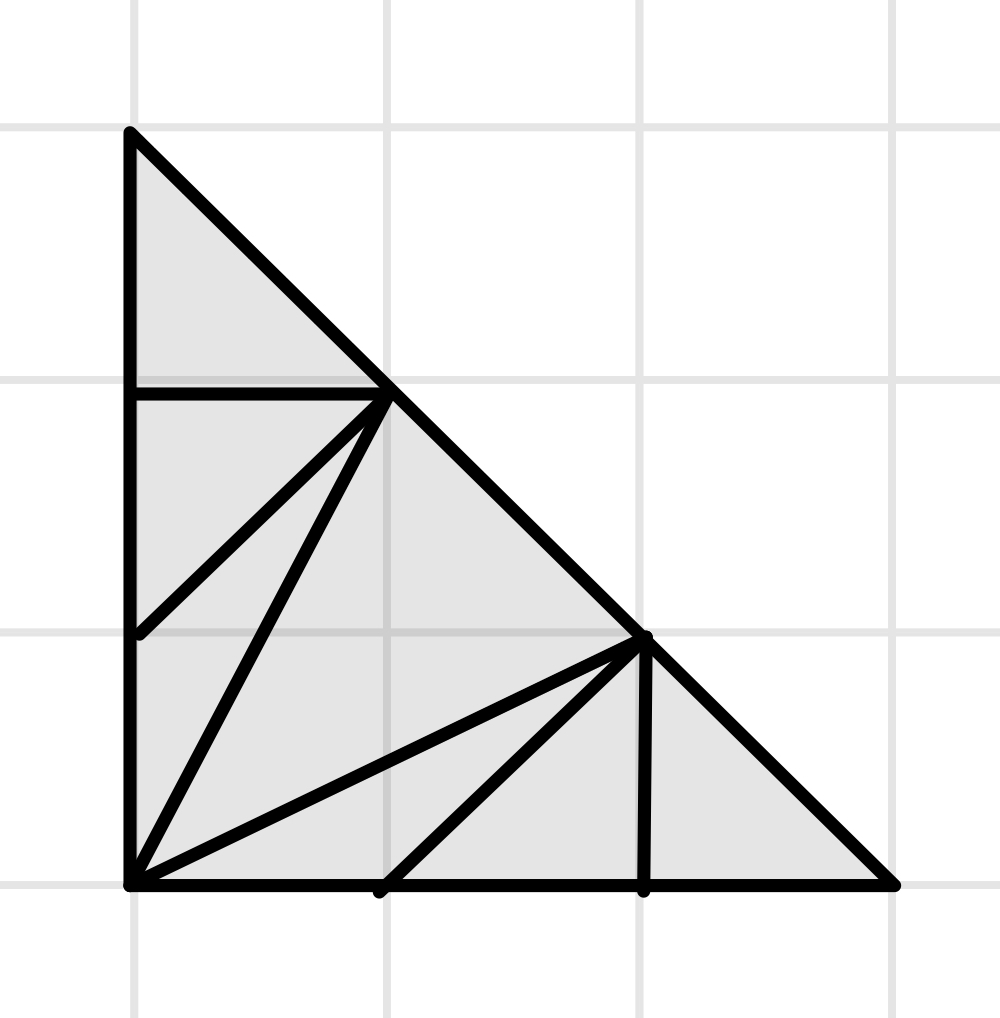}
    \end{tabular}
    \caption{Examples of tropical curves with their dual subdivision to which we compute their multiplicities.\label{fig:examplesmultiplicities}}
    \end{figure}
Figure \ref{fig:examplesmultiplicities} shows three tropical curves $A_1$, $A_2$ and $A_3$ and their dual subdivisions $S_1$, $S_2$ and $S_3$ below. 
We compute their complex, real and quadratic multiplicities.

\begin{table}[h!]
\centering
\begin{tabular}{|c|c|c|}
\hline
    $\operatorname{mult}_\C(A_1)=1$ & $\operatorname{mult}_\R(A_1)=1$ & $\operatorname{mult}^{\A^1}_k(A_1)=\langle1\rangle$\\
    \hline
    $\operatorname{mult}_\C(A_2)=4$ & $\operatorname{mult}_\R(A_2)=0$ & $\operatorname{mult}^{\A^1}_k(A_2)=2h$\\
    \hline
    $\operatorname{mult}_\C(A_3)=3$ & $\operatorname{mult}_\R(A_3)=-1$ & $\operatorname{mult}^{\A^1}_k(A_3)=h+\langle-1\rangle$\\
    \hline
\end{tabular}
\end{table}

The first tropical curve $A_1$ has a $4$-valent vertex and thus a parallelogram in the dual subdivision. Note that when computing the multiplicity of the curve, the parallelograms do not play a role. 
All edges of $S_1$ have length $1$ and the double Euclidean areas of the triangles are all equal to $1$. Furthermore, there are no interior lattice points in the triangles of $S_1$. It follows that all complex and real vertex multiplicities equal $1$ and the quadratic multiplicity of $A_1$ equals $\qinv{1}$.
The tropical curve $A_2$ in the middle has an edge of weight~$2$ corresponding to an edge of even length in $S_2$. Thus its real multiplicity is~$0$ and its quadratic multiplicity is a multiple of the hyperbolic form.
The last tropical curve $A_3$ has an interior lattice point in the triangle in the middle of its dual subdivision $S_3$. It follows that its real multiplicity is $-1$ and we get a $\langle-1\rangle$ summand in the quadratic multiplicity.
\end{ex}

Starting with an algebraic curve $C\in \Tor_\K(\Delta)$ through $p_1,\ldots, p_n$ defined by $f$, one can associate a tropical curve 
\begin{equation}
\label{eq:Kapranov}
A_f\coloneqq \overline{\{(\val(x),\val(y)): \text{ $(x,y)$ is a geometric point in $Z_f$\}}}\subset \R^2\end{equation} 
through $x_1,\ldots,x_n$ with Newton polygon $\Delta$ where $Z_f=\{f=0\}\subset \operatorname{Spec}(\K[x^{\pm1},y^{\pm1}])$. We call $A_f$ the \emph{tropical curve associated to $f$}.

\begin{figure}
    \includegraphics[width=\textwidth]{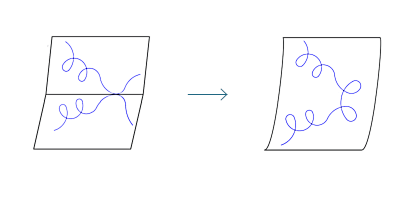}
    \caption{Deformation into nodes. \label{fig:defintonodes}}
    \end{figure}

This yields one direction of the correspondence theorem; the other direction is much more subtle.
The idea is to find all genus $g$ curves in $\Tor_\K(\Delta)$ through $p_1,\ldots,p_n$ that tropicalize to a given tropical curve $A$ through $x_1,\ldots,x_n$ and compute their quadratic weights. To do this, we need to figure out what information is needed to find these curves.
We have seen that a polynomial $f$ defining a curve in $\Tor_\K(\Delta)$ defines a tropical curve~$A_f$ with Newton polygon $\Delta$. Let $S_f\colon \Delta=\Delta_1\cup\ldots\cup\Delta_N$ be its dual subdivision. By the results in \cite[$\S3$]{NishinouSiebert}, \cite[$\S2.3$]{Shustin} or \cite[Proposition 3.5]{Tyomkin}, this defines a flat family of toric varieties
\[\pi\colon \mathcal{X}\rightarrow \A^1_k\]
such that $X_t= \Tor_{\kappa(t)}(\Delta)$ for $t\neq 0$ and 
\[X_0= \bigcup_{i=1}^N\Tor_k(\Delta_i),\]
such that when $\Delta_i$ and $\Delta_j$ have a common edge $\sigma=\Delta_i\cap\Delta_j$, then $\Tor_k(\Delta_i)$ and $\Tor_k(\Delta_j)$ are glued along $\Tor_k(\sigma)$ in $X_0$.
Furthermore, one can show that $f$ gives rise to a family of curves in $\mathcal{X}$  with special curve
\begin{equation}
\label{eq:familyofcurves}
C_0=\bigcup_{i=1}^NC_i\subset \bigcup_{i=1}^N\Tor_k(\Delta_i)=X_0,
\end{equation}
cf. \cite[Propsition 6.3]{NishinouSiebert}.
It holds that if $\Delta_i$ and $\Delta_j$ have a common edge $\sigma=\Delta_i\cap\Delta_j$, then $C_i$ and $C_j$ meet in exactly one point $z\in \Tor_k(\sigma)$ and both curves meet $\Tor_k(\sigma)$ with the same multiplicity $m=\vert\sigma\vert$. 
Furthermore, locally $z$ deforms into $m-1$ nodes in the family as illustrated in Figure \ref{fig:defintonodes}.

After a non-toric blow up of $X_0$ one gets a local equation $\mathcal{T}_{[\sigma]}f$ for the curve with these~$m-1$ nodes. Given a genus $g$ simple tropical curve $A$ with Newton polygon $\Delta$ through points $x_1,\ldots, x_n$, we are looking for a polynomial equation
    \[f(x,y)=\sum_{(i,j)\in \Delta\cap\Z^2}\widetilde{c}_{ij}x^iy^j\]
with coefficients $\widetilde{c}_{ij}$ defined over a finite (and separable) field extension of $\K$ defining a genus $g$ algebraic curve in $\Tor_\K({\Delta})$ such that 
\[f(p_1)=\ldots=f(p_n)=0,\]
and such that the associated tropical curve $A_f$ equals $A$.
Eugenii Shustin's way of restoring all such $f(x,y)$ heavily relies on the patchworking construction first used by Oleg Viro.
Roughly speaking, the input for Shustin's patchworking theorem as stated in \cite[Theorem~2.51]{tropicalBook} is 
\begin{itemize}
    \item the tropical curve $A$, its dual subdivision $S\colon\Delta=\Delta_1\cup\ldots\cup \Delta_N$ and a convex piecewise linear function $\nu:\Delta\rightarrow \R$ whose graph projects to the dual subdivision~$S$;
    \item polynomials 
    \[f_l(x,y)=\sum_{(i,j)\in \Delta_l\cap \Z^2}a_{ij}x^iy^j\]
    for $l=1,\ldots,N$, such that whenever $\Delta_l$ and $\Delta_s$ have a common edge $\sigma$ the coefficients $a_{ij}$ of $f_l$ and $f_s$ agree for each $(i,j)\in \sigma$, defining rational curves in~$\Tor_k(\Delta_l)$ for $l=1,\ldots,N$, with the same properties as the curves $C_l$ in \eqref{eq:familyofcurves};
    \item \emph{deformation patterns} $f_{[\sigma]}$ for each edge $\sigma$ of the dual subdivision $S$ of $A$ of lattice length $\vert\sigma\vert\ge 2$ with the same properties as $\mathcal{T}_{[\sigma]}f$;
    \item \emph{initial conditions} to pass through the points $p_1,\ldots,p_n$.
\end{itemize}
The output is a polynomial of the form
\begin{equation}
    f(x,y)=\sum_{(i,j)\in \Delta\cap\Z^2}(a_{ij}+c_{ij}(t))t^{\nu(i,j)}x^iy^j
\end{equation}
defining a unique genus $g$ curve $C$ in $\Tor_{\K}(\Delta)$ with $c_{ij}(t)\in \K$ having only positive exponents, i.e. $c_{ij}(0)=0$, such that $f$ matches the input data.
Shustin's proof methods are analytic and work for $k=\C$ or $k=\R$. 
Takeo Nishinou and Bernd Siebert use \emph{log deformation theory} to restore the unique algebraic curve in their correspondence theorem~\cite{NishinouSiebert}. Nishinou-Siebert’s approach must work in large positive characteristic though this is not claimed. Ilya Tyomkin generalizes the Nishinou-Siebert correspondence theorem and also states it for large positive characteristic \cite{Tyomkin}. Nishinou-Siebert's and Tyomkin's approaches are purely algebro-geometric and hence we can apply it to our set up. We use the same characteristic assumption as Tyomkin in \cite{Tyomkin}, namely that $\operatorname{char}k=0$ or~$\operatorname{char}k>\operatorname{diam}(\Delta)$. This assures that the field extensions that occur in our computations are all separable over the base field and that the quadratic weights are well defined.

Then, we show that the quadratic weight of this unique curve can be computed as the product of the quadratic weights of the curves defined by the polynomial equations $f_l$ for~$l=1,\ldots,N$, and the polynomials equations $f_{[\sigma]}$ from the input data, generalizing the real case in \cite{Shustin}. 
Shustin shows that there are finitely many choices for the input of the patchworking theorem, namely $\operatorname{mult}_\C(A)$ many choices (when counting curves defined over the algebraic closure  of $k$). 
Hence, for  these choices we need to compute the sum
\[\sum \Tr_{\kappa(C)/\K}\left(\Wel^{\A^1}_{\kappa(C)}(C)\right)\in \GW(\K)\cong \GW(k)\]
over all genus $g$ algebraic curves $C$ in $\Tor_\K(\Delta)$ through $p_1,\ldots,p_n$ with associated tropical curve $A$. That is, the sum over all curves defined by the polynomials $f(x,y)$ we constructed using the refined patchworking theorem arising from the $\operatorname{mult}_\C(A)$ choices for the input data.
We show that this sum equals the quadratically enriched multiplicity $\mult_k(A)$ defined in Definition \ref{df: enriched multiplicity}.

We proceed along the lines of \cite{Shustin}. Let $A$ be a genus $g$ tropical curve with Newton polygon $\Delta$ passing through a $\Delta$-generic configuration $x_1,\ldots,x_n$. We need to find all possibilities for the input data of the Patchworking theorem and compute the quadratic weights of the $C_l$'s and $C_{[\sigma]}'s$ for the different choices of the input data.

Shustin does this for the real case in Lemma 3.5 ($C_l$ for~$\Delta_l$ a triangle), Lemma 3.6 ($C_l$ for~$\Delta_l$ a parallelogram), and Lemma 3.9 (deformation patters $C_{[\sigma]}$) in \cite{Shustin}.
We generalize these Lemmas to our setting in Section \ref{section:quadraticenrichments}.
From this, we get the quadratically enriched multiplicity and the following theorem.

\begin{thm}[Quadratically enriched correspondence theorem]
\label{thm:QEcorrespondence}
    Assume $k$ is a field of characteristic $0$ or $\operatorname{char}k>\operatorname{diam}(\Delta)$.
    Let $\mathcal{P}=\{p_1,\ldots,p_n\}\subset \Tor_\K(\Delta)$ be a generic configuration of points such that $x_1,\ldots,x_n$ is $\Delta$-generic. Let $A$ be a genus $g$ tropical curve with Newton polygon $\Delta$ passing through $x_1,\ldots,x_n$.
    Then under the canonical isomorphism~$\GW(k)\cong \GW(\K)$ the quadratic multiplicity $\mult_k(A)$ of $A$ is mapped to
    \[\sum \Tr_{\kappa(C)/\K}\left(\Wel^{\A^1}_{\kappa(C)}(C)\right),\]
    where the sum runs over all genus $g$ curves $C$ in $\Tor_\K(\Delta)$ passing through $\mathcal{P}$ and tropicalizing to $A$. 
\end{thm}

Theorem \ref{thm:QEcorrespondence} does not require the base field to be perfect, since the fields of definition under consideration are separable field extensions. This is because they are extensions of degree less than the characteristic, in the case of positive characteristic.
However, for any field $k$ with $\operatorname{char} k\neq 2$, we can set 
\[N^{\A^1,\operatorname{trop}}_k(\Delta,g)\coloneqq\sum \mult_k(A),\]
where the sum runs over all genus $g$ tropical curves $A$ through the $\Delta$-generic configuration of points $x_1,\ldots,x_n$.

We have that if $A$ is a simple tropical curve, then its multiplicity $\mult_k(A)$ consists of sums of $\langle1\rangle$'s and $\langle-1\rangle$'s, and so it is completely determined by its rank $\operatorname{mult}_\C(A)$ and the signature $\operatorname{mult}_\R(A)$ (see Remark \ref{rmk:rankandsign}). 
A consequence of this fact, together with Mikhalkin's Correspondence theorem and the invariance of the real tropical counts in~\cite{IKSCapHar}, is that at the level of enriched tropical invariants, we have that
\[N^{\A^1,\operatorname{trop}}_k(\Delta,g)=\frac{N^{\operatorname{trop}}(\Delta,g)-W^{\operatorname{trop}}(\Delta,g)}{2}h+W^{\operatorname{trop}}(\Delta,g)\langle1\rangle\in \GW(k).\]

Assume now that the base field~$k$ is perfect of $\operatorname{char}{k}\neq 2,3$ and that~$\Tor_k(\Delta)$ is a del Pezzo surface.
This implies that~$\K$ is perfect of $\operatorname{char}{\K}\neq 2,3$ and that~$\Tor_\K(\Delta)$ is a del Pezzo surface.
Thus, by \cite{KLSWCount}, the count $N^{\A^1}_\K(\Delta,0,\mathcal{P})=N^{\A^1}_\K(\Delta)$ is independent of~$\mathcal{P}$.
In particular, the invariant $N^{\A^1}_\K(\Delta)$ equals $N^{\A^1}_\K(\Delta,0,\mathcal{Q})$ where $\mathcal{Q}$ is a point configuration that specializes to a generic configuration of $k$-points $\overline{\mathcal{Q}}$ in $\Tor_k(\Delta)$. The latter coincides with the count of rational curves $N^{\A^1}_k(\Delta)$ in $\Tor_k(\Delta)$ passing through $\overline{\mathcal{Q}}$, and therefore, we get that
\[N^{\A^1}_\K(\Delta)=N^{\A^1}_k(\Delta).\]

This implies that $\operatorname{rk} (N^{\A^1}_k(\Delta))=N(\Delta,0)$ and $\operatorname{sgn}(N^{\A^1}_k(\Delta))=W(\Delta)$ whenever $k\subset\R$. As a consequence we get the following corollary of Theorem~\ref{thm:QEcorrespondence}.

\begin{coro}
\label{coro:correspondencedelPezzo}
If $\Tor_k(\Delta)$ is a toric del Pezzo surface and $k$ is a perfect field of characteristic $0$ or $\operatorname{char}k>max\{\operatorname{diam}(\Delta),3\}$, then
   \[N^{\A^1}_k(\Delta)=\frac{N(\Delta,0)-W(\Delta)}{2}h+W(\Delta)\langle1\rangle\in \GW(k).\]
\end{coro}



As an example, the count of rational degree $d$ curves in $\mathbb{P}_k^2$ passing through a configuration of $3d-1$ general $k$-points with their quadratic weight in $\GW(k)$ equals $N^{\A^1}_k(\Delta_d)$, where~$\Delta_d=\operatorname{Conv}\{(0,0),(d,0),(0,d)\}$. These are the first values of these invariants, together with the complex and real counterparts:

\begin{table}[h!]
\begin{tabular}{ |c|c|c|c| } 
 \hline
$d$&$N(\Delta_d,0)$ &$W(\Delta_d)$ &$N^{\A^1}_k(\Delta_d)$\\
 \hline
$1$& $1$ & $1$ & $\langle 1\rangle$ \\ 
$2$& $1$ & $1$ & $\langle 1\rangle$ \\ 
$3$& $12$ & $8$ & $8\langle1\rangle+2h$ \\ 
$4$& $620$ & $240$ & $240\langle1\rangle+190h$ \\ 
 \hline
\end{tabular}
\end{table}

In \cite{caporasoharrisformula}, Caporaso and Harris provided a recursion formula for the count of complex plane curves in higher genus, generalizing Kontsevich's recursion formula~\cite{Kontsevich_1994} for the count of complex plane rational curves~$N(\Delta_d,0)$. 
Tropical counterparts of this formula allowing to compute $N(\Delta,0)$ and $W(\Delta)$ by combinatorial means were given in~\cite{GathmannMarkwigCH} and~\cite{IKSCapHar}, respectively.
In the paper \cite{ArithmeticCountOfTropCurves}, we prove a quadratically enriched version of the Caporaso-Harris formula using tropical geometry.

\subsection{Acknowledgements}
Both authors have been supported by the ERC programme QUADAG.  This paper is part of a project that has received funding from the European Research Council (ERC) under the European Union's Horizon 2020 research and innovation programme (grant agreement No. 832833).\\ 
\includegraphics[scale=0.08]{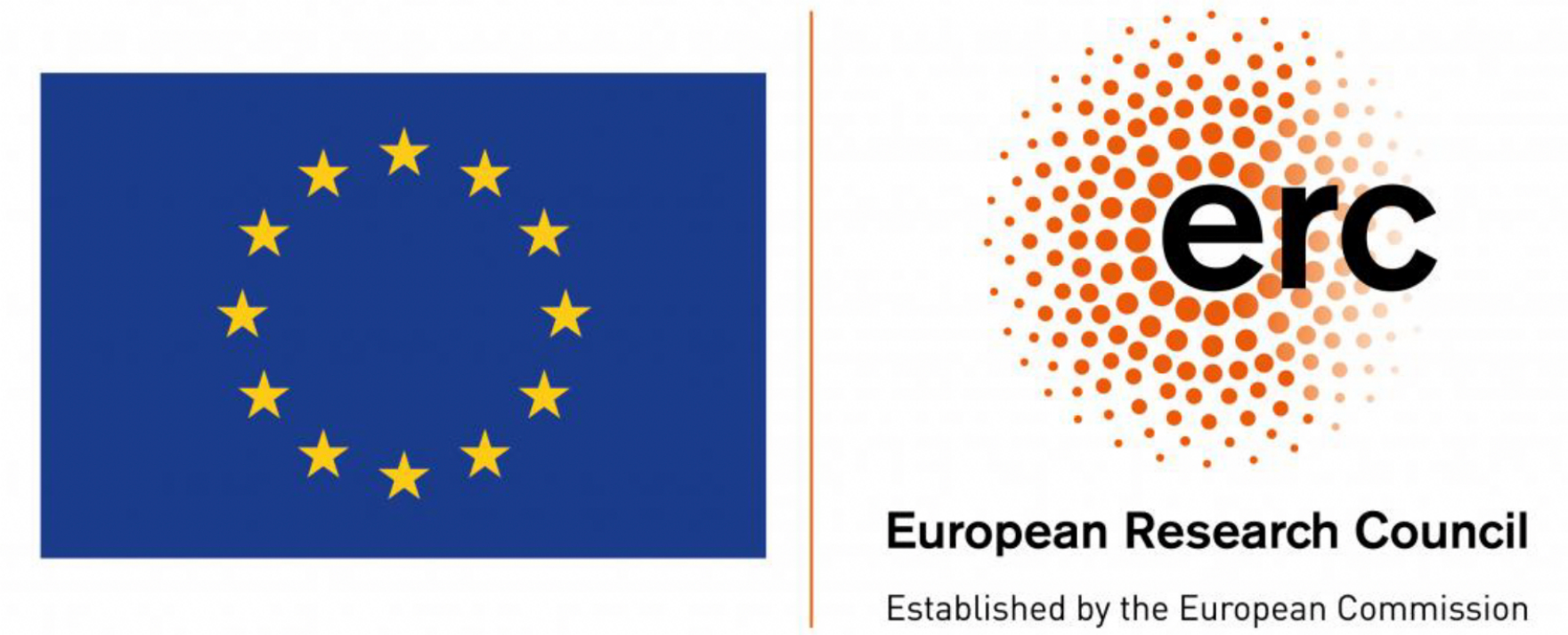}

The first author would like to thank the Universität Duisburg-Essen and
the Università degli Studi di Napoli Federico II for support.
The collaboration on this project started and was partially carried out at the Center for Advanced Study Young Fellows Project ``Real Structures in Discrete, Algebraic, Symplectic and Tropical Geometry''.

We thank Marc Levine for many discussions that really helped understanding the problem and what needs to be done; we are particularly grateful for discussions on the behavior of the invariant with respect to change of the base field and suggestions for generalizations. In addition, we thank him for corrections.
We thank Kirsten Wickelgren for clarifying comments and helpful discussion.
We want to thank Eugenii Shustin for explaining the details and ideas of his work which this paper is built on.
We also want to thank Hannah Markwig for helpful discussions on the structure of the proof of the main theorem which really helped to fill in some details.  
We also thank Felix Röhrle for helpful discussions and very good and detailed comments.
We thank Kris Shaw for helpful discussions, Jake Solomon for helpful comments and Helge Ruddat for discussions on the log approach.

Finally, we thank the anonymous referee for their helpful comments and corrections, which really improved the exposition of this article.

\section{Preliminaries and notation}
\subsection{Notation}
\label{subsection:notation}
We use the following conventions and notation.
\begin{itemize}
    \item $\Delta$ for a convex lattice polygon in $\R^2$,
    \item $k$ a field of characteristic $0$ or greater than the diameter $\operatorname{diam}(\Delta)$ of $\Delta$,
    \item by a polynomial we usually mean Laurent polynomial with exponents $(i,j)\in \Delta\cap\Z^2$ (that is, we also allow negative exponents),
    \item $\K\coloneqq k\Puiseux$ the field over Puiseux series over $k$,
    \item $A$ for a tropical curve with Newton polygon $\Delta$,
    \item $S:\Delta=\Delta_1\cup\ldots\cup\Delta_N$ for the subdivision of $\Delta$ induced by $A$,
    \item $\nu\colon \Delta\rightarrow \Z^2$ for a piecewise linear function inducing the dual subdivision $S$, that means that the graph of $\nu$ projects to $S$,
    \item $\sigma$ for an edge of $S$,
    \item $g$ for the genus,
    \item $n\coloneqq\vert \Delta\cap\Z^2\vert-1-\delta+g=\vert\partial\Delta\cap\Z^2\vert-1+g$ for the number of points in the generic configuration of points,
    \item $p_1,\ldots,p_n$ for the points in $(\K^{\times})^2$,
    \item $x_1,\ldots,x_n$ for the points in $\R^2$ with $x_i=\operatorname{val}(p_i)$ for $i=1,\ldots,n$,
    \item $\sigma_1,\ldots,\sigma_n$ for the marked edges, i.e., the edges in the dual subdivision $S$ of $A$ corresponding to arcs of $A$ with a marking $x_i$, $i\in\{1,\ldots,n\}$,
    \item for a lattice polygon $\Delta'\subset\R^2$ we write $\vert \Delta'\vert$ for its double Euclidean area,
    \item for an edge $\sigma$ of $S$ we let $\vert\sigma\vert$ be its \emph{lattice length}, that is the number of lattice points on $\sigma$ minus $1$ (we often use the letter $m$ for the lattice length),
    \item we say an edge is \emph{odd} (respectively \emph{even}) if it has odd (respectively even) lattice length,
    \item we write $V(\Delta')$ for the set of vertices of a polygon $\Delta'$ and $V(S)$ for the set of vertices of the subdivision, i.e., $V(S)=\bigcup_{l=1}^N V(\Delta_l)$,
    \item we often express a (Laurent) polynomial with coefficients in $\K$ as \\$f(x,y)=\sum_{(i,j)\in \Delta\cap\Z^2}\widetilde{c}_{ij}(t)x^iy^j
    \in \K[x,y]$, and every coefficient $\widetilde{c}_{ij}(t)$ we express as~$(a_{ij}+c_{ij}(t))t^{\nu(i,j)}$, with $a_{ij}\in k$ and such that~$\val(c_{ij})<0$,
    \item $C$ the curve defined by $f$ in $\Tor_\K(\Delta)$,
    \item $f_l(x,y)\coloneqq \sum_{(i,j)\in \Delta_l}a_{ij}x^i y^j$ for $l=1,\ldots, N$,
    \item $C_l$ the curve in $\Tor_k(\Delta_l)$ defined by $f_l$ for $l=1,\ldots,N$,
    \item for an edge $\sigma$ of $\Delta_l$, we call $f_l^{\sigma}\coloneqq\sum_{(i,j)\in \sigma\cap\Z^2}a_{ij}x^iy^j$ the \emph{truncation} of $f_l$ to $\sigma$,
    \item for common edges $\sigma=\Delta_l\cap\Delta_s$ we use the letter $z$ for the point $z\in\Tor(\sigma)\cap C_l\cap C_s$ and $m$ for the multiplicity at which the curves $C_l$ and $C_s$ meet $\Tor_k(\sigma)$,
    \item $\operatorname{Int}(\Delta)$ for the number of interior lattice points of $\Delta$,
    \item we write $E(S)$ for the set of edges of $S$ and $E^*(S)\coloneqq E(S)/\sim$
    for the set of \emph{extended edges} of $S$. That is, we identify parallel edges of a common parallelogram in $S$. More precisely, $\sigma\sim\sigma'$ iff there is a sequence $\Delta_{i_1},\ldots,\Delta_{i_l}$ of parallelograms in $S$, such that  $\sigma$ is an edge of $\Delta_{i_1}$, $\sigma'$ is an edge of $\Delta_{i_l}$ and such that $\Delta_{i_j}$ and $\Delta_{i_{j+1}}$ have a common edge; all these common edges in the sequence are parallel, $\sigma$ is the edge of $\Delta_{i_1}$ parallel to the common edge of $\Delta_{i_1}$ and $\Delta_{i_2}$ and $\sigma'$ is the edge of $\Delta_{i_l}$ parallel to the common edge of $\Delta_{i_{l-1}}$ and $\Delta_{i_l}$.
    We write $[\sigma]$ for the equivalence class of $\sigma$ in $E^*(S)$.
    See Figure \ref{fig:extended edges+curve 2} for an example for this equivalence relation.
    \item For each extended edge $[\sigma]\in E^*(S)$ of lattice length $m\ge 2$, we write $f_{[\sigma]}$ for a \emph{deformation pattern} for this extended edge and call $C_{[\sigma]}$ the curve defined by $f_{[\sigma]}$.
\end{itemize}

\subsection{Motivation from classical topology and $\A^1$-homotopy theory}
\label{subsection:motivation from classical and a1homotopytheory}
Recall from classical algebraic topology that a map $f\colon X\rightarrow Y$ of smooth, closed, connected, oriented real $n$-manifolds, has a well-defined degree. The map $f$ induces a map in $n$-th homology~$f_*\colon H_n(X)\rightarrow H_n(Y)$. The orientations of $X$ and $Y$ define isomorphisms $H_n(X)\cong \Z$ and $H_n(Y)\cong \Z$. The degree of $f$ is $\deg f=f_*(1)\in \Z$ under these isomorphisms. 
Also, recall that if $y$ is a point in the target $Y$ with finitely many points in the preimage, then the degree $\deg f$ equals the sum of local degrees
\[\deg f=\sum_{x\in f^{-1}(y)}\deg_xf.\]
If $y$ is a regular value, then $f$ is locally a homeomorphism around the preimages $x\in f^{-1}(y)$ and thus $\deg_xf\in \{\pm1\}$.

$\A^1$-homotopy theory has been developed to make the techniques from algebraic topology available for smooth algebraic varieties over a field $k$. In \cite{morel}, Fabien Morel defines the \emph{$\A^1$-degree}. An analog of the degree in classical algebraic topology, which is valued in the Grothendieck-Witt ring $\GW(k)$ of $k$, instead of the integers. We briefly introduce this ring in Subsection \ref{subsection:GW(k)}. One can use the $\A^1$-degree to get a well-defined \emph{quadratically enriched count}. That is, a count of $n$-pointed stable maps to a del Pezzo surface $\mathcal{S}$ valued in $\GW(k)$.

Let $\mathcal{S}$ be a del Pezzo surface, $k$ a perfect field of characteristic not equal to $2$ or $3$, $\beta$ an effective Cartier divisor on $\mathcal{S}$, and $\overline{\mathcal{M}}_{0,n}(\mathcal{S}, \beta)$ the moduli space of stable $n$-pointed maps to $\mathcal{S}$ of degree $\beta$.
In \cite{KLSWCount,KLSWorientation} Kass, Levine, Solomon and Wickelgren show that there is a well-defined degree of the evaluation map
\[\operatorname{ev}\colon \overline{\mathcal{M}}_{0,n}(\mathcal{S}, \beta)\rightarrow \mathcal{S}^n\]
valued in $\GW(k)$. Just like in the classical algebraic topology case, this $\A^1$-degree can be expressed as a sum of local $\A^1$-degrees
\[\deg^{\A^1}(\operatorname{ev})=\sum_{x\in \operatorname{ev}^{-1}(y)}\deg^{\A^1}_x(\operatorname{ev}).\]
For example, if $y\in \mathcal{S}^n$ is a suitable general configuration of $n$ $k$-points, then the preimage is comprised of isomorphism classes $x$ of unramified maps $f_x\colon C_x\cong\mathbb{P}^1_k\longrightarrow S$ in the moduli space whose image has only ordinary double points and contains the points in the configuration $y$.
In~\cite{KLSWCount}, it is shown  that
\[\deg^{\A^1}_x(\operatorname{ev})=\Tr_{\kappa(C_x)/k}\left(\Wel^{\A^1}_{\kappa(C_x)}(C_x)\right).\]
Hence,
\[\deg^{\A^1}(\operatorname{ev})=\sum_{x\in\operatorname{ev}^{-1}(y)}\Tr_{\kappa(C_x)/k}\left(\Wel^{\A^1}_{\kappa(C_x)}(C_x)\right)\eqqcolon N^{\A^1}_k(\mathcal{S},\beta).\]
The left hand side of this equation does not depend on the chosen point configuration of~$k$-points, which proves invariance of $N^{\A^1}_k(\mathcal{S},\beta)$.

\subsection{The Grothendieck-Witt ring}
\label{subsection:GW(k)}
We recall the basic definitions and some facts about the Grothendieck-Witt ring.

\begin{df}
\label{df:GW}
The Grothendieck-Witt ring $\GW(R)$ of a ring $R$ is the group completion of the semi-ring of isometry classes of non-degenerate symmetric bilinear forms over $R$ under the direct sum $\oplus$ and tensor product $\otimes$. 
\end{df}

We will do several explicit calculations in $\GW(R)$ mainly when $R=k$ is a field. In this case $\GW(k)$ has a nice presentation. It is generated by the isometry classes of non-degenerate symmetric bilinear forms on a $1$-dimensional $k$-vector space
\[\langle a\rangle\coloneqq [(x,y)\mapsto axy],\]
for $a\in k^{\times}$. Subject to the following relations
\begin{enumerate}
    \item $\qinv{ab^2}=\qinv{a}$,
    \item $\qinv{ab}=\qinv{a}\qinv{b}$ for $a,b\in k^\times$,
    \item $\qinv{a}+\qinv{b}=\qinv{a+b}+\qinv{ab(a+b)}$ for $a,b\in k^\times$ satisfying $a+b\neq 0$.
\end{enumerate}

\begin{df}
We write $h$ for the \emph{hyperbolic form}. That is, the form on a $2$-dimensional $k$-vector space (or free rank $2$ R-module over $R$ when $R$ is not a field), with Gram matrix \[\begin{pmatrix}0&1\\1&0\end{pmatrix}.\]
\end{df}
When $R$ is a field $k$, one can show that for  $a\in k^\times$, we have that
\begin{equation}
\label{eq:hyperbolicform}
h=\qinv{1}+\qinv{-1}=\qinv{a}+\qinv{-a}.
\end{equation}
Let $\beta:V\times V\rightarrow k$ be a non-degenerate symmetric bilinear form. The \emph{rank} of $\beta$ is the rank of the $k$-vertor space $V$. The rank extends to a homomorphism
\[\operatorname{rk}\colon \GW(k)\rightarrow \Z.\]

\begin{ex}
    Over $\C$, the rank defines an isomorphism of rings~$\operatorname{rk}\colon\GW(\C)\cong \Z$.
\end{ex}

\begin{ex}
    Over $\R$, the class of non-degenerate symmetric form is completely determined by its rank and its \emph{signature}. Recall that one can diagonalize a non-degenerate symmetric bilinear form over $\R$ such that the associated Gram matrix is diagonal with only $1$'s and $-1$'s on the diagonal. The signature of the form equals the number of $1$'s minus the number of $-1$'s.
\end{ex}

\begin{rmk}
\label{rmk:GW(K)}
Let $k$ be a field of characteristic not equal to $2$ and $\K=k\Puiseux$ the field of Puiseux series over $k$.
    Then the map
    \[\K^\times/(\K^\times)^2\longrightarrow k^\times/(k^\times)^2\colon a_0t^{q_0}+a_1t^{q_1}\ldots \longmapsto a_0,\]
    is a bijection (we assumed the first coefficient $a_0\neq0$). This bijection induces an isomorphism of rings (cf.  \cite[Theorem 4.7]{mps})
    \[\GW(\K)\cong \GW(k).\]
\end{rmk}
A non-degenerate symmetric bilinear form $\beta:V\times V\rightarrow R$ over a ring $R$ is \emph{split} if there exists a submodule $N\subset V$ such that $N$ is a direct summand of $V$ and $N$ is equal to its orthogonal complement $N^{\bot}$.
We say that two non-degenerate symmetric bilinear forms~$\beta\colon V\times V\rightarrow R$ and~${\beta'\colon V'\times V'\rightarrow R}$ are \emph{stably equivalent} if there exist split symmetric bilinear forms $s$ and $s'$ such that $\beta\oplus s\cong\beta'\oplus s'$.
\begin{df}
The \emph{Witt ring} of $R$ is the set of classes of stably equivalent non-degenerate symmetric bilinear forms with addition the direct sum $\oplus$ and multiplication the tensor product $\otimes$.
\end{df}

In this paper, we will do some of the computations in the Witt ring instead of the Grothendieck-Witt ring. We do this without loss of information because of the following remark.

\begin{rmk}
\label{rmk:GWandW}
If $k$ is a field of characteristic different than $2$, then 
the split non-degenerate symmetric bilinear forms are exactly the multiples of the hyperbolic form $h$.
Recall that in this case, we have that $\langle a\rangle +\langle-a\rangle=h$ for any unit $a$. Hence, 
\[\operatorname{W}(k)=\frac{\GW(k)}{\Z\cdot h}.\]
More generally, if $R$ is local and $2$ is invertible, then an element of $\GW(R)$ is completely determined by its rank and its image in $\W(R)$. In case $R=\R$, the image of a non-degenerate symmetric bilinear in $\W(\R)\cong \Z$ coincides with its signature.
\end{rmk}

Assume that $R$ is a commutative ring. We are particularly interested in the case when~$R$ is a finite étale $k$-algebra.
For a finite projective $R$-algebra $L$ one can define the trace map~$\Tr_{L/R}:L\rightarrow R$ that sends $b\in L$ to the trace of the multiplication endomorphism~$m_b(x)=b\cdot x$.
If the algebra~$L$ is étale over $R$, then this induces \emph{trace maps}~$\Tr_{L/R}\colon \GW(L)\rightarrow \GW(R)$ and~$\Tr_{L/R}\colon \operatorname{W}(L)\rightarrow \operatorname{W}(R)$, which send the class of a bilinear form $\beta\colon V\times V\longrightarrow L$ over $L$ to the form
\[V\times V\xrightarrow{\quad\beta\quad}L\xrightarrow{\;\operatorname{Tr}_{L/R}\;}R\]
over $R$.
We will compute several trace forms in the proof of our main result. So we collect some facts about the trace form here.
Let $\etale$ be a finite \'{e}tale $R$-algebra. 
\begin{enumerate}
    \item If $R=k$ is a field, then $\etale=L_1\times \ldots\times L_s$ for some finite separable field extensions~$L_1,\ldots,L_s$ of~$k$ and the trace map $\Tr_{\etale/k}:\etale\rightarrow k$ equals the sum of  field traces $\Tr_{\etale/k}=\sum_{i=1}^s\Tr_{L_i/k}$.
    \item $\Tr_{\etale/R}$ is $R$-linear.
    \item If $F$ is a finite \'{e}tale $\etale$-algebra, then 
    \[\Tr_{F/R}=\Tr_{\etale/R}\circ\Tr_{F/\etale}.\]
    \item Let $\etale$ be a finite \'{e}tale $R$-algebra of rank $m$.
Then $\Tr_{\etale/R}(h)=m\cdot h$.
\end{enumerate}

The following Proposition can be found in \cite[Proposition 2.13]{JPP}.

\begin{prop}
\label{prop:traces}
Let $L$ be a finite \'{e}tale $k$-algebra and let $\etale=\frac{L[x]}{(x^m-D)}$, for some $D\in L^\times$. Further, assume that $\operatorname{char}k$ does not divide $m$. Then for $a\in L$, we get that
\begin{enumerate}
    \item $\Tr_{\etale/L}(\qinv{ m\cdot a})=\begin{cases}
\qinv{ a} +\frac{m-1}{2}h & m\text{ odd,} \\
\qinv{ a}+\qinv{  a\cdot D}+\frac{m-2}{2}h &  m \text{ even.}
\end{cases}$
    \item $\Tr_{\etale/L}(\qinv{ m\cdot a\cdot x})=\begin{cases}
\qinv{ a\cdot D} +\frac{m-1}{2}h & m\text{ odd,} \\
\frac{m}{2}h &  m \text{ even.}
\end{cases}$
\end{enumerate}
\end{prop}

\subsection{Toric geometry}
Tropical geometry provides a  language for studying toric varieties. Tropical varieties can be seen as the tropicalizations of their corresponding toric varieties. The tropicalization process assigns to each point in the toric variety a tropical point, capturing some information of the original variety. 
Moreover, toric morphisms naturally induce morphisms on the tropical side. This correspondence between toric morphisms and tropical morphisms enables us to translate the algebraic count to a combinatorial count. 
We will detail the tropicalization in Section~\ref{section:tropicalization}.

The toric variety $\Tor_k(\Delta)$ associated to a convex lattice polyhedron $\Delta$ is the toric variety associated with the fan~$\Sigma$ consisting of the cones~$\sigma$ over the proper faces of~$\Delta$.
For each cone~$\sigma$, its affine toric variety is the spectrum of the semigroup algebra $k[S_{\sigma}]$, where
\[S_{\sigma}=\{u\in\operatorname{Hom}(\Z^n,\Z)\mid u(v)\geq0, \forall v\in\sigma\}.\]
These affine toric varieties are glued by the morphisms induced by the inclusion of the algebras $k[S_{\sigma}]\subset k[S_{\tau}]$, whenever~$\tau$ is a face of~$\sigma$.
This construction gives rise to a variety with a torus action, where the torus orbits can be described by the boundary of $\Delta$.

The geometry and combinatorics of the toric variety is hence determined by their associated polyhedra $\Delta$.
The cones in~$\Sigma$ correspond bijectively to the orbits of the torus action on $\Tor_k(\Delta)$.
A morphism between two toric varieties $\Tor_k(\Delta)$ and $\Tor_k(\Delta')$ is called toric if it is induced by a linear map $f$ of the lattice such that the image of every cone of $\Sigma$ is contained in a cone of $\Sigma'$.

As a basic example, let us consider the cone $\sigma=\R_{\ge 0}$. Hence the semigroup $S_{\sigma}=\N$ is the polynomial ring in one variable $k[x]$ and we have that $\Tor_k(\R_{\ge 0})=\A^1_k$ by taking its spectrum.

A key example of this construction is the toric surface associated to the lattice triangle~$\Delta_d=\operatorname{Conv}\{(0,0),(d,0),(0,d)\}$, which is the projective plane. The affine toric surfaces associated to the cones coincide with the standard affine charts, and moreover, its boundary divisors are the coordinate lines. This construction is independent of the size of the triangle $d$, however this will play a role as the degree of the curves in the enumeration.

\subsection{Useful results}
Lastly, we recall some useful well-known results we will use in this paper.
We start by recalling Pick's theorem, which relates the lattice points of a lattice polygon to its area.

\begin{thm}[Pick's theorem]
\label{thm:Pick}
    Let $P$ be a lattice polygon in $\R^2$. Let $A$ be its area, $b$ be the number of lattice points on the boundary and $i$ be the number of interior lattice points. Then it holds that
    \[A=i+\frac{b}{2}-1.\]
\end{thm}

For some of our computations, we need the following identity.
\begin{lm}
\label{lm:sineIdentity} 
Let $m$ be a positive integer that is not divisible by the characteristic of the field $k$. Let $\zeta_m$ be a primitive $m$-root of unity, possibly in a field extension $k(\zeta_m)$. Then, we have that
\[\prod_{j=1}^{m-1} \left( 1-\zeta_{m}^j\right)=m\in k.\]
\end{lm}
\begin{proof}
Let $p(x)=\prod_{j=1}^{m-1}  \left( x-\zeta_{m}^j\right)$ be the monic polynomial whose roots are the non trivial $m$-roots of unity. We have that $p(x)(x-1)=x^m-1$. Hence, 
\[p(x)=\frac{x^m-1}{x-1}=\sum_{j=0}^{m-1}x^j,\] 
and our assertion follows from evaluating at $x=1$.
\end{proof}

\section{Tropicalization, tropical limit and refined tropical limit}
\label{section:tropicalization}

\subsection{Tropicalization}
\label{subsection:tropicalization}
\begin{figure}
    \begin{tabular}{ccc}
    \includegraphics[scale=0.12]{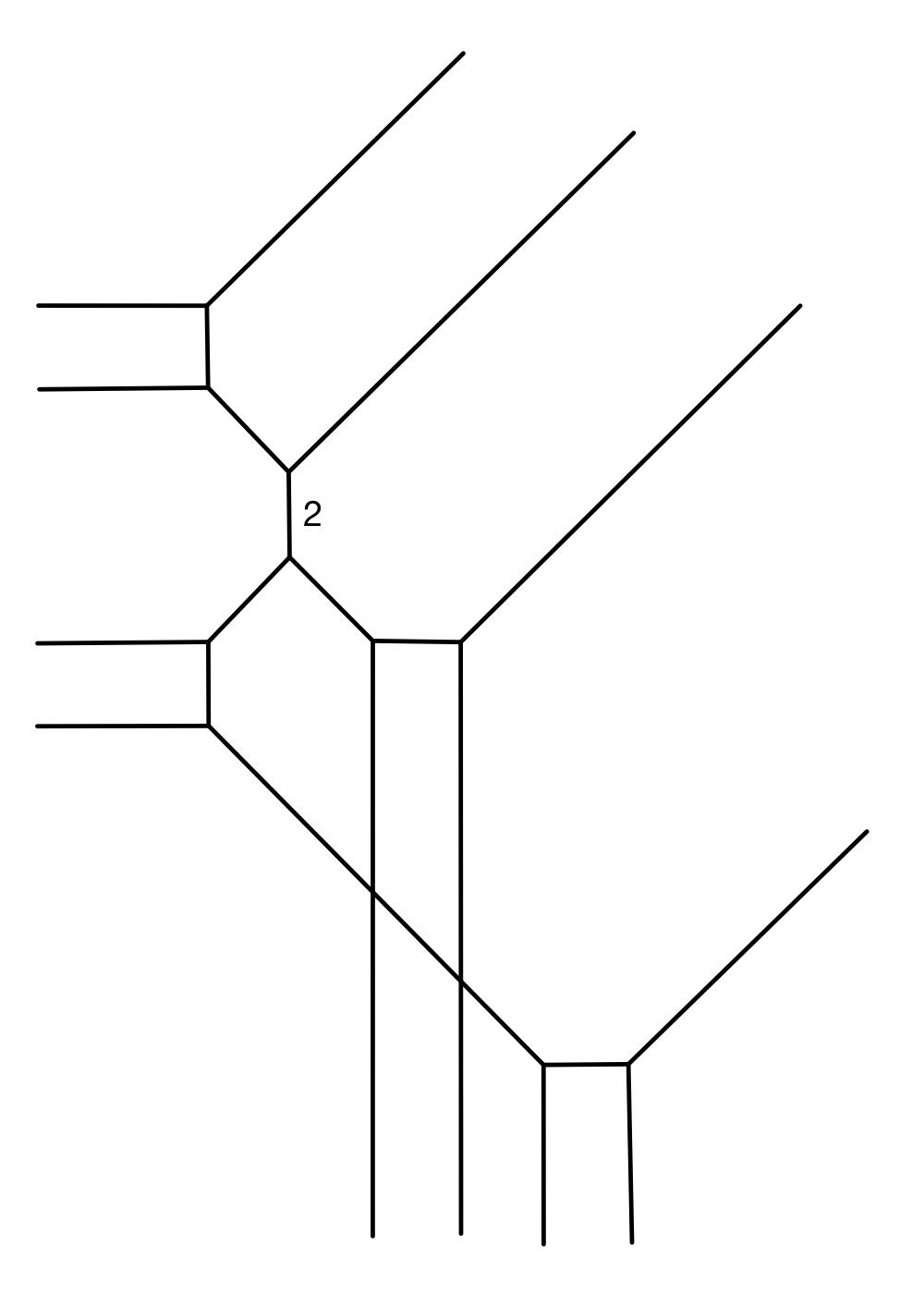}&
    \includegraphics[scale=0.12]{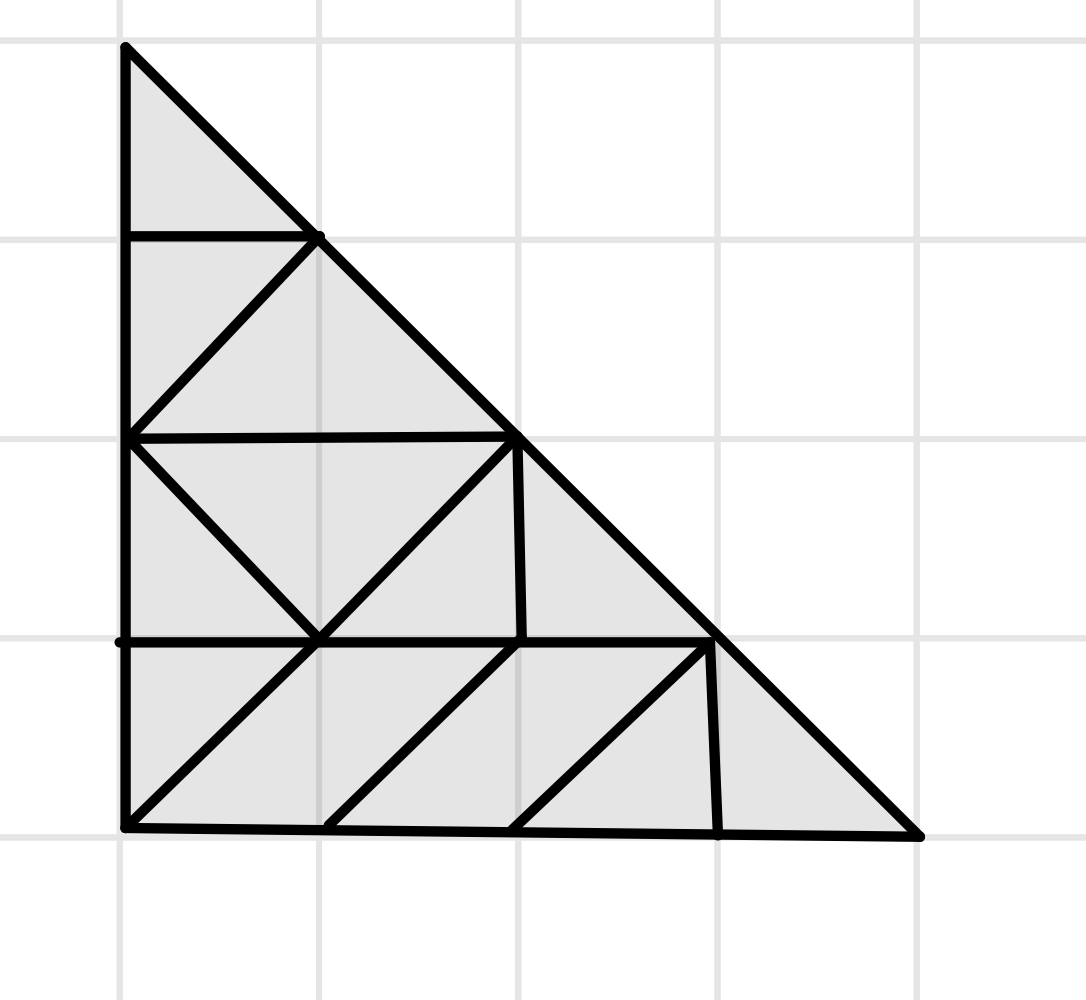}
    \end{tabular}
    \caption{A simple tropical curve and its dual subdivision.\label{fig:dualsubdiv}}
    \end{figure}

We briefly recall the definitions of a tropical curve, its dual subdivision and tropicalization.
Let $\mathbb{T}\coloneqq \R\cup\{-\infty\}$ be the tropical semifield with the two tropical operations
tropical addition
\[``x+y"\coloneqq \max\{x,y\}\]
and tropical multiplication
\[``x\cdot y"\coloneqq x+y.\]
A \emph{tropical (Laurent) polynomial} in two variables $x$ and $y$ is of the form
\[f(x,y)=``\sum a_{ij}x^iy^j"=\max\{a_{ij}+ix+jy\},\]
with $a_{ij}\in \mathbb{T}$ and $a_{ij}\neq -\infty$, for only finitely many $(i,j)\in \Z^2$. Note that this defines a piecewise linear function $\R^2\rightarrow \R$. Here, we also allow negative exponents $i,j\in \Z$.
The tropical vanishing locus $V^{\operatorname{trop}}(f)$ of a tropical polynomial is the non-smooth locus of the piecewise linear function $f$; or alternatively, the set of points $(x,y)\in \R^2$ where the maximum is attained at least twice
\[V^{\operatorname{trop}}(f)\coloneqq\{(x,y)\in\R^2:\text{maximum of $f$ at $(x,y)$ is attained at least twice}\}.\]
These tropical vanishing loci of tropical polynomials are piece-wise linear graphs with weights on edges defined as follows. For a given segment where the maximum of the polynomial is attained twice, namely by monomials~$a_{ij}x^{i}y^{j}$ and~$a_{i'j'}x^{i'}y^{j'}$, its weight is given by the maximum of the $\gcd(i'-i,j'-j)$'s where the maximum ranges over all pairs of $2$-tuples $\{(i,j),(i',j')\}$ corresponding to pairs of monomials where the maximum is attained on the given segment. More precisely, one attaches the following \emph{weight $\omega(e)$} to an edge $e$ of the tropical curve
\[\omega(e)\coloneqq \max\{\gcd(\vert i-i'\vert,\vert j-j'\vert)\},\]
where the maximum runs over all pairs $(i,j)$ such that $f(x,y)=a_{ij}+ix+jy$ for all points ${(x,y)\in e}$.
We only write weights greater than $1$. So edges of weight $1$ have no labelling.
\begin{ex}
\label{ex:tropicalline}
    We start with an easy example, namely a tropical line, that is the tropical vanishing locus of a degree $1$ tropical polynomial. Let 
    \[f(x,y)=``1+2x+(-1)y"=\max\{1, x+2,y-1\}.\]
    Then the tropical vanishing locus of $f$ is displayed in Figure \ref{fig:tropical line}.
    \begin{figure}
    \includegraphics[scale=0.15]{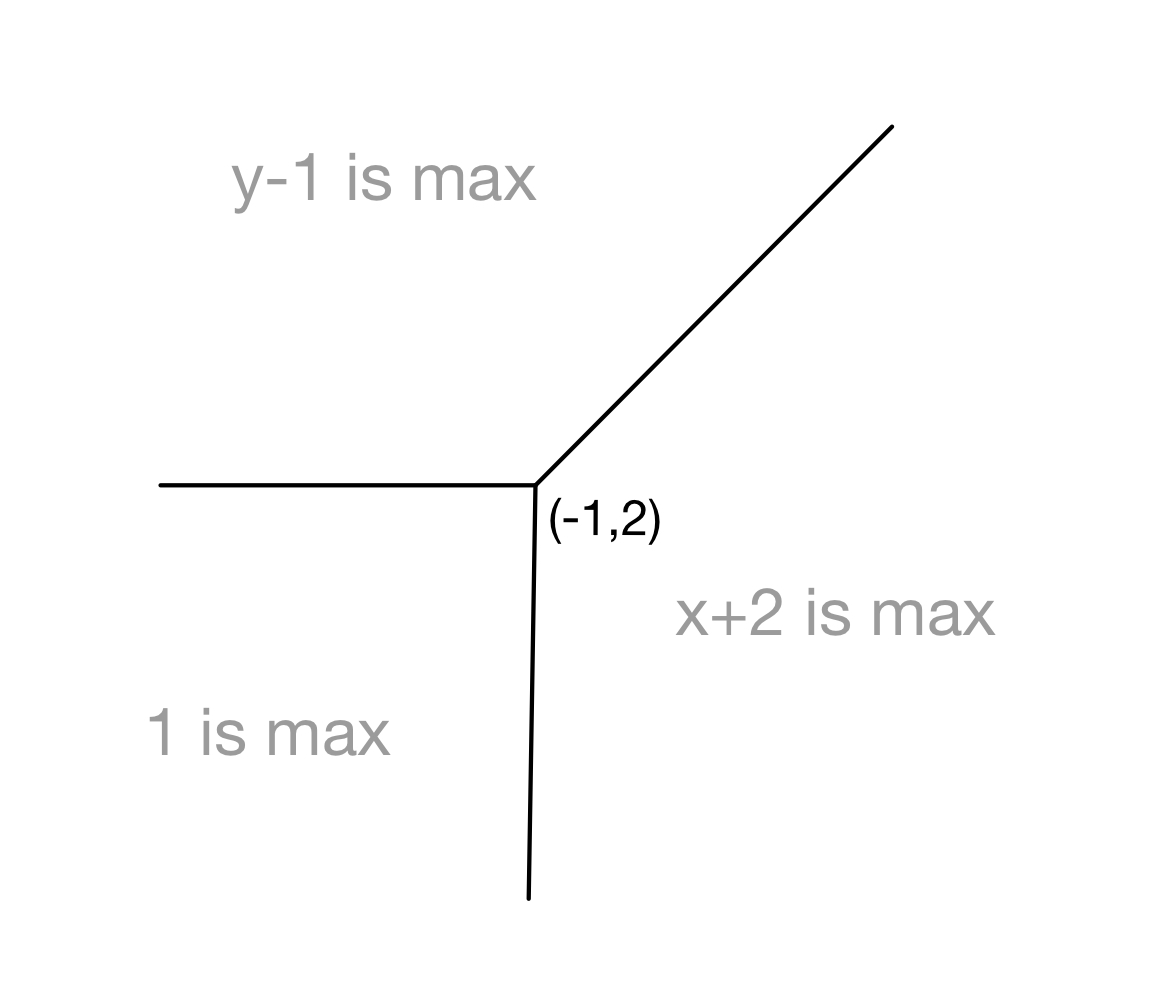}
    \caption{A tropical line.\label{fig:tropical line}}
    \end{figure}
    In this example all the weights of the edges are equal to $1$.
\end{ex}
\begin{ex}
    The first picture in Figure \ref{fig:dualsubdiv} shows a tropical curve with one edge of weight~$2$.
\end{ex}
The following definition coincides with the definition of a tropical curve in $\R^2$ defined algebraically as the tropical vanishing locus of a tropical polynomial above. It is known in the literature as an \emph{embedded tropical curve}.
\begin{df}
\label{df:tropcurve}
A \emph{tropical curve} $C$ is a finite weighted graph $(V,E,\omega)$ embedded in $\R^n$, where the set~$E=E^{\infty}\cup E^{\circ}$ is the disjoint union of univalent edges $E^{\infty}\subset V$ and non-directed edges $E^{\circ}\subset \{e\subset V\mid \Card(e)=2\}$, such that every edge $l\in E^{\infty}$ embeds into a ray of an integer line, i.e., a line given by $\bar{a}_e\cdot t+\bar{b}_e$ with $\bar{a}_e\in\Z^n\setminus\{\bar{0}\},\bar{b}_e\in\Q^n$; every edge $e\in E^{\circ}$ embeds into a segment of the graph of an integer line; and every vertex $v\in V$ satisfies
the \emph{balancing condition}
\begin{equation}
\label{eq:balancing condition}
\sum_{e\in E, v\in e} \omega(e)\cdot \mathrm{u}_e=0,
\end{equation}
where $\mathrm{u}_e=\displaystyle\frac{\pm 1}{\gcd(\bar{a}_e)}\bar{a}_e$ is oriented outwards from the vertex $v$, and $\omega:E\lra\Z$ is the non-negative \emph{weight function}.
We call $\bar{a}_e$ a \emph{director} vector of $e$ and $\mathrm{u}_e$ a \emph{primitive} vector of $e$ at $v$.
When drawing a tropical curve we write the weights not equal to $1$ next to the edges.
The \emph{genus} of a tropical curve is the first Betti number of the graph $(V,E,\omega)$ (before the embedding).
\end{df}

\begin{df}
The \emph{degree} of the tropical curve~$C$ is the multiset of primitive vectors associated to its legs $\{\mathrm{u}_l\mid l\in E^{\infty}\}$ counted with weights.
\end{df}

The Newton polygon of a tropical curve defined by a tropical polynomial $f$ is the convex hull in $\R^2$ of the indices $(i,j)\in \Z^2$ for which the coefficients $a_{ij}\neq -\infty$, i.e.,
\[\Delta(f)\coloneqq \operatorname{Conv}\{(i,j)\in \Z^2:a_{ij}\neq-\infty\}.\]
The Newton polygon can be obtained from the degree of the curve, it equals a polygon whose dual fan has rays in the directions of the elements of the degree with lattice lengths given by the multiplicities of the elements of the degree.
Henceforth, we fix a convex lattice polygon $\Delta\subset \R^2$.
We recall how to obtain a tropical curve from an algebraic curve in $\Tor_\K(\Delta)$ where $\K=k\Puiseux$ is the field of Puiseux series over $k$.
Recall that there is a non-Archimedean valuation on $\K$
\[\val:\K\rightarrow\R\cup\{-\infty\}\colon \widetilde{a}=a_0t^{q_0}+a_1t^{q_1}\ldots\mapsto -q_0.\]
To a polynomial \begin{equation}
\label{eq:f}
f(x,y)=\sum_{(i,j)\in\Delta\cap\mathbb{Z}^2}\widetilde{c}_{ij}x^iy^j\in\K[x,y],
\end{equation}
with $\widetilde{c}_{ij}\in\K$, we associate the following tropical polynomial
\[f^{\operatorname{trop}}(x,y)\coloneqq``\mkern-15mu\sum_{(i,j)\in\Delta\cap \Z^2}\val(\widetilde{c}_{ij})x^iy^j"=\max_{(i,j)\in \Delta\cap\Z^2}\{\val(\widetilde{c}_{ij})+ix+jy\}\in\mathbb{T}[x,y].\]
We use the same name for the formal variables of $f$ and $f^{\operatorname{trop}}$, but in this context, the elements of $\mathbb{T}$ should be considered as the valuations of elements of $\K$.

Let
$A_f\coloneqq V^{\operatorname{trop}}(f^{\operatorname{trop}})$
be the tropical curve defined as the tropical vanishing locus of~$f^{\operatorname{trop}}$, endowed with its weights.
The curve $A_f$ has Newton polygon $\Delta$.
We say that~$A_f$ is the \emph{tropical curve associated with $f$}, and that the curve $C$ defined by $f$ \emph{tropicalizes} to~$A_f$.
This definition of $A_f$ agrees with the definition in \eqref{eq:Kapranov} by Kapranov's theorem.
The tropicalization of a family of curves preserves the genus (cf. e.g. \cite{Mikhalkin}). In particular, if the curve $C$ defined by $f$ is rational, then the associated tropical curve $A_f$ has genus $0$.

\begin{rmk}
The tropicalization already gives one direction in the correspondence theorem: Given a genus $g$ curve $C\subset \Tor_\K(\Delta)$, tropicalizing yields a unique genus $g$ tropical curve with Newton polygon $\Delta$. So the main task in proving a tropical correspondence theorem, is to find all algebraic curves that tropicalize to a given tropical curve.
\end{rmk}

Given a tropical curve $A$ with Newton polygon $\Delta$, set \[\widetilde{\Delta}\coloneqq\operatorname{Conv}\{(i,j,-\val(\widetilde{c}_{ij}))\in \R^3: (i,j)\in \Delta\cap\mathbb{Z}^2\}\subset\R^3\]
and define the piecewise linear function
\begin{equation}
\label{eq:nu}
\nu\colon\Delta\longrightarrow \R \colon (x,y)\longmapsto\nu(x,y)\coloneqq\min\{\gamma:(x,y,\gamma)\in \widetilde{\Delta}\}.\end{equation}
The linearity domains of $\nu$ define a subdivision $S$ of $\Delta$, called the \emph{dual subdivision} (see Figure \ref{fig:dualsubdiv}).
This subdivision describes $\Delta$ as the union $\Delta=\Delta_1\cup\ldots\cup\Delta_N$, where $\nu\vert_{\Delta_i}$ is linear. The subdivision $S$ satisfies the following duality properties:
\begin{itemize}
    \item the components of $\R^2\setminus A$ are in one to one correspondence with the vertices $V(S)$ of $S$;
    \item the edges of $A$ are in one to one correspondence with the vertices of $S$, dual edges are orthogonal, and for an edge $e$ of $A$ of weight $\omega(e)$, its dual edge has lattice length $\omega(e)$;
    \item the vertices of $A$ are in one to one correspondence with the polygons $\Delta_1,\ldots,\Delta_N$, and the valency of a vertex of $A$ is equal to the number of sides of the dual polygon.
\end{itemize}

We recall the following definition from \cite{Mikhalkin}.
\begin{df}
    A tropical curve $A$ with Newton polygon $\Delta$ and dual subdivision $S$ is \emph{nodal} if
the subdivision $S$ consists of triangles and parallelograms. If additionally, all lattice points on the boundary are vertices of $S$ (or equivalently all unbounded edges of $A$ have weight $1$), the tropical curve is called \emph{simple}.

\end{df}
\begin{rmk}
\label{rmk:nu}
    Let $\nu_f\colon \Delta\rightarrow \R$ be the piecewise linear function assigned to a simple tropical curve $A_f$ defined by a polynomial $f$ as in \eqref{eq:f}. Then the coefficients $\widetilde{c}_{ij}$ of $f$ can be written as $\widetilde{c}_{ij}=(a_{ij}+c_{ij}(t))t^{\nu_f(i,j)}$ where $a_{ij}\in k$ and $c_{ij}\in \K$ with $\val(c_{ij})<0$. 
\end{rmk}
\begin{rmk}
    For a nodal tropical curve in $\R^2$ one can read of its genus as follows. It equals the first Betti number of the embedded graph minus the number of $4$-valent vertices, or equivalently, minus the number of parallelograms in the dual subdivision.
\end{rmk}

\begin{df}
    The tropical curves with a given dual subdivision $S$ of $\Delta$ are parametrized by a convex polyhedron, whose dimension is called the \emph{rank} of $A$ and is denoted by~$\operatorname{rk}(A)$. In particular, a tropical curve of degree $\Delta$, with dual subdivision $S$, is determined by~$\operatorname{rk}(A)$ points in the plane $\R^2$.
\end{df}
The rank of a tropical curve is finite since the configuration space of its vertices is finite dimensional. 
By considering deformations of the edges at infinity and taking into account these deformations get transferred through the balancing condition at every vertex, 
we have that $\operatorname{rk}(A)\ge\vert \partial\Delta\cap\Z^2\vert-1$.
The following definition is based on \cite[Definition~4.7]{Mikhalkin}.
\begin{df}
    A configuration of $n=\vert\partial\Delta\cap\Z^2\vert-1+g$ points  
    in $\R^2$ is \emph{$\Delta$-generic} if for any genus $g$ tropical curve $A$ with Newton polygon $\Delta$ passing through these points, we have that~$A$ is simple and has rank $\operatorname{rk}(A)=n$.
\end{df}
In \cite{Mikhalkin} is proved that the $\Delta$-generic configurations of real points form a dense set. However, many of the arguments hold for configurations of points in $\Q^2\subset\R^2$. In particular, Proposition 4.11 in \cite{Mikhalkin} proves that the dimension of the space of configurations~$\mathcal{P}\in\operatorname{Sym}^n(\Q^2)$ that are not $\Delta$-generic has dimension at most $2n-1$.
Hence, we can assume that the configuration of points $p_1,\ldots,p_n$ are chosen in such way that the points $x_1,\ldots,x_n\in\Q^2\subset\R^2$ are in tropical general position, and thus we can restrict our discussion to simple tropical curves.

\begin{ex}
    The Newton polygon of a tropical line as in Example \ref{ex:tropicalline} is the triangle~$\Delta_1=\operatorname{Conv}\{(0,0),(1,0),(0,1)\}$ and thus its rank equals $\vert\partial\Delta_1\cap\Z^2\vert-1=3-1=2$. So a tropical line is determined by $2$ points in $\Delta_1$-general position. Indeed, two of the rays of tropical line are fixed by the points and the third one is determined by the balancing condition \eqref{eq:balancing condition} as illustrated in Figure \ref{fig:tropicallinethrough2points}.

    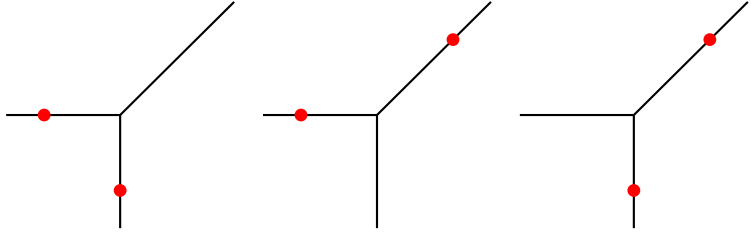
\begin{figure}
    \begin{tabular}{ccc}
    \begin{tikzpicture}
    \draw[thick] (0,0) -- (-1.5,0);
    \draw[thick, red, fill=red] (-1,0) circle (2pt);
    \draw[thick] (0,0) -- (0,-1.5);
    \draw[thick, red, fill=red] (0,-1) circle (2pt);
    \draw[thick] (0,0) -- (1.5,1.5);
\end{tikzpicture}&
    \begin{tikzpicture}
    \draw[thick] (0,0) -- (-1.5,0);
    \draw[thick, red, fill=red] (-1,0) circle (2pt);
    \draw[thick] (0,0) -- (0,-1.5);
    \draw[thick] (0,0) -- (1.5,1.5);
    \draw[thick, red, fill=red] (1,1) circle (2pt);
\end{tikzpicture}&
    \begin{tikzpicture}
    \draw[thick] (0,0) -- (-1.5,0);
    \draw[thick] (0,0) -- (0,-1.5);
    \draw[thick, red, fill=red] (0,-1) circle (2pt);
    \draw[thick] (0,0) -- (1.5,1.5);
    \draw[thick, red, fill=red] (1,1) circle (2pt);
\end{tikzpicture}
    \end{tabular}
\caption{Tropical lines determined by two points in $\R^2$.\label{fig:tropicallinethrough2points}}
    \end{figure}
\end{ex}

\subsection{Toric degenerations and the tropical limit}
\label{subsection:toricdegeneration}
We have seen that to an algebraic curve in $\Tor_\K(\Delta)$ we can assign an associated tropical curve. Now we have to figure out, how much additional information one needs, to reconstruct an algebraic curve in $\Tor_\K(\Delta)$ from a tropical curve. The first information we need is the following.
We will see that one can extend a given curve $C$ defined by $f$ to a family of curves $\widetilde{C}$ with generic fiber equal to $C$ and special fiber of the form
\[C_0=\bigcup_{l=1}^NC_l\subset \bigcup_{l=1}^N\Tor_k(\Delta_l),\]
where the $\Delta_l$ are the polygons in the dual subdivision $S_f=\Delta_1\cup\ldots\cup\Delta_N$ associated to~$C$. In this subsection we recall how to find $C_0$ and recall some of its properties from \cite{Shustin}, \cite{NishinouSiebert} and \cite{Tyomkin}.

\begin{df}
    A flat toric morphism
    $\pi:\mathcal{X}\rightarrow \A^1_k$
    is called a \emph{toric degeneration}.
\end{df}

\begin{prop}
    A tropical curve $A$ with dual subdivision $S\colon\Delta=\Delta_1\cup\ldots\cup\Delta_N$ defines a toric degeneration $\pi:\mathcal{X}\rightarrow \A^1_k$ such that
    \begin{itemize}
        \item for $t\neq 0$ the fiber $X_t=\pi^{-1}(t)\cong \Tor_{\kappa(t)}(\Delta)$,
        \item the special fiber is given by $X_0=\pi^{-1}(0)\cong \bigcup_{i=1}^N\Tor_k(\Delta_i)$ glued along their toric boundary divisors, that is, if $\Delta_i$ and $\Delta_j$ have a common edge $\sigma=\Delta_i\cap\Delta_j$, then~$\Tor_k(\Delta_i)$ and $\Tor_k(\Delta_j)$ are glued along $\Tor_k(\sigma)$.
    \end{itemize}
\end{prop}

\begin{proof}
This is \cite[$\S3$]{NishinouSiebert}, \cite[$\S2.2$]{ArguzBousseau} or \cite[Proposition 3.5]{Tyomkin} in the algebraic setting and in \cite[$\S2.3$]{Shustin} in the analytic setting.
The toric morphism is induced by the following:
After a parameter change $t\mapsto t^m$ one can assume that all exponents of the coefficients of $f$ are integers. 
    As in \cite[$\S2.3$]{Shustin}, set
    \[\overline{\Delta}\coloneqq \{(\alpha,\beta,\gamma)\in \R^3: (\alpha,\beta\in \Delta)\text{, }\gamma\ge\nu_f(\alpha,\beta)\},\]
    and let $\mathcal{X}\coloneqq \Tor_k(\overline{\Delta})$. The map $\overline{\Delta}\rightarrow \R$ defined by $(\alpha,\beta,\gamma)\mapsto \gamma$ gives rise to the map $\pi\colon \mathcal{X}=\Tor_k(\overline{\Delta})\rightarrow \Tor_k(\R_{\ge 0})=\A^1_k$ which defines a toric degeneration with the properties listed in this Proposition. 

\end{proof}

Let us express $f\in\K[x,y]$ as
\[f(x,y)=\sum_{(i,j)\in \Delta\cap \Z^2}(a_{ij}+c_{ij}(t))t^{\nu(i,j)}x^iy^j,\]
such that $a_{ij}\in k$ and $c_{ij}(t)\in \K$ such that $\val(c_{ij}(t))<0$. 
Set
\begin{equation}
\label{eq:fs}
f_l(x,y)\coloneqq \sum_{(i,j)\in \Delta_l\cap\Z^2}a_{ij}x^iy^j,
\end{equation} 
for $l=1,\ldots,N$,
and let $C_l$ be the curve in $\Tor_k(\Delta_l)$ defined by $f_l$.
Let $A_f$ be the associated tropical curve and let $\pi_f\colon \mathcal{X}\rightarrow \A^1_k$ be the associated toric degeneration.
As before, we assume that the exponents in the coefficients of $f$ are integers and hence $C$ defines a curve over $k((t))$.
\begin{prop}
\label{prop:family of curves}
    Up to a base change $\A^1\rightarrow \A^1\colon t\mapsto t^b$, the curve
    $C$, defined by~$f$ as above, gives rise to a family of curves $\widetilde{C}\rightarrow \operatorname{Spec}k[[t]]$ with generic fiber $C$, which fits into the following commutative diagram
    \[\begin{tikzcd}
  \widetilde{C} \arrow[r] \arrow[d]
    & \mathcal{X} \arrow[d] \\
  \operatorname{Spec} k[[t]] \arrow[r]
&  \A^1_k \end{tikzcd}\]
    where the bottom map is given by $k[t]\rightarrow k[[t]]$.
Furthermore, the special fiber of $\widetilde{C}$ is 
    \[C_0=\bigcup_{l=1}^NC_l\subset X_0=\bigcup_{l=1}^N\Tor_k(\Delta_l),\]
    where the curve $C_l$ is defined by the $f_l$ as defined in \eqref{eq:fs}, for $l=1,\ldots,N$.
\end{prop}
\begin{proof}
    This is \cite[Proposition 3.8]{Tyomkin},
    see also \cite[Proposition 6.3]{NishinouSiebert}, \cite[Theorem~4.24]{GrossBook} or \cite[$\S2.3$]{Shustin}. 
\end{proof}

\begin{df}
    We say that $(A_f,f_1,\ldots,f_N)$ is the \emph{tropical limit} of $f$.
\end{df}

\begin{prop}
\label{prop:properties Ci}
    Assume that $f$ as in Proposition \ref{prop:family of curves} defines a curve in $\Tor_\K(\Delta)$ through a generic configuration of $n$ points. Then,
      \begin{enumerate}
      \item the $C_l$ are all rational and have nodes as only singularities. All these nodes lie outside the toric boundaries.
        \item For each triangle $\Delta_l$ in $S_f$, the curve $C_l$ meets each of the three boundary divisors of~$\Delta_l$ in exactly one point, where it is unibranch and smooth.
        \item For each parallelogram $\Delta_l$ in $S_f$, the polynomial $f_l$ defining $C_l$ splits into a product of a monomial and binomials, i.e., it is of the form
        \[x^ry^s(\alpha x^a+\beta y^b)^p(\gamma x^c+\delta y^d)^q\]
        with $\gcd(a,b)=\gcd(c,d)=1$ and $ad\neq cd$.
        \item If $\sigma=\Delta_l\cap\Delta_s$ is a common edge, then the curves $C_l$ and $C_s$ meet in exactly one point~$z\in \Tor_k(\sigma)$ and the intersection multiplicity of $C_l$ and $C_s$ with $\Tor_k(\sigma)$ both equals~$m=\vert\sigma\vert$. This condition is called the \emph{kissing condition}.
    \end{enumerate}
\end{prop}
\begin{proof}
    Shustin does this with a topological argument \cite[p.185 $\S3.3$]{Shustin}. Nishinou-Siebert prove this in \cite[Theorem 8.3]{NishinouSiebert} in an algebraic setting by showing that the curve~$C_0$ consists of images of \emph{maximally degenerate stable maps} which are the curves with the aforementioned properties. 
\end{proof}

\subsection{Refined tropical limit}
\label{subsection:refinedtropicallimit}
We have seen that with a curve $C\subset \Tor_\K(\Delta)$, we can associate a tropical curve which defines a toric degeneration $\pi:\mathcal{X}\rightarrow \A^1_k$ and a family of curves $\widetilde{C}$ in this toric degeneration with generic fiber $C$ and special fiber $C_0$. It turns out (see for example \cite{Shustin, NishinouSiebert,ArguzBousseau}) 
that two things can happen to a node of $C$ in this degeneration process. Either the node degenerates to a node of one of the components $C_l$ of $C_0$ and this node of $C_l$ is not contained in the toric boundary of $\Tor_k(\Delta_l)$; or the node degenerates to a point $z\in C_s\cap C_l$ contained in the toric divisor $\Tor_k(\sigma)$ where $\sigma=\Delta_l\cap\Delta_s$ is the common edge. We deal with the latter case in this subsection in the following way. We perform a non-toric blow up at the point $z$. This does not affect the curve $C$, we have only performed a coordinate change of the coordinates $x$ and $y$ in $f(x,y)$ and get a different defining equation $\widetilde{f}$ for $C$. Repeating what we did in Subsections \ref{subsection:tropicalization}
and \ref{subsection:toricdegeneration} with this new equation $\widetilde{f}$, we get a different tropical curve, toric degeneration and tropical limit. In particular, we get a new triangle $\Delta_{[\sigma]}$ in the dual subdivision of the new tropical curve and a curve $C_{[\sigma]}$ such that the nodes of $C$ that degenerated to $z\in \Tor_k(\sigma)$ before, now degenerate to nodes of $C_{[\sigma]}$. 

We do all of this for two reasons. Firstly, we can find the types of nodes of $C$ that degenerated to $z$ before, by computing the types of nodes of $C_{[\sigma]}$. Secondly, we will see in Section \ref{section:patchworking} that the $C_{[\sigma]}$'s are part of the input of the refined patchworking theorem which we need for finding all curves $C$ that tropicalize to a given tropical curve $A$.

We recall Shustin's computations for the non-toric blow up and the curves $C_{[\sigma]}$.
Assume that $\Delta_l$ and $\Delta_s$ in $S$ have a common edge $\sigma=\Delta_l\cap\Delta_s$. Then $C_l$ and $C_s$ meet at a point~$z\in \Tor_k(\sigma)$. We will see that this point deforms into $\vert\sigma\vert-1$ nodes as illustrated in Figure \ref{fig:defintonodes}. As done in \cite[$\S3.5$]{Shustin}, we perform a blow up at the point $z$ and compute the tropical limit of the blown up curve. 

\begin{figure}
    \begin{tabular}{ccc}
    \includegraphics[scale=0.65]{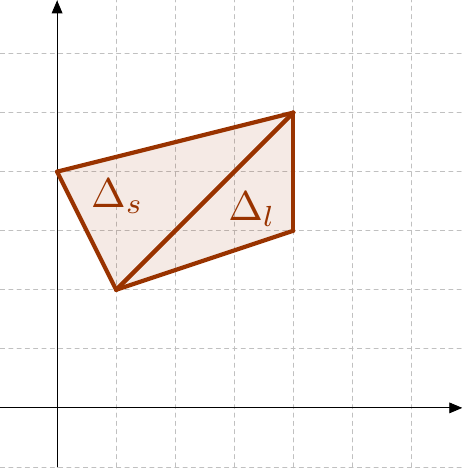}&
    \includegraphics[scale=0.65]{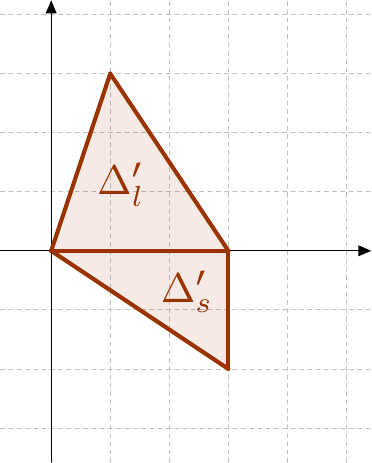}&
    \includegraphics[scale=0.65]{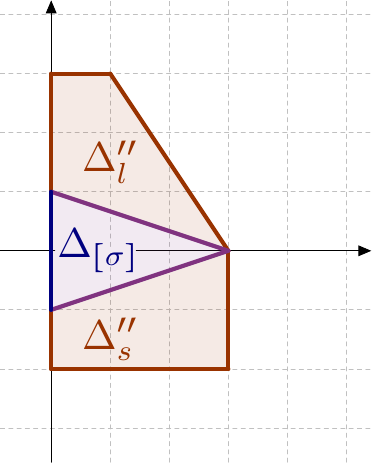}
    \end{tabular}
    \caption{Coordinate change and shift.\label{fig:defpattern}}
    \end{figure}

Assume that $\sigma=\Delta_s\cap \Delta_l$, where $\Delta_s$ and $\Delta_l$ are both triangles in $S_f$ and $\vert \sigma\vert=m\ge 2$. Let $z_{\sigma} \in \Tor_k(\sigma)$ be the point where $C_s$ and $C_l$ meet $\Tor(\sigma)$ with multiplicity $m$. 
We perform a coordinate transformation that sends $\sigma$ to  
$\sigma'\coloneqq[(0,0),(m,0)]$ and such that the transformation of $\Delta_l$ and $\Delta_s$ called $\Delta_l'$ and $\Delta_s'$ lie in the right half plane (see the transform from the first to the second picture in Figure \ref{fig:defpattern}).
We call the new coordinates~$x'$ and $y'$ and the new polynomial $f'(x',y')$.

After multiplying by a constant from $\K^\times$ we can assume that the piecewise linear function~$\nu_f'$ (see \eqref{eq:nu} and Remark \ref{rmk:nu} for the definition) is zero along $\sigma'$. Let $x'=x''+\xi$ and $y'=y''$, where $\xi$ is the $x$-coordinate of $z_\sigma'$ yielding $f''(x'',y'')\coloneqq f'(x',y')$. 
Note that the polygons $\Delta_l''$ and $\Delta_s''$ are not necessarily triangles anymore (see Figure \ref{fig:defpattern} on the right).
Let
\[\Delta_{[\sigma]}\coloneqq \operatorname{Conv}\{(0,1),(0,-1),(m,0)\},\]
that is, $\Delta_{[\sigma]}$ is the blue triangle in Figure \ref{fig:defpattern}. 

There is $\tau(t)\in \K$ with $\val \tau(t) <0$ such that $\widetilde{f}(\widetilde{x},\widetilde{y})\coloneqq f''(\widetilde{x}+\tau(t),\widetilde{y})$ does not contain the monomial $\widetilde{x}^{m-1}$. The subdivision of the Newton polygon $\widetilde{\Delta}$ of $\widetilde{f}$ contains a subdivision of the triangle $\Delta_{[\sigma]}$
and this subdivision has no vertices at the point $(m-1,0)$.
Assume~$\widetilde{f}$ is of the form
\[\sum_{(i,j)\in \widetilde{\Delta}\cap\Z^2}(b_{ij}+d_{ij})t^{\nu_{\widetilde{f}}(i,j)}x^iy^j\]
with $b_{ij}\in k$ and $d_{ij}\in \K$ with $\val(d_{ij})<0$.

\begin{df}
\label{df:defpattern}
    \item We call the polynomial 
    \[\mathcal{T}_\sigma(f)\coloneqq\sum_{(i,j)\in \Delta_{[\sigma]}}b_{ij}x^iy^j\]
    defined by the initials $b_{ij}$ of $\widetilde{f}$ for $(i,j)\in \Delta_{[\sigma]}$, the \emph{$[\sigma]$-refinement of $f$}. 
\end{df}

\begin{figure}
    \begin{tabular}{ccc}
    \includegraphics[scale=0.57]{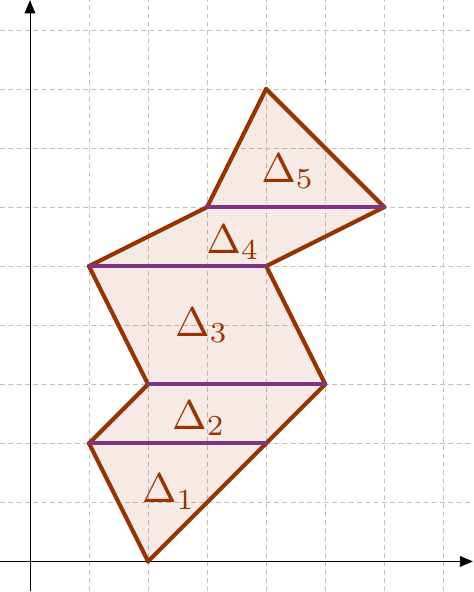}&
    \includegraphics[scale=0.57]{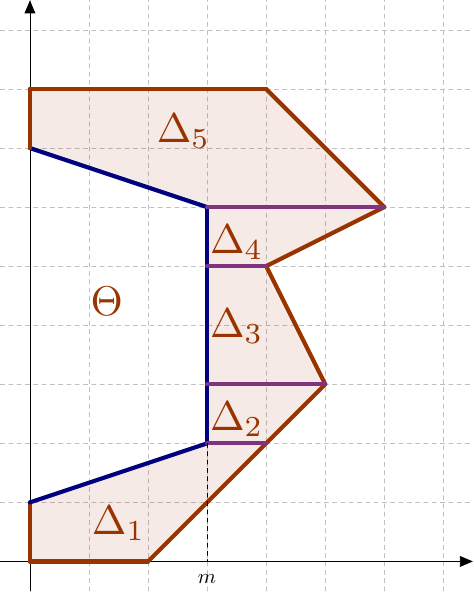}&
    \includegraphics[scale=0.57]{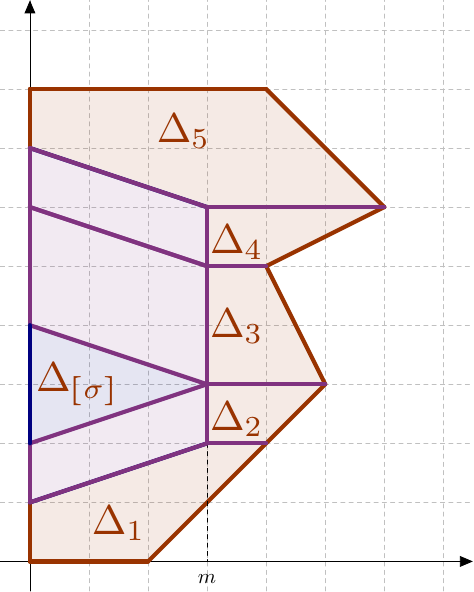}
    \end{tabular}
    \caption{Coordinate change and shift.\label{fig:defpatternwithparallelogram}}
    \end{figure}

We still need to consider the case that $\sigma$ is the edge of a parallelogram. We recall the following from \cite[$\S3.6$]{Shustin}. Assume $\vert \sigma\vert=m\ge 2$ and $\sigma$ is a common edge of~$\Delta_l$ and $\Delta_s$ out of which at least one is a parallelogram. After relabelling there is a chain $\Delta_1,\ldots,\Delta_p$ of polygons in $S_f$, such that $\Delta_1$ and $\Delta_p$ are triangles, $\Delta_2,\ldots,\Delta_{p-1}$ are parallelograms and~$\sigma_1=\Delta_1\cap\Delta_2,\ldots,\sigma_{p-1}=\Delta_{p-1}\cap\Delta_p$ are common edges which are parallel to each other. We can assume that $\nu_f$ is constant along the edges $\sigma_1,\ldots,\sigma_{p-1}$ and that after applying a coordinate change, we can achieve that the polygons all lie in the right half-plane and the edges $\sigma_1,\ldots,\sigma_{p-1}$ are horizontal. Let $f'(x',y')$ be the polynomial after this coordinate change. The truncations of $f'$ of the edges $\sigma_1,\ldots,\sigma_{p-1}$ contain a factor~${(x-\xi)^m}$ for some~$\xi\in k^\times$. Again let $f''(x'',y'')=f'(x'+\xi,y')$. Shustin shows that the polygons~$\Delta_1'',\ldots,\Delta_p''$ bound a trapezoid $\Theta$ and that $\nu_{f''}$ on $\Delta_1''\cup\ldots\cup\Delta_p''$ extends to exactly one subdivision of $\Theta$ consisting of parallelograms and one triangle
\[\Delta_{[\sigma]}=\operatorname{Conv}\{(0,d-1),(0,d+1),(m,d)\}\]
with edges parallel to the edges in $\Theta$ for some $d\in \Z$ (see Figure \ref{fig:defpatternwithparallelogram}). As before there exists~$\tau(t)$ with $\val(\tau)<0$ such that in $\widetilde{f}(\widetilde{x},\widetilde{y})=f''(x''+\tau(t),y'')$ the coefficient of~$\widetilde{x}^{m-1}\widetilde{y}^d$ vanishes.  

\begin{df}
\label{df:defpatternextended edge}
 
   We call the polynomial 
    \[\mathcal{T}_\sigma(f)\coloneqq\sum_{(i,j)\in \Delta_{[\sigma]}}b_{ij}x^iy^j\]
    defined by the initials $b_{ij}$ of $\widetilde{f}$ for $(i,j)\in \Delta_{[\sigma]}$, the \emph{$[\sigma]$-refinement of $f$}. 
\end{df}

The next lemma shows that the $[\sigma]$-refinements $\mathcal{T}_{[\sigma]}f$'s in Definitions \ref{df:defpattern} and~\ref{df:defpatternextended edge} define rational curves with~$m-1$ nodes. 
Analogously to \cite[Lemma 3.10]{Shustin}.

\begin{lm}
The curve defined by $T_{[\sigma]}(f)$ defines a rational curve in $\Tor_k(\Delta_{[\sigma]})$ with exactly~$m-1$ ordinary double points as only singularities.
\end{lm}
\begin{proof}
This follows from Proposition \ref{prop:properties Ci}: The curve defined by $\widetilde{f}$ equals the curve defined by $f$, but it defines a different toric degeneration. Proposition \ref{prop:family of curves} tells us that we get a curve in the special fiber of this new toric degeneration. In particular, the curve defined by $T_{[\sigma]}f$ in $\Tor_k(\Delta_{[\sigma]})$ is a component of this curve. By Proposition \ref{prop:properties Ci} this component is rational and nodal. The number of nodes equals the number of interior points of $\Delta_{[\sigma]}$, which is $m-1$.
\end{proof}

\begin{df}
\label{df:refinedtroplimit}
    We call the data $\left(A_f,f_1,\ldots,f_N,\{\mathcal{T}_{[\sigma]}f:[\sigma]\in E^*(A_f),\vert\sigma\vert\ge 2\}\right)$ \emph{refined tropical limit} of $f$.
\end{df}

Note that the $b_{ij}$ in the definition of $\mathcal{T}_{[\sigma]}(f)$ in both Definitions \ref{df:defpattern} and \ref{df:defpatternextended edge} on the edges $[(0,d+1),(m,d)]$ and $[(0,d-1),(m,d)]$ (where in the first case $d=0$) are determined by $f_l$, $f_s$ and the point $z$.
\begin{df}
\label{df:deformationpattern}
    A polynomial $f_{[\sigma]}$ with Newton polygon $\Delta_{[\sigma]}$ with vanishing coefficient of~$x^{m-1}y^d$ and such that the truncations to the edges $[(0,d+1),(m,d)]$ and $[(0,d-1),(m,d)]$ agree with the respective truncations of $\mathcal{T}_{[\sigma]}f$  is called \emph{deformation pattern} compatible with~$f_l$, $f_s$ and $z$.
\end{df}

\subsection{Condition to pass through a fixed point}
\label{subsection:conditionspoints}
The last information we obtain from the polynomial $f$ defining $C\subset \Tor_\K(\Delta)$ needed for the refined patchworking theorem is called \emph{refined condition to pass through $p_l$} for each $p_1,\ldots,p_n$. In what follows we write $h.o.t.$ for monomials with coefficients with valuations lower than the summands before, that is higher order terms in $t$, and $O(t)$ for monomials with coefficients with negative valuation.
We recall the following from \cite[$\S2.5.9$]{tropicalBook} and \cite[Lemma 3.12]{Shustin}.
As before, we start with
\[f(x,y)=\sum_{(i,j)\in \Delta\cap\Z^2}(a_{ij}+c_{ij}(t))t^{\nu(i,j)}x^iy^j\]
with $a_{ij}\in k$, $c_{ij}(t)\in \K$
such that $\val(c_{ij})<0$, defining a curve in $\Tor_\K(\Delta)$ through $p_1,\ldots,p_n$.
Assume that $\sigma_l$ is the edge in $S$ dual to the edge marked with $x_l=\val(p_l)$, for some~$l$ in~$\{1,\ldots,n\}$.
After a coordinate change and changing $\nu$ up to a constant in $\K^\times$, we can assume that $\sigma_l=[(0,0),(m,0)]$ with $m=\vert\sigma\vert$, that $\nu\vert_{\sigma_l}=0$, and that $\nu(i,j)>0$ for~$(i,j)\not\in \sigma_l$. Then, we can write 
\[p_l=(\xi^0_l+\xi^1_lt,\eta_l^0+\eta^1_lt)\]
for some $\xi^0_l,\eta^0_l\in k^\times$ and $\xi^1_l,\eta^1_l\in \K$ with $\val(\xi^1_l),\val(\eta^1_l)\le 0$.
Then, 
\[f(x,y)=\sum_{i=0}^mc_{i0}(t)x^i+a_{m0}(x-\xi_l^0)^m+O(t).\]
This is because the special curve $C_0$ meets $\Tor_k(\sigma)$ at $(\xi_l^0,\eta_l^0)$ with multiplicity $m$ (see Proposition \ref{prop:properties Ci}).
Setting $x'=x-\xi_l^0$, we obtain
\[f'(x',y)=\sum_{i=0}^{m-1}c'_{i0}(t)x'^i+(a_{m0}+c_{m0}(t))x'^m+O(t)\]
for some $c_{i0}'(t)\in\K$ with $\val(c_{i0}'(t))<0$ for $i=0,\ldots,m-1$. Let
\[\tau(t)\coloneqq -\frac{c'_{m-1,0}(t)}{ma_{m0}}+h.o.t.,\]
such that if $x''=x'-\tau(t)$, then
\[f''(x'',y)\coloneqq f'(x',y)\]
has a zero coefficient at $x''^{m-1}$.
Then, the polynomial equation $f''(x'',y)$ equals the sum
\[\sum_{i=0}^{m-2}c_{i0}(t)x''^i +(a_{m0}+c_{m0})x''^m+ (a_{01}+c_{01}(t))yt^r+(a_{0,-1}+c_{0,-1}(t))y^{-1}t^s+h.o.t.,\]
for some $c_{i0}(t)\in \K$ with $\val(c_{i0}(t))<0$, and some $r,s>0$ with $r\neq s$. Assume that $r<s$.
Now we plug in the new coordinates for $p_l$. Namely,
\[(t\xi_l^1-\tau+\ldots, \eta^0_l+\eta_m^1t+\ldots).\]
By equating the coefficient of the minimal power of $t$ to zero, we get the following equation 
\[\eta_l^0a_{01}t^r+a_{m0}(\xi^1_mt-\tau)^m+h.o.t.=0\]
which yields
\[\tau=\xi_l^1t-\left(-\frac{\eta_l^0a_{01}}{a_{m0}}\right)^{\frac{1}{m}}t^{\frac{r}{m}}+h.o.t.\]
Hence, the coefficients $c_{00}(t)$ and $c_{m0}(t)$ satisfy
\begin{equation}
    \label{eq:conditionlengthbigger1}
    \begin{split}
    &c_{00}(t)-\frac{a_{00}}{a_{m0}}c_{m0}(t)=\\
    &(-1)^mma_{m0}(\xi_l^0)^{m-1}\xi_l^1t+(-1)^{m-1}m(\xi_l^0)^{m-1}(-\eta_l^0a_{01}a_{m0}^{m-1})^{\frac{1}{m}}t^{\frac{r}{m}}+h.o.t.
    \end{split}
\end{equation}
Note that there are $m$ different choices for \eqref{eq:conditionlengthbigger1} corresponding to the $m$-th roots of the coefficient~$-\eta_s^0a_{01}a_{m0}^{m-1}$.
We call the choice of $m$-th root of $-\eta_s^0a_{01}a_{m0}^{m-1}$ defined by the curve $C$, a \emph{refined condition to pass through $p_l$}.
\section{From tropical to algebraic curves}
\label{section:patchworking}

Let $\mathcal{P}=\{p_1,\ldots,p_n\}$ be a generic configuration of $\K$-points in $\Tor_\K(\Delta)$ such that the points $x_1,\ldots,x_n$ with~$x_i=\val(p_i)$ for $i=1,\ldots, n$ is $\Delta$-generic. 
In Section \ref{section:tropicalization} we have seen that a curve $C\subset\Tor_\K(\Delta)$ defined by a (Laurent) polynomial~$f$ passing through $\mathcal{P}$ defines
\begin{itemize}
    \item a simple tropical curve $A_f$ with Newton polygon $\Delta$ through $x_1,\ldots,x_n$, with a dual subdivision $S_f:\Delta=\Delta_1\cup\ldots\cup\Delta_N$;
    \item rational curves $C_l\subset\Tor_k(\Delta_l)$, for $l=1,\ldots,N$, such that
    \begin{itemize}
        \item for a triangle $\Delta_l$, the curve $C_l$ is irreducible, rational, nodal and meets $\Tor(\partial \Delta_s)$ at exactly $3$ points where it is unibranch and smooth,
        \item for a parallelogram $\Delta_l$, the defining polynomial $f_l$ of $C_l$ splits into a product of a monomial and binomials;
    \end{itemize}
    \item $[\sigma]$-refinements $\mathcal{T}_{[\sigma]}(f)$ for each extended edge $[\sigma]$ with $\vert\sigma\vert\ge 2$ defining a rational curve $C_{[\sigma]}\subset \Tor_k(\Delta_{[\sigma]})$ with $\vert\sigma\vert-1$ nodes as only singularities;
    \item refined conditions to pass through the points $p_1,\ldots,p_n$.
\end{itemize}

In the following theorem we, use the non-quadratically enriched classical correspondence theorems to argue that the data above uniquely determines the curve $C\subset \Tor_\K(\Delta)$ one started with, and to quadratically enrich Shustin's \emph{refined patchworking theorem}\linebreak\cite[$\S5.1$]{Shustin}, \cite[Theorem 2.51]{tropicalBook}.  

\begin{thm}[Quadratically enriched refined patchworking theorem]
\label{thm:patchworking}
Let $k$ be a field of characteristic zero or greater than $\operatorname{diam}(\Delta)$.
    Let $p_1,\ldots,p_n\in\Tor_\K(\Delta)$ be a generic configuration of~$\K$-points that projects by the valuation $\val$ to a $\Delta$-generic configuration~$x_1,\ldots,x_n\in \Q^2$. Assume further that we are given
    \begin{enumerate}
        \item a simple irreducible genus $g$ tropical curve $A$ with Newton polygon $\Delta$ passing through~$x_1,\ldots,x_n$, with dual subdivision $S:\Delta=\Delta_1\cup\ldots\cup \Delta_N$ induced by a piecewise linear function $\nu$,
        \item for each $\Delta_l$ in $S$, a rational curve $C_l$
        defined by $\sum_{(i,j)\in \Delta_l\cap\Z^2}a_{ij}x^iy^j$, with coefficients $a_{ij}\in k$ with the properties (1)-(4) of Proposition \ref{prop:properties Ci},
        \item deformation patterns $C_{[\sigma]}$ defined by $f_{[\sigma]}$ with coefficients in $k$ compatible with the curves $C_1,\ldots, C_N$ in the sense of Definition \ref{df:deformationpattern}.
        \item refined conditions to pass through the points $p_1,\ldots,p_n$ as in \eqref{eq:conditionlengthbigger1} defined over $k$.
    \end{enumerate}
    Then there exists a unique irreducible genus $g$ nodal curve $C\in \Tor_{k\Puiseux}(\Delta)$ defined over~$k\Puiseux$ passing through~$p_1,\ldots,p_n$, defined by a polynomial of the form
    \[f(x,y)=\sum_{(i,j)\in \Delta\cap\Z^2}(a_{ij}+c_{ij}(t))t^{\nu(i,j)}x^iy^j\]
    such that $A_f=A$, the tropical limit of $C$ agrees with $C_l$ for $l=1,\ldots,N$, $\mathcal{T}_{[\sigma]}(f)=f_{[\sigma]}$ for any extended edge $[\sigma]$ with $\vert \sigma\vert\ge 2$, and such that the refined conditions to pass through the points $p_1,\ldots,p_n$ for this curve agree with the given ones.

    Furthermore, there is a bijection
    \[\phi\colon \{\text{nodes $z$ of $C$}\}\xrightarrow{1:1}\{\text{nodes of $C_l$ for $l=1,\ldots,N$ and of $C_{[\sigma]}$ for $[\sigma]\in E^*(A)$}\},\]
    such that for a node $z$ of $C$,
    the residue field $\kappa(z)$ equals $\kappa(\phi(z))\Puiseux$ where $\kappa(\phi(z))$ is the residue field of the node $\phi(z)$, 
    and such that $\operatorname{Type}(z)\mapsto \operatorname{Type}(\phi(z))$ under the bijection~$\kappa(z)\Puiseux^{\times}/(\kappa(z)\Puiseux^\times)^2\cong \kappa(z)^\times /(\kappa(z)^\times)^2$ from Remark \ref{rmk:GW(K)}.
 
\end{thm}

\begin{proof}
By the correspondence theorems in \cite{Mikhalkin} and \cite{Shustin} for $k=\C$ or $k=\R$, \cite{NishinouSiebert} for $\operatorname{char} k =0$ or \cite{Tyomkin} for $\operatorname{char} k =0$ or $\operatorname{char} k>\operatorname{diam}\Delta$, there exist exactly $\operatorname{mult}_\C(A)$ curves $C\subset \Tor_{\overline{k}\Puiseux}(\Delta)$ through $p_1,\ldots,p_n$ defined over $\overline{k}\Puiseux$ and tropicalizing to $A$, where $\overline{k}$ is an algebraic closure of $k$. We will see in Section \ref{section:correspondencetheorem} that the number of different choices for the input data in this theorem equals $\operatorname{mult}_\C(A)$. Each of the $\operatorname{mult}_\C(A)$ curves gives rise to a unique choice for the input data. Thus there exists exactly one curve~$C\subset \Tor_{\overline{k}\Puiseux}(\Delta)$ matching the given input data. Indeed this curve is defined over $\K=k\Puiseux$ since there is a unique solution.

It follows also from the non-quadratically enriched correspondence theorems that the nodes of $C$ correspond bijectively to the nodes of the $C_l$ for $l=1,\ldots,N$ and the $C_{[\sigma]}$ for $[\sigma]\in E^*(A)$ (see for example \cite[Theorem 5]{Shustin} or \cite{ArguzBousseau}).
It remains to show the statement about the types of the nodes. We show this in a very similar fashion to the computations in \cite[$\S6$]{Pauli}.

Let $z$ be a node of $C$. Then $\val(z)\in A$ either lies on an edge of $A$ or is a vertex.
We start by assuming the latter. Let $\Delta_l$ be the polygon dual to this vertex in the dual subdivision~$S$. Then $\phi(z)$ is a node of $C_l$.
After a coordinate change and after multiplying~$\nu$ by a unit in~$k\Puiseux$ we can assume that
$\nu$ equals $0$ on $\Delta_l$ and is positive outside of $\Delta_l$.
Hence, we can write $z$ as
 $z=(\alpha,\beta)=(\alpha^0+\alpha^1(t),\beta^0+\beta^1(t))$ with $\alpha^0,\beta^0\in F$ and~$\alpha^1(t),\beta^1(t)\in F\Puiseux$, with $\val(\alpha^1(t)),\val(\beta^1(t))<0$, where $F$ is some finite separable field extension $F/k$. Then,
 \[f(x,y)=\sum_{(i,j)\in \Delta_l\cap\Z^2}a_{ij}x^iy^j+h.o.t.\]
 Here, we again mean higher order terms in $t$.
 Hence, $(\alpha_0,\beta_0)$ is a singular point of $C_l$ and thus a node (we assume that the $C_l$ are all nodal). One can find $(\alpha_0,\beta_0)$ by solving 
 \[f_l(x,y)=0\text{, }\frac{\partial f_l}{\partial x}(x,y)=0\text{, }\frac{\partial f_l}{\partial y}(x,y)=0.\]
Since $(\alpha^0,\beta^0)$ is a node of $C_l$, its Hessian determinant $\det \operatorname{Hessian}(f_l)(\alpha^0,\beta^0)$ does not vanish. It follows that we get linear equations which one can use to solve for the initial coefficients in $\alpha^1 $ and $\beta^1$, and then also inductively for all other coefficients. It follows that~$\alpha^1,\beta^1$ are elements of $F\Puiseux$ (and not in some field extension of $F\Puiseux$).
We compute the type of $z$:
\begin{align*}
    \operatorname{Type}(z)=&-\det \operatorname{Hessian}(f)(z)\\
    =&-\det \operatorname{Hessian}(f_l)(z)+h.o.t.\\
    =&-\det \operatorname{Hessian}(f_l)(\alpha^0,\beta^0)+h.o.t.
\end{align*}
in $F\Puiseux^\times/(F\Puiseux^\times)^2$. Since $(\alpha^0,\beta^0)$ is a node of $C_l$, $-\det \operatorname{Hessian}(f_l)(\alpha^0,\beta^0)$ does not vanish and equals $\operatorname{Type}(z)$ in $F\Puiseux^\times/(F\Puiseux^\times)^2$ (see Remark \ref{rmk:GW(K)}). 

The type of~$\phi(z)=(\alpha^0,\beta^0)$ is $-\det \operatorname{Hessian}(f_l)(\alpha^0,\beta^0)\in F^\times/(F^\times)^2$ and thus one can identify the type of $z$ with the type of $\phi(z)$ via the bijection between $F^\times/(F^\times)^2$ and~$F\Puiseux^\times/(F\Puiseux^\times)^2$ as in Remark \ref{rmk:GW(K)}.

Now assume $z$ is a node of $C$ such that $\val(z)$ is not a vertex, but lies on an edge $\sigma$. In this case $\phi(z)$ is a node of $C_{[\sigma]}$. Let $\sigma$ be the edge of $S$ dual to this edge. After the coordinate change from Subsection \ref{subsection:refinedtropicallimit} and changing $\nu$ such that $\nu$ restricted to $\Delta_{[\sigma]}$ is zero and positive outside of $\Delta_{[\sigma]}$, we can assume that $z$ is of the form 
\[z=(\alpha^0+\alpha^1,\beta^0+\beta^1),\]
with $(\alpha^0,\beta^0)$ defined over some finite separable field extension $F/k$, and such that $\alpha^1,\beta^1$ are in~$F\Puiseux$, with $\val(\alpha^1),\val(\beta^1)<0$. All these coordinate changes, change the determinant of the Hessian by a square and the same computation as above proves the claim about the types.
 \end{proof}

\section{Quadratic enrichments}
\label{section:quadraticenrichments}

    Theorem \ref{thm:patchworking} provides a way to reconstruct a unique algebraic genus $g$ curve passing through $p_1,\ldots,p_n$ in $\Tor_\K(\Delta)$ from the data (1)-(4). In this section we find the different possibilities for the curves $C_l$ in (2) and $C_{[\sigma]}$ in (3) in Theorem \ref{thm:patchworking} explicitly, analogous to the methods in~\cite{Shustin}. Subsequently, we use these expressions to compute the quadratic weights of the curves  $C_l$ and $C_{[\sigma]}$. By the same theorem, these determine the quadratic weight of the unique curve one reconstructs from this data.

\subsection{$\Wel_{\kappa(C_l)}^{\A^1}(C_l)$ when $\Delta_l$ is a triangle}
\label{subsection:lemma3.5}
In this subsection we quadratically enrich \cite[Lemma 3.5]{Shustin}. 
Let $\Delta'\coloneqq \Delta_l$ be a lattice triangle with edges $\sigma$, $\sigma'$ and $\sigma''$.
Assume that for all $(i,j)\in (\sigma\cup \sigma')\cap\Z^2$, there are given coefficients $a_{ij}\in K'$, in a finite field extension $K'$ of our base field $k$.
We want to find all rational curves $C'$ in $\Tor_{K'} (\Delta')$ meeting each boundary divisor at exactly one $K'$-point and being defined by a polynomial
\[f'\coloneqq \sum_{(i,j)\in \Delta'\cap \Z^2}a_{ij}x^iy^j=\sum_{(i,j)\in \left(\Delta'\setminus (\sigma\cup\sigma')\right)\cap\Z^2}a_{ij}x^iy^j+\sum_{(i,j)\in (\sigma\cup \sigma')\cap\Z^2}a_{ij}x^iy^j.\]
That is, we have to find all possible coefficients  $a_{ij}$, with $(i,j)\in \left(\Delta'\setminus (\sigma\cup\sigma')\right)\cap\Z^2$, such that $f'$ defines such a curve.

\begin{df}
    A \emph{lattice preserving transformation} is an affine transformation 
    $z\mapsto Az+b$,
    where $A\in \operatorname{GL}_2(\Z)$ and $b\in \Z^2$.
\end{df}

\begin{rmk}
    A lattice preserving transformation preserves the number of interior points, boundary points and area of a lattice polygon. 
\end{rmk}
We perform the following lattice preserving transformation to our triangle $\Delta'$.
\begin{lm}
\label{lm:latticepathtrafo}
    Let $\Delta'$ be a lattice triangle and let $\sigma''$ be one of its edges. Then, there exists a lattice preserving transformation that sends $\Delta'$ to $\operatorname{Conv}\{(0,m),(p,0),(q,0)\}$, for some integers $0\le p<q\le m$, and such that the edge $\sigma''$ is sent to the edge $[(q,0), (0,m)]$.
\end{lm}

\begin{figure}
    \centering
    \begin{tabular}{ccc}
    \includegraphics{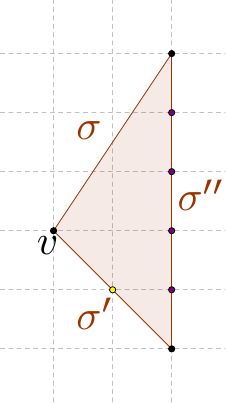}&
    \includegraphics[width=0.3\textwidth]{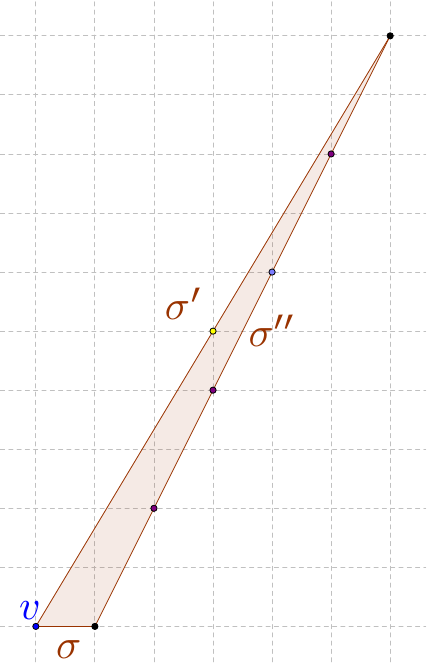}&
    \includegraphics[width=0.267\textwidth]{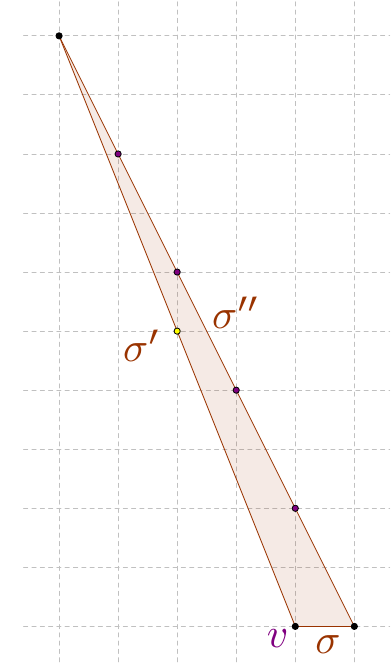}
    \end{tabular}
    \caption{Lattice preserving transformations applied to an arbitrary triangle.}
    \label{fig:latttransf}
\end{figure}

\begin{proof}
    Call $\sigma$ and $\sigma'$ the remaining edges of $\Delta'$  and assume that $\vert \sigma\vert \le \vert \sigma'\vert$. Let $v\in \Z^2$ be the common vertex of $\sigma$ and $\sigma'$. Up to a translation by $-v\in \Z^2$, we can assume that the vertex $v$ is at the origin. Let $(a, b)^{\mkern-1.5mu\mathsf{T}}\in\Z^2$
    be the primitive vector pointing into the direction of $\sigma$. Since $\gcd(a,b)=1$, there exist $n,m\in \Z$ such that $an+bm=1$.
    We perform the lattice preserving transformation defined by 
    $$\begin{pmatrix}n & m \\b & - a \end{pmatrix}\in \operatorname{GL}_2(\Z),$$ which sends the edge~$\sigma$ to the~$x$-axis. Let 
    $w=(w_1,w_2)^{\mkern-1.5mu\mathsf{T}}\in\Z^2$
    be the primitive vector pointing into the direction of the transformed edge of $\sigma'$. Up to the reflection across the $x$-axis, we can assume that $w_2>0$.
    Let $l\in\Z$ be the minimum integer such that~$w_1-kw_2$ is non-positive. We apply the transformation $$\begin{pmatrix}
        1 & -l\\ 0 & 1
    \end{pmatrix}\in \operatorname{GL}_2(\Z),$$
    which leaves the transformed edge $\sigma$ invariant, since it lies on the $x$-axis, and moves the image of $\sigma'$ to 
    $\vert \sigma'\vert(w_1-lw_2, w_2)^{\mkern-1.5mu\mathsf{T}}$.
    Next, we shift the triangle by $\vert \sigma'\vert(lw_2-w_1)$ to the right.
    We have that $p\coloneqq\vert \sigma'\vert(lw_2-w_1)\ge0$, $q\coloneqq p+\vert \sigma\vert> p$ and~$m\coloneqq\vert \sigma'\vert w_2\ge q$ since otherwise $\vert \sigma\vert>\vert \sigma'\vert(w_1-(l-1)w_2)\ge\vert \sigma'\vert$. See Figure~\ref{fig:latttransf} for an example of this algorithm.
\end{proof}

\begin{rmk}
    \label{rmk:doubleareaandlatticelengths}
    For the triangle $\Delta'=\operatorname{Conv}\{(0,m),(p,0),(q,0)\}$ as in the Lemma \ref{lm:latticepathtrafo} above we have that its double Euclidean area equals
    \[\vert \Delta'\vert=m\cdot(q-p),\]
    and the lattice lengths of the edges equal
    \begin{align*}\vert\sigma\vert&=\vert[(p,0),(q,0)]\vert=(q-p),\\
    \vert\sigma'\vert&=\vert [(0,m),(p,0)]\vert=\gcd(p,m)\eqqcolon d_p,\\
    \vert\sigma''\vert&=\vert [(0,m),(q,0)]\vert=\gcd(q,m)\eqqcolon d_q .\end{align*}
\end{rmk}

A lattice preserving transformation defines a coordinate change of the coordinates $x$ and $y$ of a curve in $\Tor_{K'}(\Delta')$. By \cite[$\S2$]{LevineWelschinger} a coordinate change does not change the type of a node.
Hence, we can assume that $\Delta'=\operatorname{Conv}\{(0,m),(p,0),(q,0)\}$. 
Up to a reparametrization of the domain, we can assume that the three $K'$-points where $C'$ meets the toric boundary are $0$, $1$ and $\infty$. By choosing suitable affine coordinates on the torus, the parametrization of the curve~$C'$ is given by
\begin{align}
\label{eq:parametrization}
    \begin{split}
    x(t)&=\alpha t^m,\\
    y(t)&=\beta t^p(t-1)^{q-p}.
    \end{split}
\end{align}
Then finding all possibilities for the $a_{ij}$ with $(i,j)\in \left(\Delta'\setminus(\sigma\cup\sigma')\right)\cap \Z^2$ amounts to finding all possibilities for $\alpha$ and $\beta$ in \eqref{eq:parametrization}.

\begin{lm}
\label{lm:F'}
    The different choices for $\alpha$ and $\beta$ in \eqref{eq:parametrization} correspond to the solutions of 
    \begin{equation}
    \label{eq:equationsalphabeta}
    a\alpha+\epsilon_1=0\text{, }\quad b\beta^{\frac{m}{d_p}}(-1)^{\frac{m(q-p)}{d_p}}+\epsilon_2\alpha^{\frac{p}{d_p}}=0
    \end{equation}
    for some fixed $a,b,\epsilon_1,\epsilon_2\in K'^\times$.
    In particular, the $K'$-algebra defining these different choices is given by (note that $\alpha$ already lies in $K'$)
    \begin{equation}
    \label{eq:F'}
    F'=\frac{K'[\beta]}{\left(b\beta^{\frac{m}{d_p}}(-1)^{\frac{m(q-p)}{d_p}}+\epsilon_2(\frac{-\epsilon_1}{a})^{\frac{p}{d_p}}\right).}
    \end{equation}
\end{lm}
\begin{proof}
    This follows the ideas of the proof of \cite[Lemma 3.5]{Shustin}. We generalize them so that they work over an arbitrary base field. 
    The curve $C'$ has defining equation
    \[f'=\sum_{(i,j)\in \Delta'\cap\Z^2}a_{ij}x^iy^j\]
    with given coefficients $a_{ij}$ for $(i,j)\in \sigma\cup\sigma'=[(p,0),(q,0)]\cup [(0,m),(p,0)]$. Without loss of generality, we can assume that $a_{p0}=1$. 
    By Proposition \ref{prop:properties Ci}, the curve $C'$ should meet each $\Tor_{K'}(\sigma)$ and $\Tor_{K'}(\sigma')$ in one point, with multiplicities $\vert\sigma\vert$ and $\vert \sigma'\vert$, respectively. Hence, the trunctions $f'^{\sigma}$ and $f'^{\sigma'}$ of the defining polynomial $f'$ to $\sigma$ and $\sigma'$ equal
        \[f'^{\sigma}=x^p(ax+\epsilon_1)^{q-p} \text{ and } f'^{\sigma'}=(by^{\frac{m}{d_p}}+\epsilon_2x^{\frac{p}{d_p}})^{d_p},\]
    respectively, for some $a,b,\epsilon_1,\epsilon_2\in K'^\times$ satisfying $\epsilon_1^{q-p}=\epsilon_2^{d_p}=1$.
    In the given parametrization, the curve $C'$ meets $\Tor_{K'}(\sigma)$ at $t=1$. Hence, 
    \[f'^{\sigma}(x(1),y(1))=\alpha^p(a\alpha+\epsilon_1)^{q-p}\]
    has to vanish with order $q-p=\vert\sigma\vert$. This only happens if $a\alpha+\epsilon_1=0$, which gives the first equation in \eqref{eq:equationsalphabeta}. 
    Observe that
    \[f'^{\sigma'}(x(t),y(t))=t^{mp}
    \left(b\beta^{\frac{m}{d_p}}(t-1)^{\frac{(q-p)m}{d_p}}+\epsilon_2\alpha^{\frac{p}{d_p}}\right)^{d_p}.\]
    The curve $C'$ meets the toric divisor $\Tor_{K'}(\sigma')$ at $t=0$ with vanishing order $d_p=\vert\sigma'\vert$ and therefore we get that $b\beta^{\frac{p}{d_p}}(-1)^{\frac{(q-p)m}{d_p}}+\epsilon_2\alpha^{\frac{p}{d_p}}=0$, which is the second equation in \eqref{eq:equationsalphabeta}.
\end{proof}

Next we compute the quadratic weight $\Wel^{\A^1}_{F'}(C')$ of the curve $C'$ defined by the parametrization $\eqref{eq:parametrization}$, with $\alpha$ and $\beta$ as in Lemma \ref{lm:F'} in $\GW(F')$.

\begin{lm}
\label{lm: Welschinger for 3.5}
All curves parametrized by \eqref{eq:parametrization} with $\alpha$ and $\beta$ satisfying \eqref{eq:equationsalphabeta} are nodal and meet the toric boundary $\Tor_{K'}(\partial \Delta')$ where it is unibranch and smooth.
Furthermore,
    \[\Wel^{\A^1}_{F'}(C')=\qinv{ (-1)^{\operatorname{Int}(\Delta')} m^{d_q+d_p+q-p}(q-p)^{d_p+d_q}d_p^{d_p}d_q^{d_q}}\in\GW(F'),\]
    where $\operatorname{Int}(\Delta')$ is the number of interior points in $\Delta'$.
\end{lm}
\begin{proof}
    In order to compute the quadratic weight of the parameterized rational curve~$C'$, we compute the product over all nodes of the curve of the negative determinant of the Hessian of a defining polynomial at every nodal point, up to squares. In order to compute the local contribution at every nodal point, we use the fact that this determinant can be computed up to squares directly from the parametrization.

Namely, let us assume that there are parameters~$t$ and~$s$ defining a non-degenerate double nodal point~$p$. The linear change of variable given by the tangent vectors $(x'(t),y'(t))$ and $ (x'(s),y'(s))$ of the parametrization at $p$, allows us to have an explicit local equation of the form~$\gamma u v+\mathrm{h.o.t.}$, where $\gamma\in k$ and $u,v$ are local coordinates. Therefore, minus the determinant of the Hessian of a defining equation at the node $p$ equals
\[\left(x'(t)y'(s)-x'(s)y'(t)\right)^{2},\]
up to squares in~$k$.
Hence, using the parametrization in \eqref{eq:parametrization}, we have that the local contribution at a node $(x(t),y(t))=(x(s),y(s))$ equals
\begin{equation}\label{eq:hesspar}
    \left[
\frac{mx(t)}{t} \frac{y(s)(qs-p)}{s(s-1)}-\frac{mx(s)}{s} \frac{y(t)(qt-p)}{t(t-1)}
\right]^2=\left(m(q-p)x(t)y(t)\right)^2\left[\frac{t-s}{t(t-1)s(s-1)}
\right]^2.
\end{equation}
Now, let us compute the coordinates of the nodes of~$C'$. The parameters $s$ and $t$ must satisfy the equations $(x(t),y(t))=(x(s),y(s))$ with $t\neq s$. Hence, we have that $s=t\mu$, where $\mu$ is a non-trivial $m$-th root of unity, and 
\[\frac{1-t}{1-\mu t}=\mu^{\frac{p}{q-p}}\nu,\]
where $\nu$ is a $(q-p)$-th root of unity. Put $\zeta=\mu^{\frac{p}{q-p}}\nu$.
Then, the parameters where there is a node are the solutions to the equations
\begin{equation}\label{eq:nodespar}
    t=\frac{1-\zeta}{1-\zeta\mu},\quad s=t\mu,\quad \zeta=\mu^{\frac{p}{q-p}}\nu,
\end{equation}
where $\mu^m=\nu^{q-p}=1$ and $\mu\neq1,\zeta\neq1,\zeta\mu\neq1$. Hence, the solutions are given by 
\begin{align} \label{eq:mupar}
    \mu=&\,\zeta_m^{i}     & i=1,2,\dots, m-1,\\ \label{eq:zetapar}
    \zeta=&\,\zeta_{m(q-p)}^{ip+jm} & ip+jm\not\equiv 0 \mod m(q-p),\\ \label{eq:zetamupar}
    \zeta\mu=&\,\zeta_{m(q-p)}^{iq+jm} &iq+jm\not\equiv 0 \mod m(q-p),
\end{align}
where $j=1,2,\dots, q-p$, and $\zeta_m$ and $\zeta_{m(q-p)}$ are primitive $m$-th and $m(q-p)$-th roots of unity, respectively. However, this counts twice the number of nodes since the branches of a node are unordered. 

Set $m_p=\frac{m}{d_p}$ and $m_q=\frac{m}{d_q}$. Let us remark that $ip+jm\equiv 0 \mod m(q-p)$ if and only if ${i\frac{p}{d_p}+jm_p\equiv 0} \mod m_p(q-p)$ and that if $i\frac{p}{d_p}+jm_p\equiv 0 \mod m_p(q-p)$, then~$i\frac{p}{d_p}\equiv 0 \mod m_p$. Since $\gcd(\frac{p}{d_p},m_p)=1$, we have that $i\equiv 0 \mod m_p$. So there are exactly~$d_p-1$ solutions to~$ip+jm\equiv 0 \mod m(q-p)$. Namely when $i=m_pi'$, where~$i'=1,2,\dots,d_p-1$ and $j$ is the only solution to $j\equiv-i'\frac{p}{d_p}$ in range.
Analogously, there are $d_q-1$ solutions to~$iq+jm\equiv 0 \mod m(q-p)$, given by $i=m_qi'$, where~$i'=1,2,\dots,d_q-1$ and $j$ is the only solution to $j\equiv-i'\frac{q}{d_q}$ in range. Thus, there are $\frac{(m-1)(q-p)-(d_p+d_q-2)}{2}$ distinct nodes. By Pick's theorem, this corresponds to~$\operatorname{Int}(\Delta')$ the number of interior lattice points in the triangle $\Delta'$.

Finally, we can compute the quadratic weight of the curve~$C'$ by multiplying  over the nodes of~$C'$ the expression in~\eqref{eq:hesspar}, or equivalently, by multiplying the square root of this expression over all solutions of~\eqref{eq:nodespar}, multiplying it by $-1$ for each node due to the fact that we interchange the roles of~$s$ and~$t$ in one of the instances of~\eqref{eq:nodespar} that describe the same node. Thus, the quadratic weight is the class of the quadratic form given by
\begin{equation}\label{eq:welnodes}
\displaystyle\prod_{\text{nodes of }C}\left[\frac{t-s}{t(t-1)s(s-1)}\right]^2=
(-1)^{\operatorname{Int}(\Delta')}\prod_{\mathrm{sols. to\eqref{eq:nodespar}}} \frac{t-s}{t(t-1)s(s-1)},
\end{equation}
since the product over the nodes of $C'$ of the factor $(x(t)y(t))^2$ is a square in~$k$.
For a particular solution of~\eqref{eq:nodespar}, the contribution from~\eqref{eq:hesspar} equals
\begin{equation}\label{eq:localcontpar}
    \frac{t-s}{t(t-1)s(s-1)}=\frac{(1-\zeta\mu)^3}{\zeta\mu(1-\zeta)(1-\mu)}.
\end{equation}
Let us compute the product over all solutions of~\eqref{eq:nodespar} of each factor in the right-hand side of~\eqref{eq:localcontpar}. To do that, we will multiply over all instances of the range $i=1,2,\dots,m-1$, $j=1,2,\dots,q-p$ and we will divide the exceeding terms that do not contribute according to the consideration given by each case.

For the factor $1-\zeta\mu$ in the numerator, let us observe that only multiples of $d_q$ can be obtained as integer combinations $iq+jm$, and each multiple $ld_q$ can be obtained $d_q$ times if $ld_q\not\equiv 0 \mod m$ and can be obtained $d_q-1$ times if $ld_q\equiv 0 \mod m$. After dividing by~$d_q$, the product runs for $l=1,2,\dots, m_q(q-p)-1$ and we get that
\begin{equation}\label{eq:factonezetamu}
    \prod_{\mathrm{sols. to\eqref{eq:nodespar}}}(1-\zeta\mu)=
    \frac{\left[\displaystyle\prod_{l=1}^{m_q(q-p)-1}(1-\zeta_{m_q(q-p)}^{l})\right]^{d_q}}
    {\displaystyle\prod_{i=1}^{d_p-1}\left(1-\zeta_{m(q-p)}^{im_p(q-p)}\right) \cdot \prod_{j=1}^{q-p-1}\left(1-\zeta_{m_q(q-p)}^{jm_q}\right)}
    =\frac{\left[\displaystyle m_q(q-p)\right]^{d_q}}
    {\displaystyle d_p \cdot (q-p)},
\end{equation}
where $\zeta_{m_q(q-p)}$ is a primitive $m_q(q-p)$-th root of unity. The numerator has an excess of terms corrected by the second factor of the denominator corresponding to the multiples that appear only $d_p-1$ instead of $d_p$, and the first factor of the denominator corresponds to elements $i,j$ such that $\zeta=1$ which do not contribute to our product. One can check that there are $d_q(m_q(q-p)-1)$ factors in the numerator and $[(q-p)-1+d_p-1]$ factors in the denominator, which account for twice the number of nodes of~$C'$.
A direct application of Lemma~\ref{lm:sineIdentity} yields the second equality in~\eqref{eq:factonezetamu}.

For the factor $\zeta\mu$ in the denominator we can use the same set of indices for the products in the equation. We get that
\begin{equation}\label{eq:factzetamu}
    \prod_{\mathrm{sols. to\eqref{eq:nodespar}}}(\zeta\mu)=
    \frac{\left[\displaystyle\prod_{l=1}^{m_q(q-p)-1}(\zeta_{m_q(q-p)}^{l})\right]^{d_q}}
    {\displaystyle\prod_{i=1}^{d_p-1}\left(\zeta_{m(q-p)}^{im_p(q-p)}\right) \cdot \prod_{j=1}^{q-p-1}\left(\zeta_{m_q(q-p)}^{jm_q}\right)}
    =\displaystyle \zeta_{m(q-p)}^{m(q-p)\left[\frac{m(q-p)-d_p-d_q-(q-p)+2}{2}\right]}=1,
\end{equation}
since~$m(q-p)-d_p-d_q-(q-p)+2\equiv0\mod2$ for all parities of~$m, p$ and~$q$.

For the factor $1-\zeta$ in the denominator the computation is analogous to the one of $1-\zeta\mu$ by switching the roles of~$p$ and~$q$ due to the nature of the description in~\eqref{eq:zetapar} and~\eqref{eq:zetamupar}. Thus, we get that
\begin{equation}\label{eq:factonezeta}
    \prod_{\mathrm{sols. to\eqref{eq:nodespar}}}(1-\zeta)=
    \frac{\left[\displaystyle\prod_{l=1}^{m_p(q-p)-1}\left(1-\zeta_{m_p(q-p)}^{l}\right)\right]^{d_p}}
    {\displaystyle\prod_{i=1}^{d_q-1}\left(1-\zeta_{m(q-p)}^{im_q(q-p)}\right) \cdot \prod_{j=1}^{q-p-1}\left(1-\zeta_{m_p(q-p)}^{jm_p}\right)}
    =\frac{\left[\displaystyle m_p(q-p)\right]^{d_p}}
    {\displaystyle d_q \cdot (q-p)}.
\end{equation}
Lastly, for the factor $1-\mu$, for every instance of $i=1,2,\dots,m-1$, there are $q-p$ solutions for $j$, and we just need to correct for the solutions to $\zeta=1$ and $\zeta\mu=1$. Namely,
\begin{equation}\label{eq:factmu}
    \prod_{\mathrm{sols. to\eqref{eq:nodespar}}}(1-\mu)=
    \frac{\left[\displaystyle\prod_{i=1}^{m-1}\left(1-\zeta_{m}^{i}\right)\right]^{q-p}}
    {\displaystyle\prod_{i=1}^{d_p-1}\left(1-\zeta_{m}^{i}\right) \cdot \prod_{i=1}^{d_q-1}\left(1-\zeta_{m}^{i}\right)}
    =\frac{\displaystyle m^{q-p}}{d_pd_q},
\end{equation}
by applying Lemma~\ref{lm:sineIdentity}.
The theorem follows from equations~\eqref{eq:welnodes} to~\eqref{eq:factmu} by considering terms up to squares in~$k$.
\end{proof}
Let us summarize and reformulate what we did in this subsection: Let $\Delta'$ be a lattice triangle with edges $\sigma$, $\sigma'$ and $\sigma''$.
\begin{enumerate}
    \item We have found all the (finitely many) possibilities for the rational curve $C'$ in $ \Tor_{K'}(\Delta')$ with prescribed coefficients in $a_{ij}\in K'$ for all $(i,j)\in (\sigma\cup \sigma')\cap\Z^2$ such that
       $C'$  meets each boundary divisor of $\Tor_{K'}(\partial \Delta')$ in precisely one $K'$-point with the correct multiplicity. All these curves have only nodes as singularities and meet the toric boundary at non-singular, unibranch points.
    \item The $K'$-algebra defined by the different choices of $C'$ is $F'$ defined in Lemma \ref{lm:F'}.
    \item It follows from Lemma \ref{lm: Welschinger for 3.5} and Remark \ref{rmk:doubleareaandlatticelengths} that 
    the quadratic weight of such a curve $C'$ equals
    \begin{equation}
        \label{eq:quadraticweightofC'triangle}
        \Wel^{\A^1}_{F'}(C')=\qinv{ (-1)^{\operatorname{Int}(\Delta')}\frac{\vert \Delta'\vert^{\vert \Delta'\vert}}{\vert\sigma\vert^{\vert\sigma\vert}\cdot\vert \sigma'\vert^{\vert \sigma'\vert} \cdot\vert\sigma''\vert^{\vert\sigma''\vert}} }\in \GW(F').
    \end{equation}

\end{enumerate}

A part of the proof of our main theorem, the correspondence theorem, will be to compute the trace form $\Tr_{F'/K'}(\Wel^{\A^1}_{F'}(C'))\in \GW(K')$ in the case that all the edges of $\Delta'$ are odd.
Pick's theorem \ref{thm:Pick} implies that if all three edges of a lattice triangle are odd, then the double area of the triangle is also odd. It follows that $\dim_{K'}F'=\frac{\vert\Delta'\vert}{\vert\sigma\vert\vert\sigma'\vert}$ is odd.
Hence, by Proposition \ref{prop:traces} we have that
\begin{equation}
    \label{eq:triangleoddlengths}
    \begin{aligned}
    \Tr_{F'/K'}\left(\Wel_{F'}^{\A^1}(C')\right)&= \Tr_{F'/K'}\left(\qinv{(-1)^{\operatorname{Int}(\Delta')}\frac{\vert\Delta'\vert}{\vert\sigma\vert\vert\sigma'\vert\vert\sigma''\vert}}\right)\\
    &=\qinv{(-1)^{\operatorname{Int}(\Delta')}\vert \sigma''\vert}+\frac{\frac{\vert\Delta'\vert}{\vert\sigma\vert\vert\sigma'\vert}-1}{2}h\in \GW(K').
    \end{aligned}
\end{equation}
 
\subsection{$\Wel_{\kappa(C_l)}^{\A^1}(C_l)$ when $\Delta_l$ is a parallelogram}
In this subsection 
we compute the quadratic weights of the rational curve $C'\coloneqq C_l$, where $\Delta'\coloneqq\Delta_l$ is a parallelogram in $S$. 
We assume we know the coefficients of the defining polynomial $f'$ on two adjacent edges and assume that these coefficients are defined over some finite field extension $K'$ of $k$.
Recall from Proposition \ref{prop:properties Ci} that the $f'$ is of the form $x^sy^r(\alpha x^a +\beta y^b)^p(\gamma x^c+\delta y^d)^q$. Shustin shows that the choice of $\alpha,\beta,\gamma,\delta$ is uniquely determined by the known coefficients and that $\alpha,\beta,\gamma,\delta\in K'^\times$ \cite[Step 2 in~$\S 3.7$]{Shustin}. Next, we compute the quadratic weight of this unique curve, that is we provide is a quadratic enrichment of \cite[Lemma 3.6]{Shustin} for the primitive case.

\begin{prop}
\label{thm:parallelogram}
    Given integers $a,b,c,d$ such that $\gcd(a,b)=\gcd(c,d)=1$ and $ad\neq bc$, and given non-zero coefficients $\alpha,\beta,\gamma,\delta\in K'^{\times}$, the curve defined by
    \[(\alpha x^a+\beta y^b)(\gamma x^c+\delta y^d)=0\]
    has $\vert ad-bc \vert $ nodes.
    None of the nodes lies on the toric boundary, and the quadratic weight of the curve equals $\langle 1\rangle\in \GW(K')$.
\end{prop}

\begin{proof}
    Let $G_1\coloneqq \alpha x^a+\beta y^b$ and $G_2\coloneqq \gamma x^c+\delta y^d$ and $G\coloneqq G_1\cdot G_2$.
    Then
    \[G_x\coloneqq \frac{\partial G}{\partial x}=\alpha ax^{a-1}\cdot G_2+\gamma c x^{c-1}\cdot G_1\]
    and 
    \[G_y\coloneqq \frac{\partial G}{\partial y}=\beta by^{b-1}\cdot G_2+\delta d y^{d-1}\cdot G_1.\]
    For $G$, $G_x$ and $G_y$ to vanish simultaneously, we need $G_1$ and $G_2$ to vanish, or equivalently
    \begin{equation}
    \label{eq:condition nodes}
    \alpha x^a=-\beta y^b\text{ and }\gamma x^c=-\delta y^d.
    \end{equation}
    These equations have $\vert ad-bc\vert$ zeros assuming $x\neq 0$ and $ y\neq 0$ (see for example \cite{Sturmfels} or \cite{JPP}).
    In order to compute $-\det \operatorname{Hessian}(G)$ at a node $z$, put $A\coloneqq \alpha a x^{a-1}\delta d y^{d-1}$ and put~$B\coloneqq \beta b y^{b-1}\gamma c x^{c-1}$.
    Then, 
    \[G_{xx}=\frac{\partial G_x}{\partial x}=\alpha a(a-1) x^{a-2}G_2+2\alpha a x^{a-1}\gamma cx^{c-1}+\gamma c(c-1) x^{c-2}G_1,\]
    \[G_{yy}=\frac{\partial G_y}{\partial y}=\beta b(b-1) y^{b-2}G_2+2\beta b y^{b-1}\delta d y^{d-1}+\delta d(d-1) y^{d-2}G_1,\]
    \[G_{xy}=G_{yx}=A+B.\]
    By \eqref{eq:condition nodes}, we have that at the nodes, these functions evaluate to
    \[G_{xx}(z)=2\alpha a x^{a-1}c\gamma x^{c-1}
    \text{ and }
    G_{yy}(z)=2\beta b y^{b-1}d\delta y^{d-1}.\]
    Note that $G_{xx}\cdot G_{yy}(z)=4AB$. Hence, we get that
    \begin{align*}
        -\det \operatorname{Hessian}(z)=-(4AB-(A^2+2AB+B^2))=(A-B)^2,
    \end{align*}
    which is a square.
    Let $E$ be the $K'$-algebra defined by the nodes. Then
    \[\Wel_{K'}(C')=\qinv{ N_{E/K'}((A-B)^2)}=\qinv{ (N_{E/K'}(A-B))^2}=\qinv{1}\in \GW(K')\]
    where $C'$ is the curve defined by $G$.
\end{proof}
If the curve $C_l$ is defined by $(\alpha x^a +\beta y^b)^p(\gamma x^c+\delta y^d)^q$, with~$pq>1$, then its singularities are not nodal. However, its  components are branches that can be deformed independently, creating a grid consisting of $pq$ components intersecting as in the primitive case. Hence, the contribution to the quadratic weight is the product of the contribution of each two connected components intersecting. Since the quadratic weight for the primitive case is~$\qinv{1}$, the quadratic weight for the general case is a power of $\qinv{1}$ which equals $\qinv{1}$.
\subsection{$\Wel_{\kappa(C_{[\sigma]})}^{\A^1}(C_{[\sigma]})$ when $C_{[\sigma]}$ is a deformation pattern}
In this subsection we quad\-ratically enrich \cite[Lemma 3.9]{Shustin}. That is, for an extended edge $[\sigma]$ of $S$ with lattice length~$m=\vert\sigma\vert\ge 2$, we find all possible deformation patterns $C_{[\sigma]}$ and compute their quadratic weights. 
Before proving our main result in this subsection, we recall some facts about the Chebyshef polynomials of the first kind, since these will define our deformation patterns. 
\begin{df}
    The \emph{$m$-th Chebyshef polynomial} $T_m$ of the first kind is the polynomial
    \[\displaystyle T_{m}(x)=
    \frac{(x+\sqrt{x^2-1})^m+(x-\sqrt{x^2-1})^m}{2}.
    \]
\end{df}

Let us remark that the polynomials $T_m$ have coefficients in $\Z$, and that the polynomials satisfy $T_m(-x)=(-1)^mT_m(x)$, i.e., the polynomial $T_m$ is even (odd) if $m$ is even (odd, respectively).

\begin{lm}
\label{lm:chebychef}
The first derivative $T_m'(x)$ has $m-1$ simple zeros at \[\gamma_l\coloneqq\frac{\zeta_{2m}^l+\zeta_{2m}^{-l}}{2},\]
    where $\zeta_{2m}$ is a primitive $2m$-th root of unity and $l=1,\dots, m-1$.
    Moreover, the values at these zeros of the polynomial and its second derivative are 
    \[T_m(\gamma_l)=(-1)^l \quad\text{ and }\quad
    T_m''(\gamma_l)=(-1)^l\frac{4m^2}{(\zeta_{2m}^l-\zeta_{2m}^{-l})^2}.\]
\end{lm}

\begin{proof}
The first derivative $T_m'$ is a polynomial of degree $m-1$. The computation of the derivative yields.
    \[\displaystyle T_{m}'(x)=m
    \frac{(x+\sqrt{x^2-1})^m-(x-\sqrt{x^2-1})^m}{2\sqrt{x^2-1}}.
    \]
In order to evaluate this expression at $\gamma_l$, one would need to choose a square root of the expression     
    \[\gamma_l^2-1=\frac{(\zeta_{2m}^l+\zeta_{2m}^{-l})^2}{4}-1=\frac{(\zeta_{2m}^l-\zeta_{2m}^{-l})^2}{4}.\]
However, choosing the opposite square root leave the value of the polynomial $T'_m$ invariant since exchanging the roles of the summand of the numerator cancel the sign in the denominator. Hence, we have that
    \[\displaystyle T_{m}'(\gamma_l)=m
    \frac{\left(\zeta_{2m}^l\right)^m-\left(\zeta_{2m}^{-l}\right)^m}{2\sqrt{\gamma_l^2-1}}=m
    \frac{\left(-1\right)^l-\left(-1\right)^{-l}}{2\sqrt{\gamma_l^2-1}}=0,
    \]
because $\zeta_{2m}^m=-1$. Since all $\gamma_l$, $l=1,\dots,m-1$ are different, they are the zeros of the polynomial~$T_m'$.
In order to compute the value of $T_m$ at $\gamma_l$, let us remark that due the symmetry of expression, the value is invariant with respect to the choice of square root. Therefore,
    \[\displaystyle T_{m}(\gamma_l)=\frac{\left(\zeta_{2m}^l\right)^m+\left(\zeta_{2m}^{-l}\right)^m}{2}=\frac{\left(\zeta_{2m}^m\right)^l+\left(\zeta_{2m}^{-m}\right)^l}{2}=\left(-1\right)^l.
    \]
Lastly, the same arguments holds for the second derivative $T_m''$, which yield
    \[\displaystyle T_{m}''(\gamma_l)=m^2
    \frac{\left(\gamma_l+\sqrt{\gamma_l^2-1}\right)^m+\left(\gamma_l-\sqrt{\gamma_l^2-1}\right)^m}{2\left(\sqrt{\gamma_l^2-1}\right)^2}
    =m^2
    \frac{4\left(-1\right)^l}{\left(\zeta_{2m}^l-\zeta_{2m}^{-l}\right)^2}.
    \]

\end{proof}
Let $m=\vert\sigma\vert$.
Recall from Definition \ref{df:deformationpattern} that a deformation pattern for $\sigma$ is a rational curve $C_{[\sigma]}$ in $\Tor_k(\Delta_{[\sigma]})$ with $m-1$ nodes as only singularities, where the triangle $\Delta_{[\sigma]}$ is given by the convex hull $\operatorname{Conv}\{(0,0),(0,2),(m,1)\}$. The curve is defined by a polynomial 
\[f_{[\sigma]}=ay^2+byg(x)+c,\]
with 
\[g(x)=a_mx^m+a_{m-2}x^{m-2}+\ldots+a_1x+a_0,\]
with prescribed non-zero coefficients $a$, $b$ and $c$ defined over a finite field extension $F$ of~$k$ and such that $a_{m-1}=0$, that is the coefficient of $yx^{m-1}$ vanishes.
We find all possibilities for the remaining coefficients~$a_i$ in the following theorem and compute the quadratic weights of $C_{[\sigma]}$ for these choices.
\begin{prop}
\label{thm:quadratic 3.9}
\begin{enumerate}
    \item For fixed $a,b,c\in F^\times$ and an integer $m\ge 2$, let
    \[G(x,y)=ay^2+byg(x)+c\]
    with \[g(x)=a_mx^m+a_{m-2}x^{m-2}+\ldots+a_1x+a_0.\]
    There are exactly $m$ choices for the coefficients $a_{m-2},\ldots,a_0$ and one fixed $a_m\in F^\times$ such that $G$ defines a curve $C_{[\sigma]}$ which has exactly $m-1$ singularities which are nodes. 
    \item Let $L_{[\sigma]}$ be the $F$-algebra defined by these choices. Then 
    \[L_{[\sigma]}=\begin{cases}
        \frac{F[u]}{(u^m-1)}& \text{$m$ is odd,}\\
        \frac{F[u]}{(u^m-ac)} & \text{$m$ is even.}
    \end{cases}\]
    \item The quadratic weight of $C_{[\sigma]}$ equals
    \[\Wel_{L_{[\sigma]}}^{\A^1}(C_{[\sigma]})=\begin{cases} \langle 1\rangle\in \GW(L_{[\sigma]})& \text{if $m$ is odd,}\\
    \langle -u_{[\sigma]}\rangle\in \GW(L_{[\sigma]}) & \text{if $m$ is even.}
    \end{cases}\]
\end{enumerate}

\end{prop}
\begin{proof}
    This follows the ideas of the proof of \cite[Lemma 3.9]{Shustin}.
    We are looking for a polynomial $g(x)=a_mx^m+a_{m-2}x^{m-2}+\ldots+a_0$ such that
    \begin{equation}
        G(x,y)=ay^2+byg(x)+c=0 \label{eq:F=0}
    \end{equation}
    defines a curve with exactly $m-1$ nodes as only singularities. 
    In other words, we are looking for $g(x)$ such that there are exactly $m-1$ simple solutions to
    \begin{align}
    \begin{split}
        G(x,y)&=ay^2+byg(x)+c=0\\
        G_x(x,y)&=byg'(x)=0\\
        G_y(x,y)&=2ay+bg(x)=0.
    \end{split}    
    \end{align}
    Here, we write $G_x$ for the partial derivative of $G$ with respect to $x$ and $G_y$ for the partial derivative of $G$ with respect to $y$.
    
    If $y$ vanishes, then $G(x,y)=0$ cannot hold. Thus asking that $G_x(x,y)=0$ amounts to asking that $g'(x)=0$. From the equality $G_y(x,y)=0$, it follows that
    $y=-{bg(x)}/{2a}$ 
    and plugging this into $G(x,y)=0$ yields 
    $g(x)^2={4ac}/{b^2}$. 
    Thus we are looking for $g(x)$ such that the system
    \begin{equation}
    \label{eq:g}
        g(x)^2=\frac{4ac}{b^2},\quad g'(x)=0
    \end{equation}
    has exactly $m-1$ simple solutions. By \cite{caporasoharrisparameterspaces}, see also \cite{caporasoharrisformula}, there are exactly $m$ such polynomials $g$ with the same leading coefficient $a_m\in F^\times$. We study these solutions according to the parity of~$m$.

    \textbf{$m$ odd}: If $m$ is odd, the $m$ different choices of $g$ are 
    \begin{equation}
    \label{eq:g odd}
    g(x)=\frac{2\sqrt{ac}}{b}T_m\left(\mu\sqrt{ac}x\right),\end{equation}
    one for each $m$-th root of unity $\mu$. Here, $T_m$ is the $m$-th Chebyshef polynomial defined above. $T_m$ is an odd polynomial of degree $m$. So in particular, the coefficient of $x^{m-1}$ is zero. It follows directly from Lemma \ref{lm:chebychef} that the $m$ different choices of $g$ in \eqref{eq:g odd} satisfy the conditions in \eqref{eq:g}.
    Furthermore, notice that since $T_m$ is an odd polynomial, the power of $\sqrt{ac}$ in each summand is even and the square root disappears. It follows that the
    $F$-algebra defined by the different choices of $G$ is
    $L={F[u]}/{(u^m-1)}$.

    We compute the quadratic weight of the curves defined by these equations $G$'s in the ring $\GW(L)$. The nodes correspond to the $m-1$ zeros of $g'(x)$, i.e. the zeros of $T_m'\left(\mu\sqrt{ac}x\right)$. 
    By Lemma \ref{lm:chebychef} these are exactly the values of $x$ where $x$ satisfies $\mu\sqrt{ac} x=\gamma_l=(\zeta_{2m}^l+\zeta_{2m}^{-l})/{2}$ where $\zeta_{2m}$ is a primitive $2m$-th root of unity, for $l=1,\ldots,m-1$. 
    Let us denote by $x_l$ these $x$-coordinates and by $y_l$ the corresponding $y$-coordinate of the nodes, for $l=1,\ldots,m-1$, respectively.
    Consequently, $G_{xy}(x_l,y_l)=bg'(x_l)$ vanishes and the determinant of the Hessian at a node equals 
    \begin{align*}
    G_{xx}(x_l,y_l)\cdot G_{yy}(x_l,y_l)&=by_lg''(x_l)\cdot 2a\\
    &=by_l\cdot \frac{2}{b}\left(\sqrt{ac}\right)^3\mu^2T_m''\left(\mu\sqrt{ac}x_l\right)\cdot 2a\\
    &=4ay_l \left(\sqrt{ac}\right)^3\mu^2T_m''\left(\mu\sqrt{ac}x_l\right)
    \end{align*}
    Since $G_y(x_l,y_l)=0$, we have that 
    \[y_l=-\frac{bg(x_l)}{2a}=-\frac{b}{2a}\sqrt{ac}T_m\left(\mu\sqrt{ac}x_l\right),\]
    and thus,
    \begin{align*}G_{xx}(x_l,y_l)\cdot G_{yy}(x_l,y_l)&=4a\cdot \left(-\frac{b}{2a}\sqrt{ac}T_m\left(\mu\sqrt{ac}x_l\right)\right)\cdot \left(\sqrt{ac}\right)^3\mu^2T_m''\left(\mu\sqrt{ac}x_l\right)\\
    &= -2b(ac)^2\mu^2T_m\left(\mu\sqrt{ac}x_l\right)T_m''\left(\mu\sqrt{ac}x_l\right)
    .\end{align*}
    Therefore, the quadratic weight of the curve defined by $G$ equals
    \begin{align*}
    &\qinv{\prod_{l=1}^{m-1}-G_{xx}(x_l,y_l)\cdot G_{yy}(x_l,y_l)}\\
    =&\qinv{\prod_{l=1}^{m-1}-\left(-2b(ac)^2\mu^2T_m\left(\mu\sqrt{ac}x_l\right)T_m''\left(\mu\sqrt{ac}x_l\right)\right)}\\
    =&\qinv{\left(2b(ac)^2\mu^2\right)^{m-1}\cdot \prod_{l=1}^{m-1}T_m\left(\gamma_l\right)T_m''\left(\gamma_l\right)}\\
    \overset{\ref{df:GW}}{=}&\qinv{\prod_{l=1}^{m-1}T_m(\gamma_l)T_m''(\gamma_l)}
    \overset{\ref{lm:chebychef}}{=}\qinv{\prod_{l=1}^{m-1}(-1)^l(-1)^l\frac{4m^2}{(\zeta_{2m}^l-\zeta_{2m}^{-l})^2}}
    =\qinv{\prod_{l=1}^{m-1}(\zeta_{2m}^l-\zeta_{2m}^{-l})^2}\\
    =&\qinv{\prod_{l=1}^{m-1}\zeta_{2m}^{2l}(1-\zeta_{2m}^{-2l})^2}
    =\qinv{\prod_{l=1}^{m-1}\zeta_m^l(1-\zeta_m^l)^2}
    \overset{\ref{lm:sineIdentity}}{=}\qinv{\zeta_m^{\frac{m(m-1)}{2}}m^2}
    =\qinv{1}\in \GW(L).
    \end{align*}

    \textbf{$m$ even}:
    For $m$ even, we have the following $m$ different choices for $g$.
    \begin{equation}
        \label{eq:g even}
        g(x)=\pm\frac{2\sqrt{ac}}{b}T_m\left(\mu_\pm\cdot\sqrt[2m]{ac}\cdot x\right)
    \end{equation}
    where $\mu_+$ is an $m$-th root of $1$ and $\mu_-$ is an $m$-th root of $-1$. Note that all choices for~$g$ have the same leading coefficient $a_m$. Also note that since $T_m$ is an even polynomial it holds that $T_m(x)=T_m(-x)$ and thus there are indeed exactly $m$ different choices for~$g$, namely ${m}/{2}$ for the $m$ different choices for~$\mu_+$ and ${m}/{2}$ for the $m$ different choices for~$\mu_-$. Again Lemma \ref{lm:chebychef} implies that all choices for $g$ satisfy \eqref{eq:g}.
    Note that the $\mu_+$ and $\mu_-$ are all $2m$-th roots of unity. Since $T_m$ is an even polynomial and thus all monomials have even degree,
    the $F$-algebra defined by the different $G$'s equals $L={F[u]}/{(u^m-ac)}$ where~$u=\mu_\pm^2\sqrt[m]{ac}$. 
    
    We compute the quadratic weight of the curve defined by $G$ in $\GW(L)$. Recall that the nodes correspond to the $m-1$ zeros of 
    \[g'(x)=\pm\frac{2\sqrt{ac}}{b^2}\mu_\pm\sqrt[2m]{ac}\cdot T_m'(\mu_\pm\sqrt[2m]{ac}x)\]
    and thus the $x$-coordinate at the nodes satisfies $\mu_\pm\sqrt[2m]{ac}x=\gamma_l$ by Lemma \ref{lm:chebychef}. Again, we call these $x$-coordinates $x_l$ for $l=1,\ldots,m-1$ and the corresponding $y$-coordinates~$y_l$, respectively.
    Just like in the $m$ odd case, $G_{xy}$ vanishes at the nodes and hence the determinant of the Hessian of $G$ at the nodes equals
    \begin{align*}G_{xx}(x_l,y_l)\cdot G_{yy}(x_l,y_l)&=by_lg''(x_l)\cdot 2a\\
    &=by_l\left(\pm \frac{2\sqrt{ac}}{b} \mu_\pm^2\sqrt[m]{ac}\cdot T_m''(\mu_\pm\sqrt[2m]{ac}x_l)\right)\cdot 2a.
    \end{align*}
    Solving for $y_l$ using $G_y(x_l,y_l)=0$ yields
    \[y_l=-\frac{bg(x)}{2a}=\mp\frac{b}{2a}\cdot \frac{2\sqrt{ac}}{b} T_m(\mu_\pm\sqrt[2m]{ac}x_l)=\mp \sqrt{\frac{c}{a}}T_m(\mu_\pm \sqrt[2m]{ac}x_l).\]
    Hence,
    \begin{align*}
    (G_{xx}\cdot G_{yy})(x_l,y_l)&=b
    \left(\mp \sqrt{\frac{c}{a}}T_m(\mu_\pm \sqrt[2m]{ac}x_l)\right)
    \left(\pm \frac{2\sqrt{ac}}{b} \mu_\pm^2\sqrt[m]{ac}\cdot T_m''(\mu_\pm\sqrt[2m]{ac}x_l)\right)\cdot 2a\\
    &=-4ac\mu_\pm^2\sqrt[m]{ac}\cdot T_m(\mu_\pm \sqrt[2m]{ac}x_l)\cdot T_m''(\mu_\pm\sqrt[2m]{ac}x_l).
    \end{align*}
    Recall that in $L$ we have $u=\mu_\pm^2\sqrt[m]{ac}$.
    Hence, the quadratic weight of the curve defined by $G$ equals
    \begin{align*}
    &\mkern-18mu\qinv{\prod_{l=1}^{m-1}-G_{xx}(x_l,y_l)G_{yy}(x_l,y_l)}=\qinv{\prod_{l=1}^{m-1}-\left(-4ac\cdot u\cdot T_m(\mu_\pm \sqrt[2m]{ac}x_l)\cdot T_m''(\mu_\pm\sqrt[2m]{ac}x_l)\right)}\\
    =&\qinv{(4ac)^{m-1}u^{m-1}\prod_{l=1}^{m-1} T_m(\mu_\pm \sqrt[2m]{ac}x_l)\cdot T_m''(\mu_\pm\sqrt[2m]{ac}x_l)}\\
    =&\qinv{u\prod_{l=1}^{m-1}T_m(\gamma_l)T_m''(\gamma_l)}
    \overset{\ref{lm:chebychef}}{=}\qinv{u\prod_{l=1}^{m-1}(-1)^l(-1)^l\frac{4m^2}{(\zeta_{2m}^l-\zeta_{2m}^{-l})^2}}
    =\qinv{u\prod_{l=1}^{m-1}(\zeta_{2m}^l-\zeta_{2m}^{-l})^2}\\
    =&\qinv{u\prod_{l=1}^{m-1}\zeta_{2m}^{2l}(1-\zeta_{2m}^{-2l})^2}
    \mkern-3mu{=}\mkern-3mu\qinv{u\prod_{l=1}^{m-1}\zeta_m^l(1-\zeta_m^l)^2}
    \mkern-3mu\overset{\ref{lm:sineIdentity}}{=}\mkern-3mu
    \qinv{u\zeta_m^{\frac{m(m-1)}{2}}m^2}
    =\qinv{-u}\in \GW(L).
    \end{align*}
\end{proof}
Lastly, part of the proof of our quadratic enrichment of the correspondence theorem will be to compute $\Tr_{L_{[\sigma]}/F}(\Wel^{\A^1}_{L_{[\sigma]}}(C_{[\sigma]}))\in \GW(F)$. If $m$ is odd, we have that
\begin{equation}
\label{eq:defpattertraceoddcase}
\Tr_{L_{[\sigma]}/F}(\Wel^{\A^1}_{L_{[\sigma]}}(C_{[\sigma]}))=\Tr_{L_{[\sigma]}/F}(\langle1\rangle)=\frac{m-1}{2}h+\langle m\rangle
\end{equation}
by Proposition \ref{prop:traces}. If $m$ is even, then
    \begin{equation}
    \label{eq:defpatterntraceeven}
        \Tr_{L_{[\sigma]/F}}(\Wel_{L_{[\sigma]}}^{\A^1}(C_{[\sigma]}))=\Tr_{L_{[\sigma]/F}}(\qinv{-u_{[\sigma]}})=\frac{m}{2}h
    \end{equation}
    by Proposition \ref{prop:traces} and \eqref{eq:hyperbolicform}.
\section{The quadratically enriched correspondence theorem}
\label{section:correspondencetheorem}
\begin{df}[Quadratically enriched multiplicity]
    Let $A$ be a simple tropical curve with dual subdivision $S$, then we define its \emph{quadratically enriched multiplicity} to be
    \begin{equation}
    \label{eq:qW}
    \mult_k(A)\coloneqq\begin{cases} 
    \displaystyle
    \frac{\prod\vert\Delta'\vert}{2}\cdot h &\text{$S$ has an even edge,}\\[5pt]
    \displaystyle
    \frac{\prod\vert\Delta'\vert -1}{2}\cdot h +\qinv{(-1)^{\operatorname{Int}(S)}}&\text{$S$ has only odd edges,}
    \end{cases}
    \end{equation}
    where the products range over all triangles $\Delta'$ in $S$, $\vert \Delta'\vert$ denotes the double Euclidean area of the triangle and $\operatorname{Int}(S)$ is the total number of interior points of the triangles (not the parallelograms) of $S$.
\end{df}
\begin{thm}[Quadratically enriched correspondence theorem]
\label{thm:qect}
Let $k$ be a field of characteristic $0$ or greater than $\operatorname{diam}(\Delta)$ and let $\K=k\Puiseux$.
    Fix a generic configuration of~$n$ $\K$-points $p_1,\ldots,p_n$ in $\Tor_\K(\Delta)$ such that the configuration $x_1,\ldots,x_n$ with $x_i=\val(p_i)$, for~$i=1,\ldots, n$, is $\Delta$-generic.
    If $A$ is a genus~$g$ tropical curve with Newton polygon $\Delta$ passing through the configuration $x_1,\ldots,x_n$,
    then
    \[\mult_k(A)=\sum \Tr_{\kappa(C)/\K}\left(\Wel^{\A^1}_{\kappa(C)}(C)\right)\]
    under the isomorphism of $\GW(k)\cong\GW(\K)$ in Remark \ref{rmk:GW(K)},
    where the sum runs over all genus $g$ curves $C$ in $\Tor_\K(\Delta)$ passing through $p_1,\ldots,p_n$ that tropicalize to $A$.    
\end{thm}

\begin{coro}
\label{coro:NintermsofWandN}
Let~$k$ be a perfect field such that $\operatorname{char}k=0$ or $\operatorname{char}k>\max\{\operatorname{diam}(\Delta),3\}$ and such that $k$ is big enough such that a generic configuration of $n$ $k$-points in $\Tor_k(\Delta)$ exists.
Furthermore, let $\K=k\Puiseux$.
If the toric variety associated with $\Delta$ is a del Pezzo surface, and thus~$N^{\A^1}_\K(\Delta,0,\mathcal{P})=N^{\A^1}_\K(\Delta)$ is invariant, we get 
\[N^{\A^1}_k(\Delta)=N^{\A^1}_\K(\Delta)=N^{\A^1,\operatorname{trop}}_k(\Delta,0)\coloneqq \sum \operatorname{mult}_k^{\A^1}(A)\in \GW(k),\]
where the sum runs over all genus $0$ tropical curves $A$ through a fixed $\Delta$-generic configuration of points. 
In particular, we have that
\[N^{\A^1}_k(\Delta)=N^{\A^1,\operatorname{trop}}_k(\Delta,0)=\frac{N(\Delta,0)-W(\Delta)}{2}\cdot h+W(\Delta)\cdot \qinv{1}.\]
\end{coro}
\begin{rmk}
    One has to be careful when the base field $k$ is finite. It might happen that there does not exist a generic configuration of $n$ $k$-points in $\Tor_k(\Delta)$ since there are not enough $k$-rational points. However,  the field of Puiseux series $\K$  contains arbitrarily big configurations of points in general position.
\end{rmk}
\begin{proof}[Proof of Corollary \ref{coro:NintermsofWandN}]
    In case the toric variety associated with $\Delta$ is del Pezzo, it is shown in  \cite{KLSWCount,KLSWorientation} (see also Subsection \ref{subsection:motivation from classical and a1homotopytheory}) that $N^{\A^1}_\K(\Delta,0,\mathcal{P})$ is independent of the choice of point configuration $\mathcal{P}$. Thus $N^{\A^1}_\K(\Delta,0,\mathcal{P})$ equals $N^{\A^1}_\K(\Delta,0,\mathcal{P}')$ where $\mathcal{P}'$ is a generic configuration of $\K$-points in $\Tor_\K(\Delta)$ that specializes to a generic configuration~$\overline{\mathcal{P}'}$ of $k$-points in $\Tor_k(\Delta)$.
    It follows that under the canonical isomorphism~$\GW(\K)\cong \GW(k)$ in Remark \ref{rmk:GW(K)}, $N^{\A^1}_\K(\Delta,0,\mathcal{P})=N^{\A^1}_\K(\Delta,0,\mathcal{P}')$ maps to $N^{\A^1}_k(\Delta,0,\overline{\mathcal{P}'})$ which is again independent of the choice of point configuration. Hence $N^{\A^1}_k(\Delta)=N^{\A^1}_\K(\Delta)=N^{\A^1,\operatorname{trop}}_k(\Delta,0)$.

    Note that the rank of $N^{\A^1,\operatorname{trop}}_k(\Delta,0)$ equals the tropical complex count $N^{\operatorname{trop}}(\Delta,0)$ and the signature of $N^{\A^1,\operatorname{trop}}_{\R}(\Delta,0)$ equals the tropical real count $W^{\operatorname{trop}}(\Delta,0)$. Since the quadratically enriched tropical multiplicity $\operatorname{mult}_k^{\A^1}(A)$ of a simple tropical curve $A$ has only  $\qinv{1}$'s and $\qinv{-1}$'s as summands, it is completely determined by its rank and signature, and thus the same holds for $N^{\A^1,\operatorname{trop}}_k(\Delta,0)$. It follows that  
    \[N^{\A^1,\operatorname{trop}}_k(\Delta,0)=\frac{N^{\operatorname{trop}}(\Delta,0)-W^{\operatorname{trop}}(\Delta,0)}{2}\cdot h+W^{\operatorname{trop}}(\Delta,0)\cdot \qinv{1}.\]
    By Mikhalkin's non-quadratically enriched correspondence theorem, we get that
    \[N^{\A^1}_k(\Delta)=\frac{N(\Delta,0)-W(\Delta)}{2}\cdot h+W(\Delta)\cdot \qinv{1}.\]
\end{proof}

\subsection{Restoring the algebraic curve from the tropical data}
\label{subsection:restoringthecurve}
We start by recalling Shustin's reconstruction of the algebraic curves through a generic configuration~$\mathcal{P}$ of $n$ $\K$-points  in $\Tor_\K(\Delta)$ that tropicalize to a given tropical curve $A$ \cite[$\S3.7$]{Shustin} (see also \cite[$\S2.5$]{tropicalBook}).
Let $S:\Delta=\Delta_1\cup\ldots\cup\Delta_N$ be the dual subdivision of $A$.
In Section~\ref{section:quadraticenrichments} we found all possible choices for the curves $C_l$, for $l=1,\dots,N$ and deformation patterns~$C_{[\sigma]}$ for each extended edge $[\sigma]\in E^*(S)$ of lattice length at least $2$.
We will find a finite dimensional étale $k$-algebra~$\etale$ defined by all of these choices and the choices for the conditions to pass through a fixed point.
Using the refined patchworking Theorem \ref{thm:patchworking} one can find a unique curve $C$ over $\etale\Puiseux$ and we use the results from Section \ref{section:quadraticenrichments} to compute the quadratic weight $\Wel^{\A^1}_{\etale\Puiseux}(C)\in \GW(\etale\Puiseux)\cong\GW(\etale)$ of this curve. Then we show that $\operatorname{mult}_k^{\A^1}(A)=\Tr_{\etale/k}\left(\Wel_{\etale\Puiseux}(C)\right)$ in $ \GW(k)$. 

\begin{rmk}
    By our assumption on the characteristic of the field $k$ ($\operatorname{char}k=0$ or $\operatorname{char}k$ is big enough, so that all of the polynomials that occur in the construction of the $k$-algebra~$\etale$ are of degree smaller than $\operatorname{char}k$), the finite dimensional $k$-algebra $\etale$ is \`{e}tale and thus is isomorphic to the product of finite separable field extensions
    \[\etale\cong \kappa_1\times\ldots \times \kappa_s.\]
    Geometrically, this means that there exist exactly $s$ curves $C^{(1)},\ldots,C^{(s)}$ tropicalizing to~$A$ and the ground field of $C^{(i)}$ is  $\kappa_i\Puiseux$ for $i=1,\ldots,s$. In particular, 
    \[\Tr_{\etale\Puiseux/k\Puiseux}\left(\Wel_{\etale\Puiseux}^{\A^1}(C)\right)=\sum_{i=1}^s\Tr_{\kappa_i\Puiseux/\K}\left(\Wel^{\A^1}_{\kappa_i\Puiseux}(C^{(i)})\right)=N^{\A^1}_{\K}(\Delta,0,\mathcal{P}).\]
\end{rmk}

To find the $k$-algebra $\etale$ and the curve $C$ over $\etale\Puiseux$ we will use the following proposition.
\begin{prop}[\cite{ShustinArxivLongVersion}]
\label{prop:Shustin's graph construction}
Let $A$ be a simple irreducible curve of degree $\Delta$ passing through a $\Delta$-generic configuration of points $x_1,\ldots,x_r\in\R^2$.
\begin{enumerate}
    \item Let $G_0$ be the union of extended edges in $E^*(A)$ passing through $x_1,\ldots,x_n$. Then, the vertices of $G_0$ different from the four-valent vertices of $A$ have valency at most~$2$. 
    \item The subgraph $S_0$ of the dual subdivision $S$ consisting of the dual edges of $G_0$ is connected and contains all vertices $V(S)$ of $S$.
    \item 
    We construct subgraphs $G_i$ of $A$ inductively.
    Let 
    \[M_i\coloneqq\{ \epsilon\in [\epsilon']\in E^*(A), \epsilon'\in E(A\setminus G_i)\text{, $\epsilon'$ has a vertex which is bivalent vertex of $G_i$}\}\]
    and set
    \[G_{i+1}\coloneqq G_i\cup \bigcup_{\epsilon\in M_i}\epsilon.\] 
    Then at each step, it holds that for the $\epsilon\in M_i$, $\epsilon$ has exactly one vertex which is a bivalent vertex of $G_i$.
    Furthermore, $A=G_m$ for some positive integer $m$.
\end{enumerate}
\end{prop}

\begin{proof}
The proposition follows from the $\Delta$-genericity of the configuration of points. For the first assertion, if a vertex $v$ of $G_0$ has valency three, any of its adjacent edges is determined by the other two due to the balancing condition; hence, either the same curve is determined by $r-1$ points or there exists a positive dimensional family of tropical curves passing through the points $x_1,\ldots,x_r$ since the number of independent constraints is lower than the rank. In both cases the rank of the curve is off by one. 
For the second claim, a vertex $v\in V(S)\setminus V(S_0)$ is dual to a connected component of $\R^2\setminus A$ whose boundary is comprised entirely by unmarked edges. Hence, it is possible to continuously deform the length of the edges on the boundary creating a one dimensional family of tropical curves of the same degree passing through the points $x_1,\ldots,x_r$. In a similar fashion, if the graph~$S_0$ is disconnected, then there exists an edge $e\in E(S)\setminus E(S_0)$ adjacent to two connected components of $S_0$. The edge~$e$ is dual to a edge that can be deformed in a positive dimensional family of curves of the same degree passing through the points $x_1,\ldots,x_r$, contradicting the genericity of the configuration of points.
For the last statement, notice that the position of the remaining edges in $E^{*}(A)$ is prescribed by $G_0$. Iterating this argument, if an edge $e\in E^{*}(A\setminus G_i)$ is adjacent to two bivalent vertices of $G_i$, then the number of the same curve is determined by $r-1$ points, since by the balancing condition, at least one edge of $G_0$ containing a point of the configuration is determined by neighboring edges.
\end{proof}

\begin{ex}
\label{ex:graphconstruction}
    We illustrate the content of the technical Proposition \ref{prop:Shustin's graph construction} in an example in Figure \ref{fig:graphconstruction}.
    For the marked tropical curve $A$ in Figure \ref{fig:curve+subdivision} we get the graphs $G_0$, $G_1$, $G_2$ and $G_3=A$.
     \begin{figure}
    \begin{tabular}{cc}
    \includegraphics[scale=0.4]{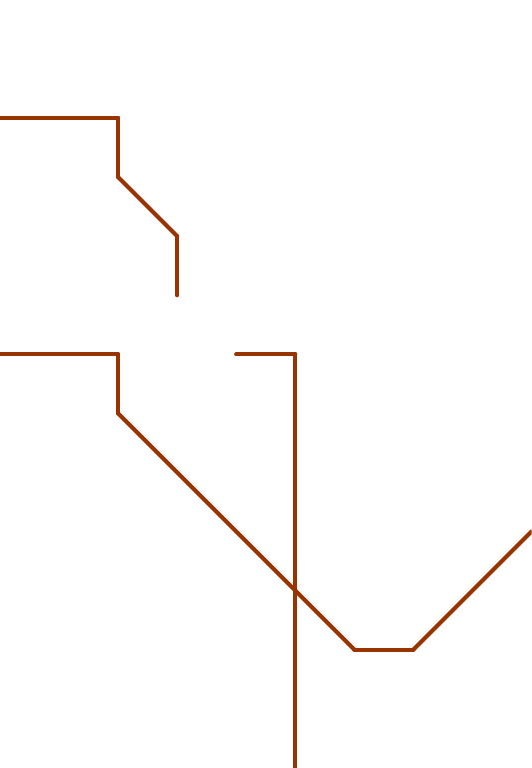}
    &
    \includegraphics[scale=0.4]{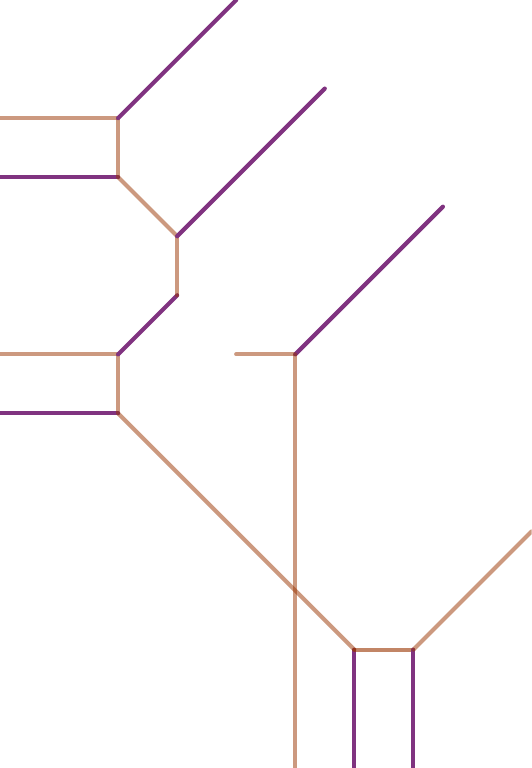}
    \\[12pt]
    \includegraphics[scale=0.4]{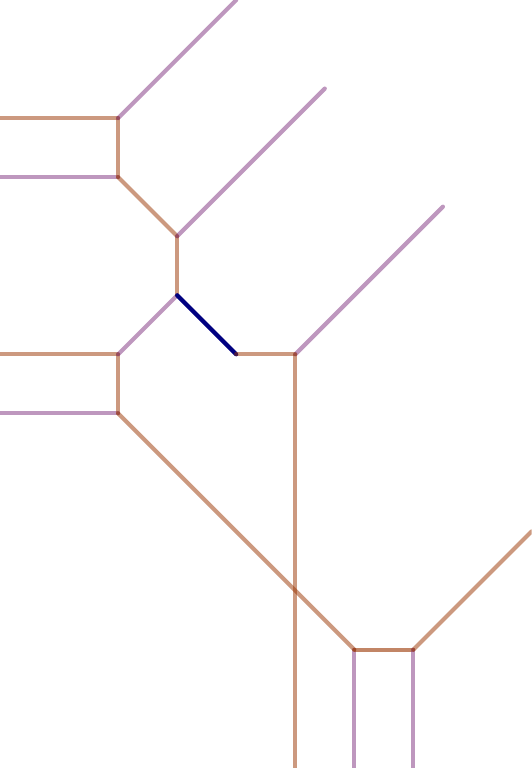}
    &
    \includegraphics[scale=0.4]{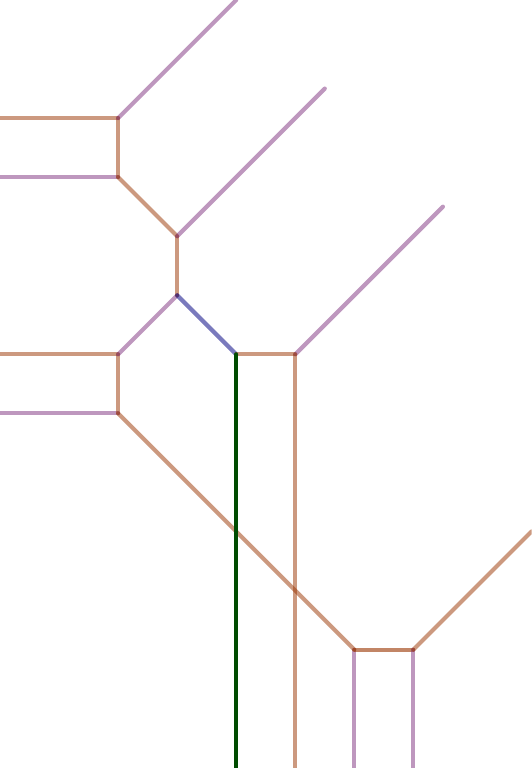}
    \end{tabular}
    \caption{Example of graphs $G_0$, $G_1$, $G_2$ and $G_3$.\label{fig:graphconstruction}}
    \end{figure}
\end{ex}

We start with the following set up.
Assume that we are given a genus $g$ tropical curve~$A$ with Newton polygon $\Delta$ passing through $\Delta$-generic points $x_1,\ldots, x_n\in \Q^2$. Since the configuration $x_1,\ldots,x_n$ is $\Delta$-generic, the curve $A$ has to be simple. 
Assume that the configuration $p_1,\ldots,p_n\in \Tor_\K(\Delta)$ of $\K$-points is generic, such that none of the points lies on the toric boundary and such that $\val(p_i)=x_i$ for $i=1,\ldots,n$.
The algebraic curves in $\Tor_\K(\Delta)$ through $p_1,\ldots,p_n$ that tropicalize to $A$ are defined by polynomials of the form
\begin{equation}
\label{eq:equationcurve}
f(x,y)=\sum_{(i,j)\in \Delta\cap\Z^2}(a_{ij}+c_{ij}(t))t^{\nu(i,j)}x^iy^j\end{equation}
with $a_{ij}\in K$, where $K/k$ is some finite separable field extension and $c_{ij}(t)\in K\Puiseux$ with $\val(c_{ij}(t))<0$.
To reconstruct these algebraic curves, we use the refined patchworking Theorem~\ref{thm:patchworking}. So we need to find all the possibilities for the coefficients $a_{ij}$'s and the deformation patterns $f_{[\sigma]}$.
See Example~\ref{ex:reconstruction} for a detailed example of the whole process.
We proceed in the following steps.
\begin{enumerate}
    \item The tropical curve $A$ defines a dual subdivision 
    \[S:\Delta=\Delta_1\cup\ldots\cup \Delta_N\]
    of $\Delta$ and a piecewise linear function 
    \[\nu\colon \Delta\rightarrow \R\]
    (see Subsection \ref{subsection:tropicalization}).
        
    \item 
    Next, we find all possibilities for the curves $C_l$ for $\Delta_l\in S$. That is, the curves defined by
    \begin{equation}
    \label{eq:fl}
    f_l(x,y)=\sum_{(i,j)\in \Delta_l\cap\Z^2}a_{ij}x^i y^j.\end{equation}
    So we have to determine all possible choices for the coefficients $a_{ij}$ in 
\eqref{eq:equationcurve}.
    \begin{claim}
        Let $(i_0,j_0)$ be a vertex of a marked edge of $S$. We set $a_{i_0j_0}=1$. Then
        asking that 
        \[f(p_1)=\ldots=f(p_n)=0\]
        determines unique $a_{ij}\in k^{\times}$ for $(i,j)\in V(S)$ and for $(i,j)\in \sigma_l$ for the extended marked edges $\sigma_1,\ldots,\sigma_n$.
    \end{claim}
    \begin{proof}
        This follows directly from the computations in \cite[Step 2 in $\S3.7$]{Shustin} which we briefly recall.

        Assume $\sigma_l=[w_1^{(l)}, w_2^{(l)}]$ is the edge dual to the marked edge of $A$ with the point $x_l$.
        Without loss of generality we can assume that $\nu\vert_{\sigma_l}=0$ and strictly positive elsewhere. Then the corresponding point $p_l$ is of the form
        $p_l=(\xi, \eta)$ with~$\xi =\xi_0+ h.o.t.$ and $\eta=\eta_0+h.o.t.$ with~$\xi_0,\eta_0\in k^\times$.
        Applying a suitable coordinate change, we can also assume that $\sigma_l=[(0,0),(m,0)]$ where $m=\vert \sigma_l\vert$. 
        The point condition $f(p_l)=0$ implies that 
                \begin{equation}
                \label{eq:f(pl)}
        a_{00}+a_{10}\xi_0+\ldots+a_{m0}\xi_0^m=0.
        \end{equation}
        If $\vert \sigma_l\vert=m=1$, the equation~\eqref{eq:f(pl)}
        corresponds to a linear equation in $a_{w_1^{(l)}}$ and~$a_{w_1^{(l)}}$ before the coordinate change.

        If $\sigma_l=\Delta_r\cap\Delta_s$, then $m=\vert\sigma_l\vert$ can be greater than $1$. Recall from Proposition~\ref{prop:properties Ci} that the curves $C_r\subset \Tor_k(\Delta_r)$ and $C_s\subset \Tor_k(\Delta_s)$ should meet along the toric boundary $\Tor_k(\sigma)$ at $(\xi_0,\eta_0)$ with multiplicity $m$.
        Consequently, the polynomial $a_{00}+a_{10}x+\ldots+a_{m0}x^m$ has to be of the form
        \begin{equation}
            a_{m0}(x-\xi_0)^m
        \end{equation}
        which leads to $m$ linear equations in $a_{00},a_{10},\ldots,a_{m0}$, or equivalently, to $m$ linear equations in the $a_{ij}$ for the $(i,j)\in \sigma_l\cap\Z^2$ in the coordinates we started with.

        For each parallelogram $\Delta_l$ in the dual subdivision with vertices $w_1$, $w_2$, $w_3$ and~$w_4$, the coefficients $a_{w_1}$, $a_{w_2}$, $a_{w_3}$ and $a_{w_4}$ have to satisfy 
        \[a_{w_1}a_{w_3}=a_{w_2}a_{w_4},\]
        since $f_l$ is of the form described in Proposition \ref{prop:properties Ci} (3).

        By Proposition \ref{prop:Shustin's graph construction} (2) the subgraph of the dual subdivision $S$ consisting of all extended marked edges is connected and contains all vertices $V(S)$ of $S$.
        So setting one coefficient $a_{i_0j_0}=1$, the linear equations described above and the equations for the vertices of the parallelograms of $S$ yield an unique coefficient $a_{ij}$ on $V(S)$ and on the marked edges $\sigma_i$ for $i=1,\ldots,n$.
    \end{proof}
    \item 
    Let $\Delta_l$ be a triangle in $S$ with edges $\sigma_l$, $\sigma_l'$ and $\sigma_l''$ and assume that we have found the coefficients $a_{ij}\in K_l$ of $f_l$ for $(i,j)\in (\sigma\cup\sigma')\cap\Z^2$ where $K_l$ is some finite $k$-algebra. We want to find the remaining coefficients $a_{ij}$ of $f_l$ in \eqref{eq:fl}.
    Let $F_l\coloneqq F'$ be the $K_l$-algebra defined in Lemma \ref{lm:F'}. 
    Then Lemma \ref{lm:F'} (or alternatively \cite[Lemma 3.5]{Shustin}) tells us that there are
    \[\dim_{K_l}F_l=\frac{\vert \Delta_l\vert }{\vert \sigma_l\vert\cdot \vert \sigma_l'\vert}\]
    possibilities for the curves $C_l$ in $\Tor_k(\Delta_l)$, or equivalently, for the remaining coefficients $a_{ij}$. By Lemma \ref{lm:F'}, we have that 
    \begin{equation}
    \label{eq:Fl}
        F_l= \frac{K_l[\beta_l]}{\left(\beta_l^{\frac{\vert \Delta_l\vert}{\vert\sigma_l\vert\cdot\vert\sigma_l'\vert}}+d_l\right)},
    \end{equation}
    with $d_l\coloneqq(-1)^{\frac{m(q-p)}{d_p}}\left(\frac{-\epsilon_1}{a}\right)^{\frac{p}{d_p}}\frac{\epsilon_2}{b}\in K_l^\times$ (see \eqref{eq:F'}), is the $K_l$-algebra defined by the different possibilities.
    
    \item Now we use the graph construction in Proposition \ref{prop:Shustin's graph construction} to find all possibilities for a global compatible choice of $\bigcup_{l=1}^NC_l\subset\bigcup_{l=1}^N\Tor_F(\Delta)$, that is we find compatible choices for $a_{ij}$ with $(i,j)\in \Delta\cap \Z^2$ where $F$ is some finite $k$-algebra, using the previous step for each triangle in $S$. 
    Let $E^*(S)$ the \emph{set of extended edges} defined in Subsection \ref{subsection:notation}.
    Note that $E^*(S)$ is the disjoint union of the set of \emph{marked extended edges} $E^*_{\text{marked}}(S)$ and the set of \emph{unmarked extended edges} $E^*_{\text{unmarked}}(S)$.
    Similarly, we define an equivalence class on the set of edges of $A$: We say two edges are equivalent if the dual edges are. Figure \ref{fig:extended edges+curve 2} illustrates this definition in an example.
    By Proposition \ref{prop:Shustin's graph construction}(1) in the dual subdivision $S$ each triangle $\Delta_l$ has at most $2$ marked edges (this is because every $3$-valent vertex of $A$ has at most valency $2$ in~$G_0$).
    
    We start with the triangles in $S$ with two marked edges, that is all triangles dual to a $2$-valent vertex in $G_0$, where $G_0\subset A$ is the subgraph defined in Proposition~\ref{prop:Shustin's graph construction}, and find all possible choices for the coefficients $a_{ij}$ by applying Lemma \ref{lm:F'} for these triangles. In particular, we get a finite choice for the $a_{ij}$ for the third unmarked edges. 
    In the next step, we proceed with the triangles in $S$ which now have prescribed coefficients on two edges, that is exactly the triangles dual to a $2$-valent vertex of $G_1$ where $G_1\subset A$ is the subgraph defined in Proposition \ref{prop:Shustin's graph construction}(2). 
    Again, we find the remaining coefficients $a_{ij}$ exactly as in the previous step, in particular for the third edges.
    We continue until we have found all possible choices for the coefficients $a_{ij}$, for each $(i,j)$ a lattice point in a triangle of $S$. 
    Finally, by \cite[$3.7$ Step $2$]{Shustin}, for a parallelogram $\Delta_l$ in $S$ for which we know the coefficients on $2$ of its edges, there is a unique choice for the remaining coefficients $a_{ij}$ of $f_l$.

    We want to describe the $k$-algebra $F$ defined by all possibilities of the choices of coefficients $a_{ij}$. We define $F$ recursively:
    Let
    \[F^{(0)}\coloneqq \frac{k[\{\beta_l\}_{G_0}]}{\{(\beta_l^{\frac{\vert\Delta_l\vert}{\vert\sigma_l\vert\vert\sigma'\vert}}-d_l)\}_{G_0}}\]
    where we adjoin a variable $\beta_l$ for each triangle $\Delta_l$ dual to a bivalent point in $G_0$. Then we set 
    \[F^{(i)}\coloneqq \frac{F_{i-1}[\{\beta_l\}_{G_i}]}{\{(\beta_l^{\frac{\vert\Delta_l\vert}{\vert\sigma_l\vert\vert\sigma'\vert}}-d_l)\}_{G_i}}\]
    where we adjoin a variable $\beta_l$ for each triangle $\Delta_l$ dual to a new bivalent point in~$G_i$,
    for $i\ge 1$. In both definitions above, the element $d_l$ is as in \eqref{eq:Fl}.
    We have seen that for some $m\ge 0$ we have $G_m=A$. Set
    \[F\coloneqq F^{(m)}.\]
    The $k$-dimension of $F$ is
    \[\dim_kF=\prod_{\text{$\Delta_l$ triangle in $S$}}\frac{\vert\Delta_l\vert}{\vert\sigma_l\vert\vert\sigma_l'\vert}.\]
    Any marked interior (that is not a boundary edge of $\Delta$) edge appears twice in the product while any unmarked interior edge appears once. Since $A$ is simple, the lengths of the boundary edges are all equal to $1$ and one gets 
    \begin{align*}\dim_kF&=\prod_{\Delta_l\text{ a triangle in }S}\vert \Delta_l\vert\cdot\prod_{[\sigma]\in E^*_{\text{marked}}}\vert\sigma\vert^{-2}\cdot \prod_{[\sigma]\in E^*_{\text{unmarked}}}\vert\sigma\vert^{-1}\\
    &=\prod_{\Delta_l\text{ a triangle in }S}\vert \Delta_l\vert\cdot \prod_{\sigma\in E^*(A)}\vert \sigma\vert^{-1}\cdot \prod_{i=1}^n\vert \sigma_i\vert^{-1}.\end{align*}
    This equals the number of possibilities for the curve $\bigcup_{l=1}^NC_l\subset \bigcup_{l=1}^N\Tor_k(\Delta_l)$.
    
    \item Now, we have found all possibilities for the polynomials $f_1,\ldots,f_N$ defyning the curves~$C_l\subset \Tor_{K_l}(\Delta_l)$. In order to use the refined algebraic patchworking Theorem~\ref{thm:patchworking}, we need to choose a deformation pattern for each extended common edge~$\sigma=\Delta_l\cap\Delta_j$ of lattice length $\ge 2$. 
    
    By Lemma 3.9 in \cite{Shustin} or Proposition \ref{thm:quadratic 3.9} there are~$\vert \sigma \vert$ choices for a deformation pattern for an extended edge $\sigma$ in $S$ 
    and the $F$-algebra defined by these~$\vert\sigma\vert$ choices equals 
    \begin{equation}
        \label{eq:Lsigma}
        L_{[\sigma]}=\begin{cases}
        \frac{F[u_{\sigma}]}{\left(u_{\sigma}^{\vert\sigma\vert}-1\right)} & \text{if $\vert \sigma\vert$ is odd,}\\
        \frac{F[u_{\sigma}]}{\left(u_{\sigma}^{\vert\sigma\vert}-a_\sigma c_\sigma\right)} & \text{if $\vert \sigma\vert$ is even,}
        \end{cases}
    \end{equation}
    where $a_\sigma c_\sigma$ is as in Proposition \ref{thm:quadratic 3.9}.
    Set 
    \begin{equation}
    \label{eq:L}
        L\coloneqq\frac{F\left[\{u_{\sigma}:[\sigma]\in E^*(A)\}\right]}{\left(\{u_{\sigma}^{\vert\sigma\vert}-b_\sigma:\text{$[\sigma]\in E^*(A)$}\}\right)},
    \end{equation}
    where $b_\sigma=1$ for $\vert\sigma\vert $ odd and $b_\sigma=a_\sigma c_\sigma$ when $\vert \sigma\vert$ is even.
    That is, the $F$-algebra~$L$ is  defined by all possible choices for the deformation patterns for all extended edges.
    
    \item Finally, there are additional $\vert \sigma_i\vert$ choices for the refined data to pass through the points $p_1,\ldots,p_n$ for each marked edge $\sigma_i$ as explained in Subsection \ref{subsection:conditionspoints}. 
    The $L$-algebra defined by the $\vert\sigma_i\vert$ choices is of the form 
    \begin{equation}
        \etale_i\coloneqq \frac{L[w]}{(w^{\vert\sigma_i\vert}-d_i)}
    \end{equation}
    for $d_i=-\eta_i^0m_ia_{01}a_{m_i0}^{m_I-1}\in L^\times$ (see \eqref{eq:conditionlengthbigger1}).
    Set
    \begin{equation}
    \label{eq:E}
        \etale\coloneqq \frac{L[w_1,\ldots,w_n]}{(w_1^{\vert\sigma_1\vert}-d_1,\ldots,w_n^{\vert\sigma_n\vert}-d_n)}.
    \end{equation}
    
    \item By the refined Patchworking Theorem \ref{thm:patchworking}, the data of the $a_{ij}$ (Step (4)), deformation patterns (Step (5)) and refined conditions to pass through a point (Step (6)) define a unique curve $C$ that tropicalizes to $A$. Note that he $k$-algebra $\etale$ from step (6) has dimension
    \begin{align*}\dim_k\etale&=\dim_kF\cdot \dim_FL\cdot \dim_L\etale\\
    &= \left(\prod_{\text{triangles }\Delta'} \vert \Delta'\vert\cdot \prod_{\sigma\in E^*{A}}\vert\sigma\vert^{-1}\cdot \prod_{i=1}^r\vert\sigma_i\vert^{-1}\right)\cdot \left(\prod_{\sigma\in E^*{A}}\vert\sigma\vert\right)\cdot \left(\prod_{i=1}^r\vert\sigma_i\vert\right)\\
    &=\prod_{\text{triangles }\Delta'}\vert\Delta'\vert \end{align*}
     which equals $\mult_\C(A)$. This is expected, since when $k$ is algebraically closed, the~$k$-algebra $\etale$ equals the product of $\dim_k\etale$ copies of $k$. For each copy of~$k$ there is a curve tropicalizing to $A$. In other words, this is the content of a correspondence theorem over an algebraically closed field.
\end{enumerate}

We illustrate the whole reconstruction for the tropical curve from Example \ref{ex:graphconstruction}.
\begin{ex}
    \label{ex:reconstruction}
    Consider the tropical curve $A$ in Figure \ref{fig:curve+subdivision} with Newton polygon $\Delta$ given by $\operatorname{Conv}\{(0,0),(4,0),(0,4)\}$ and markings $x_1,\ldots,x_{11}$. 

    \begin{figure}
    \begin{tabular}{cc}
    \includegraphics[scale=0.15]{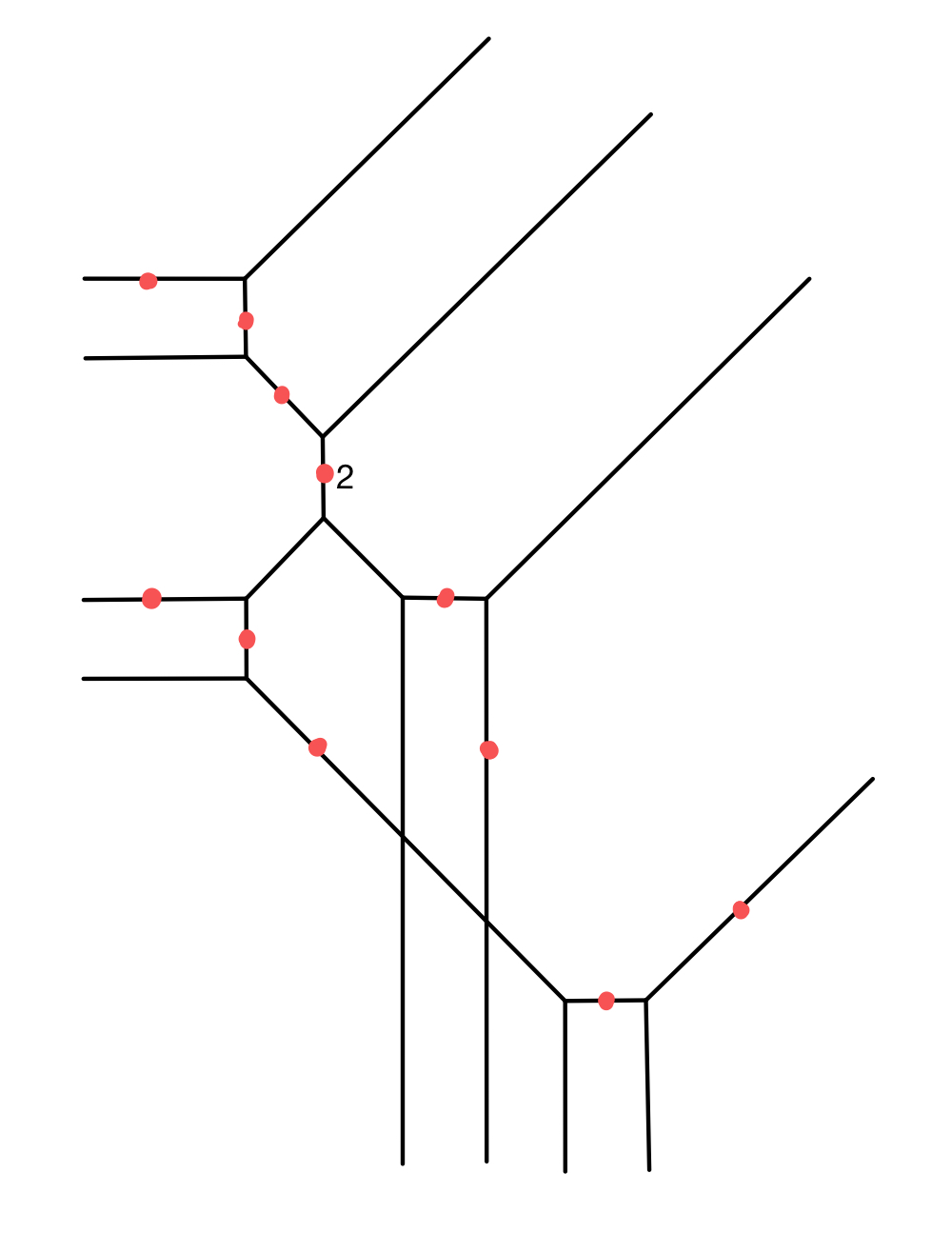}&
    \includegraphics[scale=0.15]{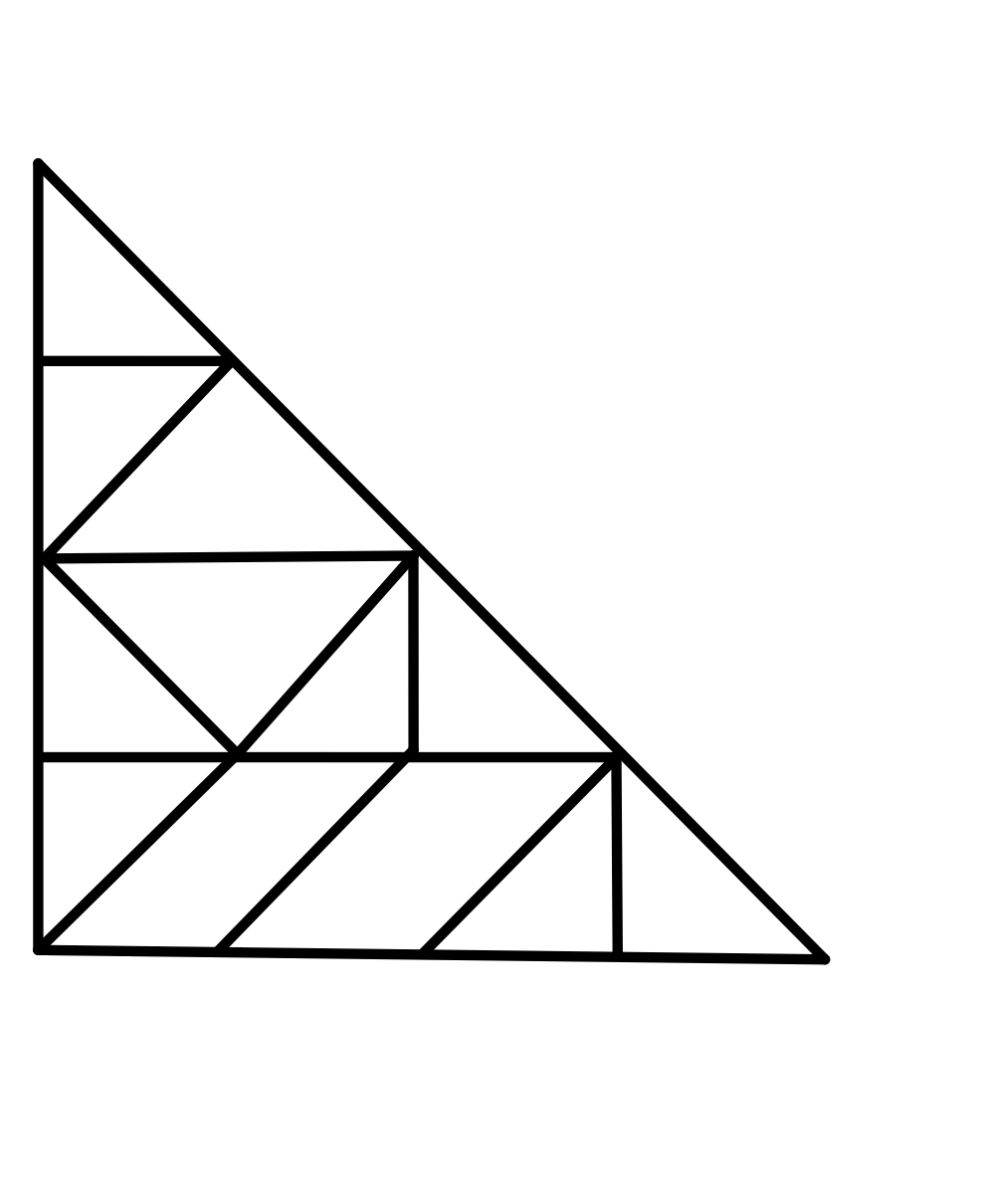}
    \end{tabular}
    \caption{A marked tropical curve and its dual subdivision.\label{fig:curve+subdivision}}
    \end{figure}
    We aim to reconstruct all rational curves in $\Tor_{\K}(\Delta)$, i.e., degree $4$ rational curves in $\mathbb{P}^2_\K$, that go through given points $p_1,\ldots,p_{11}$ with $\val(p_i)=x_i$ and which tropicalize to~$A$. 
    Figure \ref{fig:extended edges+curve 2} shows which edges in the subdivision in Figure \ref{fig:curve+subdivision} are identified by the equivalence relation $\sim$ that identifies \emph{extended edges}. 
    \begin{figure}
    \begin{tabular}{cc}
    \includegraphics[scale=0.15]{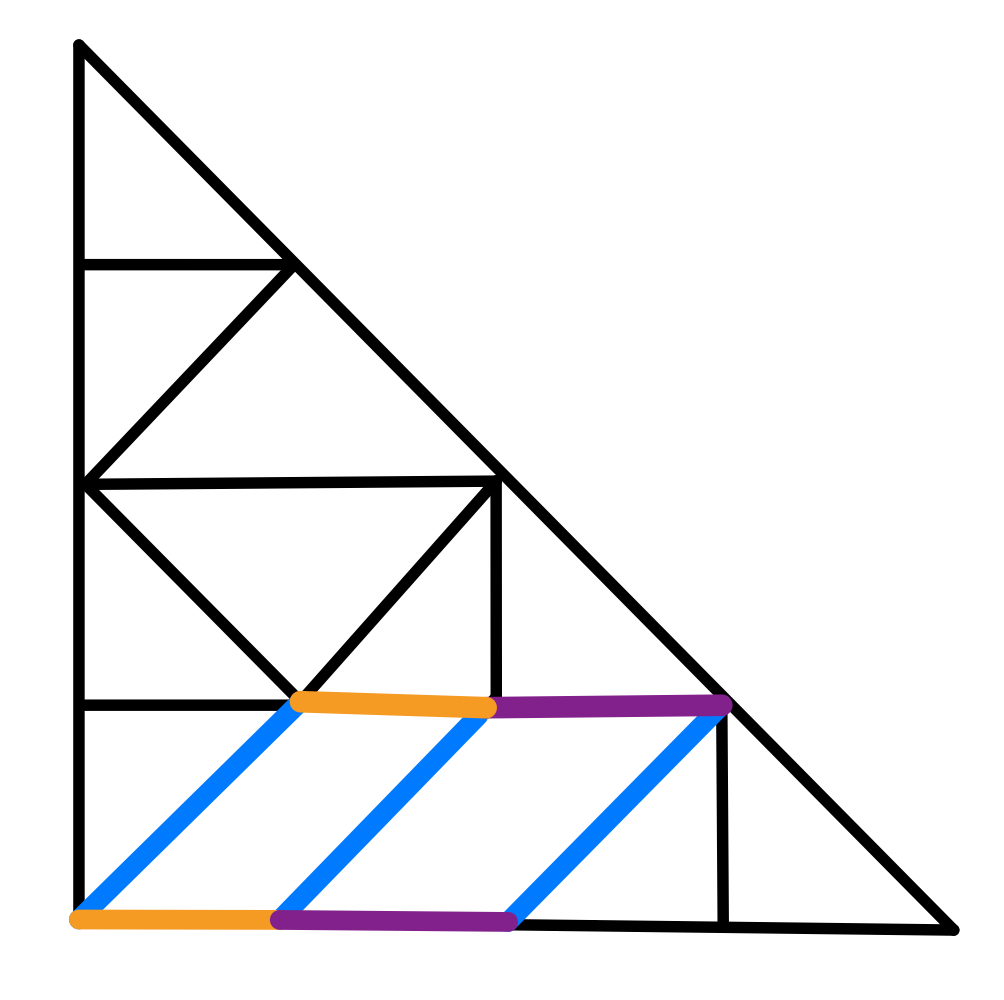}&
    \includegraphics[scale=0.15]{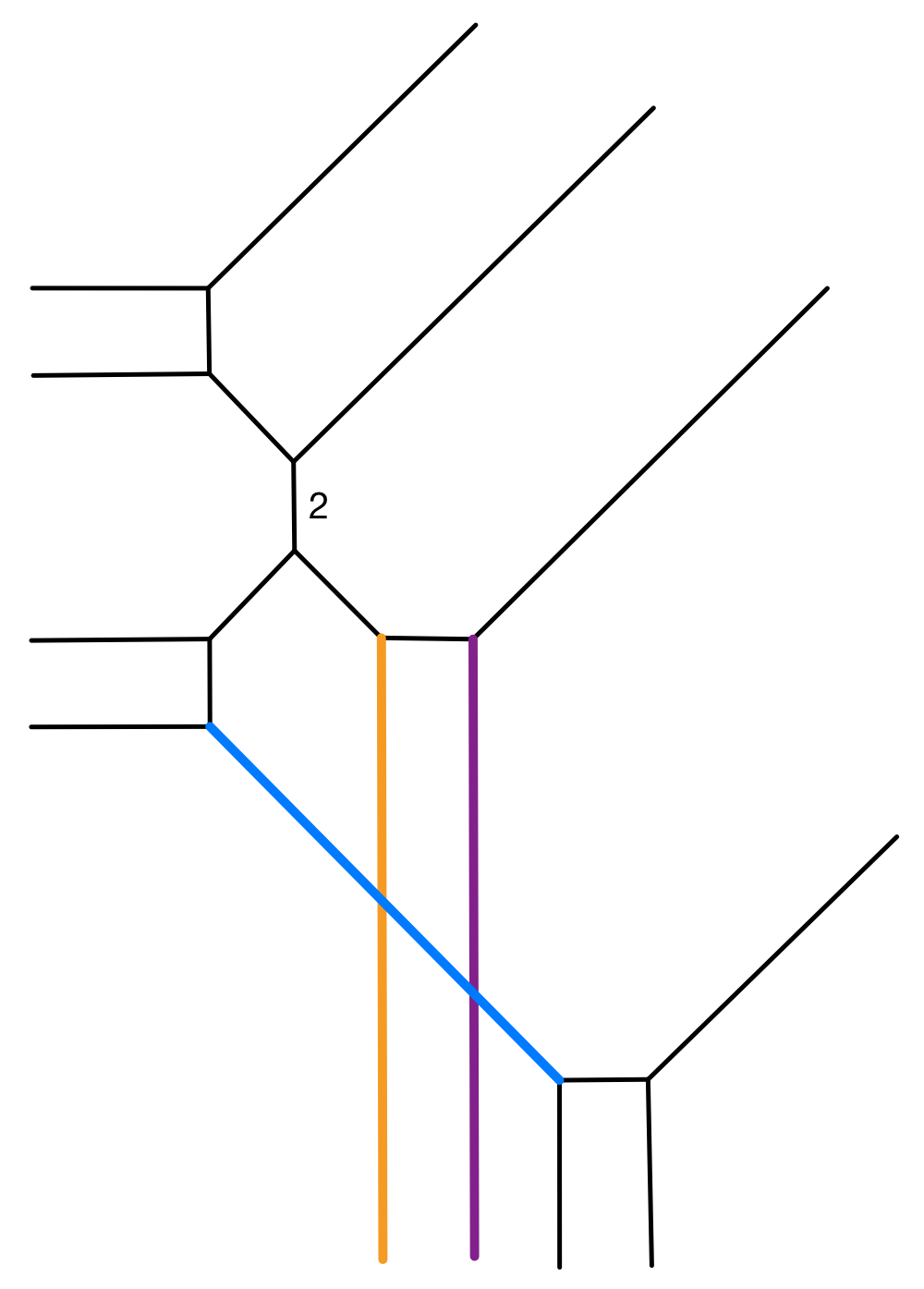}
    \end{tabular}
    \caption{Extended edges.\label{fig:extended edges+curve 2}}
    \end{figure}
    \begin{figure}
    \begin{tabular}{cc}
    \includegraphics[scale=0.2]{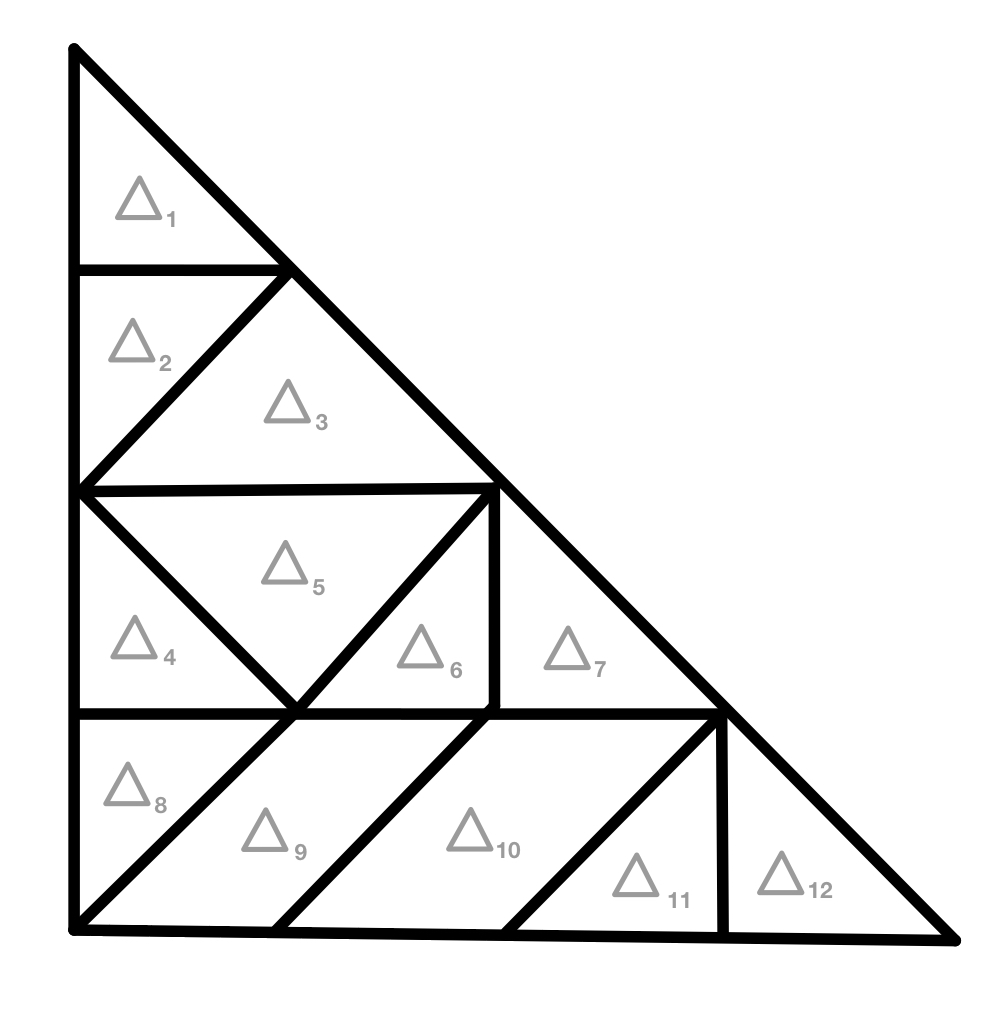}&
    \includegraphics[scale=0.19]{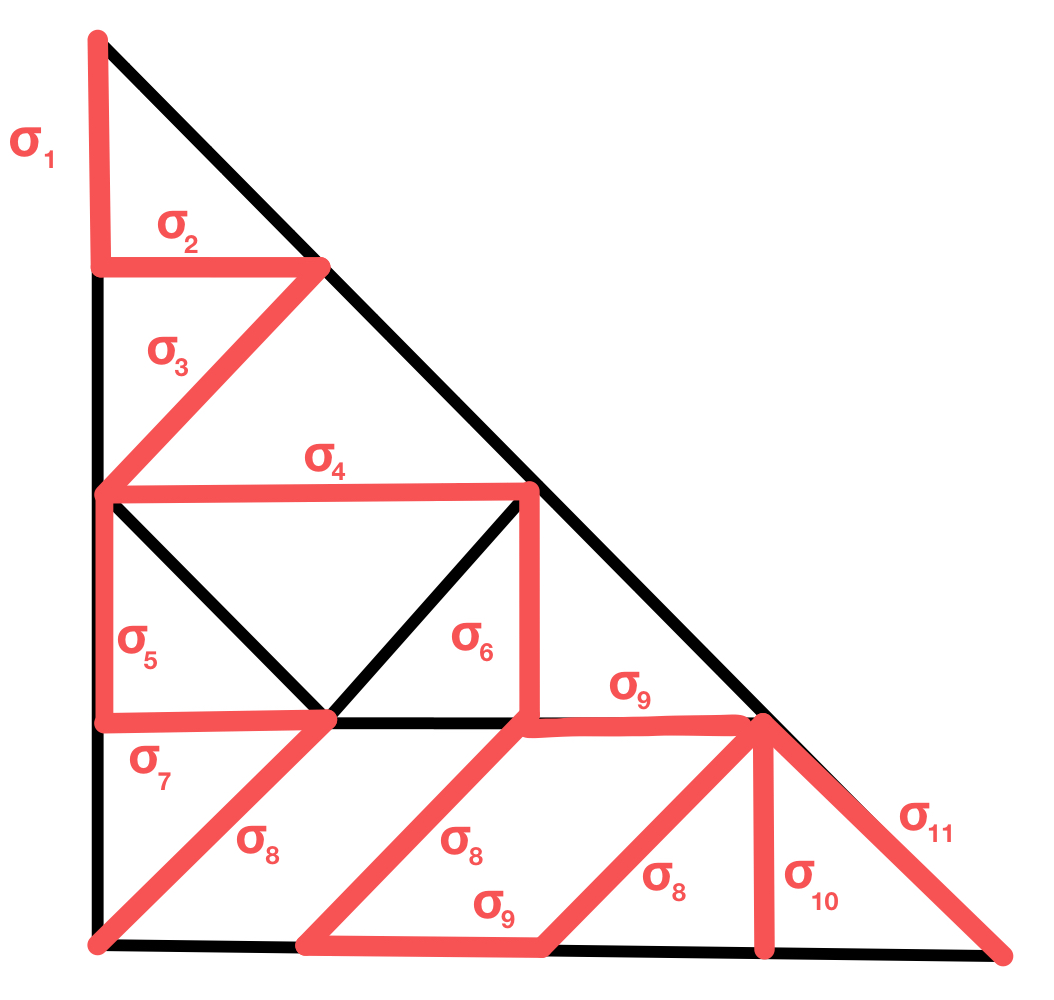}
    \end{tabular}
    \caption{Labelled dual subdivision and marked (extended) edges.\label{fig:marked edges}}
    \end{figure}
        In Figure \ref{fig:marked edges}
    we color edges dual to an edge with a marking $x_i$ and call them $\sigma_i$, for $i=1,\ldots,11$. Here, we give equivalent edges the same name.

    By Step 2 we can set $a_{40}=1$, and from this we get unique coefficients $a_{ij}\in k^\times$, for all vertices $(i,j)\in \left(\bigcup_{l=1}^{11} \sigma_l\right)\cap\Z^2$.
    
    In Figure \ref{fig:reconstruction} we illustrate how one can find all possibilities for the remaining $a_{ij}$ step by step as in step 4. First we find all $a_{ij}$ for the triangles $\Delta_l$ with two marked edges $\sigma_l$ and~$\sigma_l'$ (these are exactly the triangles dual to the $2$-valent vertices of $G_0$ in Example \ref{ex:graphconstruction}). By \cite[Lemma 3.5]{Shustin} there are 
    \begin{equation}
    \label{eq:reconstructionexample1}
    \frac{\vert \Delta_1\vert}{\vert \sigma_1\vert\vert\sigma_2\vert}\cdot \frac{\vert \Delta_2\vert}{\vert \sigma_2\vert\vert\sigma_3\vert}\cdot \frac{\vert \Delta_3\vert}{\vert \sigma_3\vert\vert\sigma_4\vert}\cdot \frac{\vert \Delta_4\vert}{\vert \sigma_5\vert\vert\sigma_7\vert}\cdot \frac{\vert \Delta_7\vert}{\vert \sigma_6\vert\vert\sigma_9\vert}\cdot\frac{\vert \Delta_8\vert}{\vert \sigma_7\vert\vert\sigma_8\vert}\cdot \frac{\vert \Delta_{11}\vert}{\vert \sigma_8\vert\vert\sigma_{10}\vert}\cdot \frac{\vert \Delta_{12}\vert}{\vert \sigma_{10}\vert\vert\sigma_{11}\vert}\end{equation}
    possibilities for the coefficients $a_{ij}$ in the red shaded area in the first picture in Figure \ref{fig:reconstruction}. In particular, we know all possibilities for the choices of $a_{ij}$ for $(i,j)\in \tau_1$, the blue edge in the second picture in Figure \ref{fig:reconstruction}. Hence, we can use \cite[Lemma 3.5]{Shustin} and get 
    \begin{equation}
    \label{eq:reconstructionexmple2}
    \frac{\vert \Delta_5\vert}{\vert \sigma_4\vert\vert\tau_1\vert}\end{equation}
    possible choices for the~$a_{ij}$ with $(i,j)\in \Delta_{5}$ (note that $\Delta_5$ is the only triangle in $S$ dual to a~$2$-valent vertex of $G_1$ which was not a $2$-valent vertex of $G_0$ in Example \ref{ex:graphconstruction}), illustrated by the blue shaded triangle in the second picture. Lastly, we find all possibilities for the~$a_{ij}$ with $(i,j)\in \Delta_6$ (also here $\Delta_6$ is the only triangle of $S$ which is dual to a $2$-valent vertex of~$G_2$ which was not $2$-valent in $G_1$), illustrated by the purple shaded triangle on the right. There are 
     \begin{equation}
     \label{eq:reconstructionexample3}
     \frac{\vert \Delta_6\vert}{\vert \sigma_6\vert\vert\tau_2\vert}\end{equation}
     possibilities.
     There is a unique choice for the $a_{ij}$ where the $(i,j)\in \Delta_{9}\cup \Delta_{10}$, i.e., lies on a parallelogram, since these are uniquely determined by the coefficients at three of their vertices.
    \begin{figure}
    \begin{tabular}{ccc}
    \includegraphics[scale=0.13]{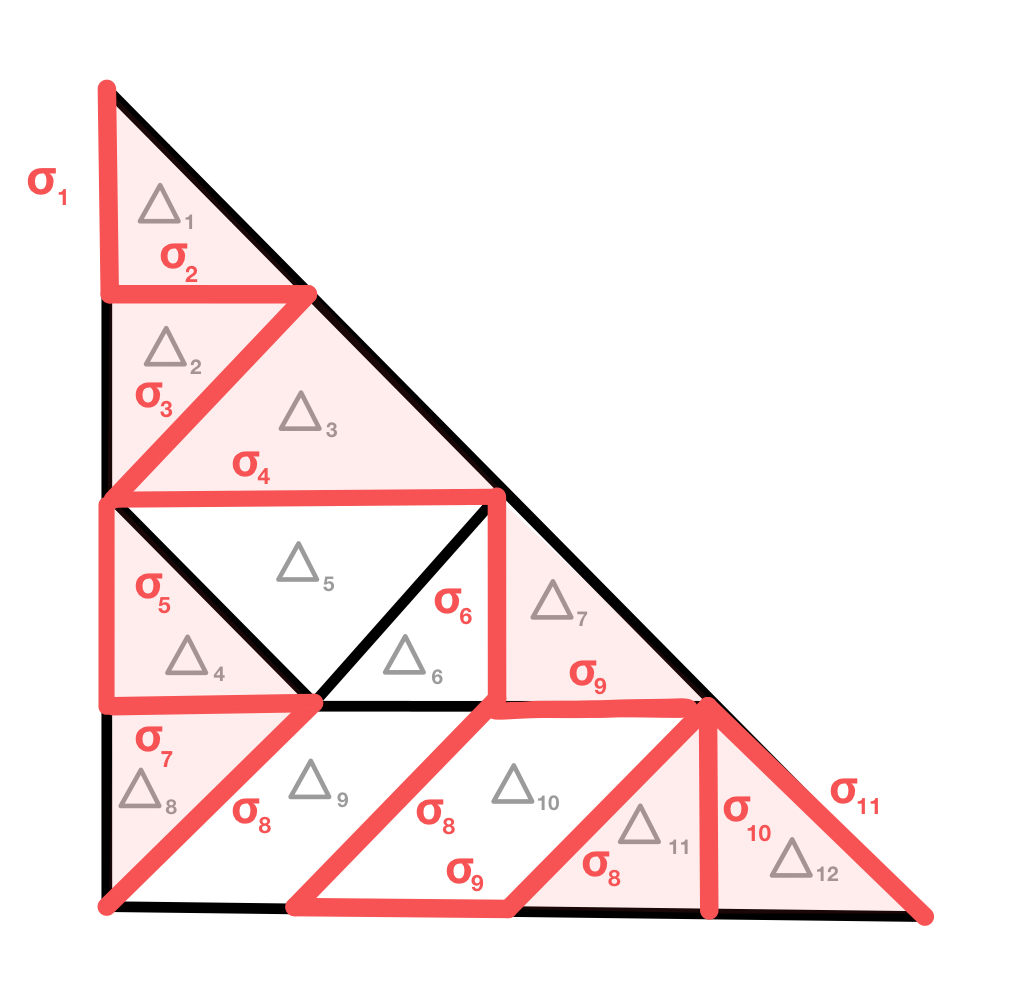}&
    \includegraphics[scale=0.13]{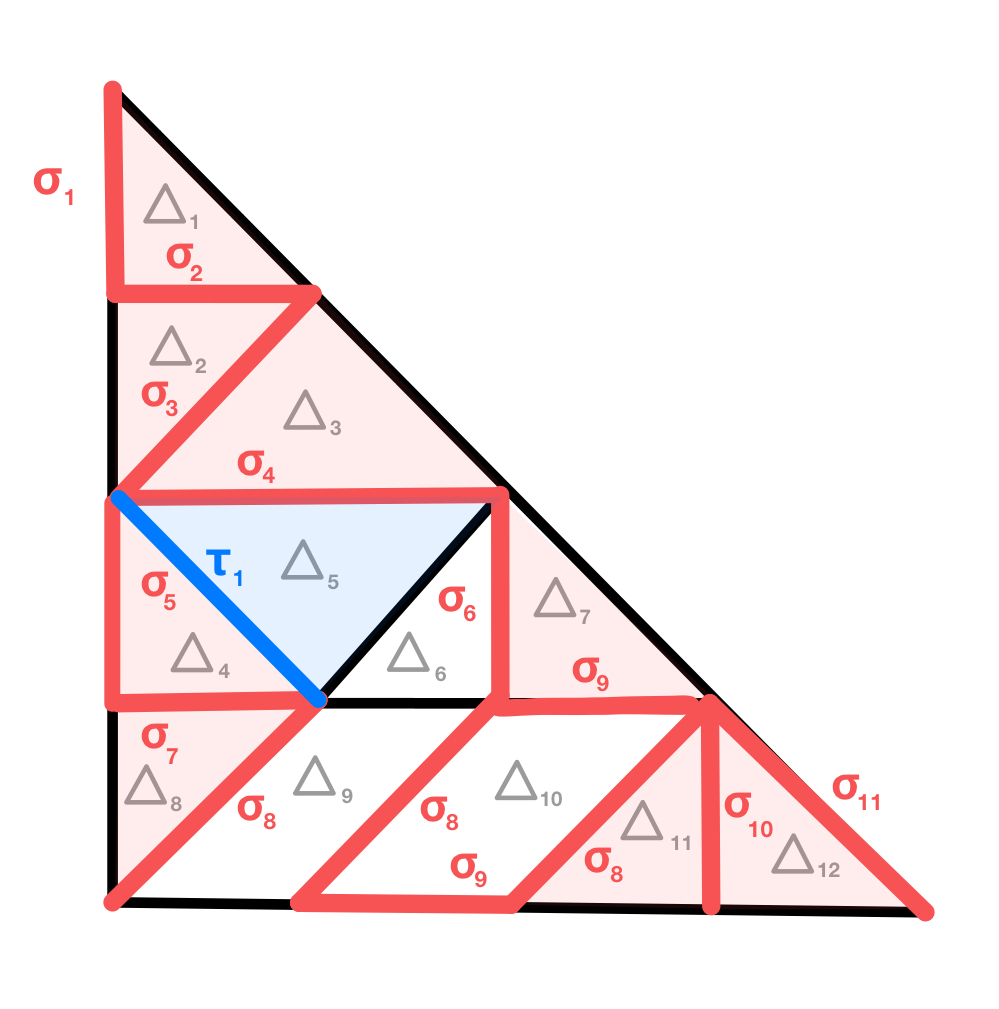}&
    \includegraphics[scale=0.13]{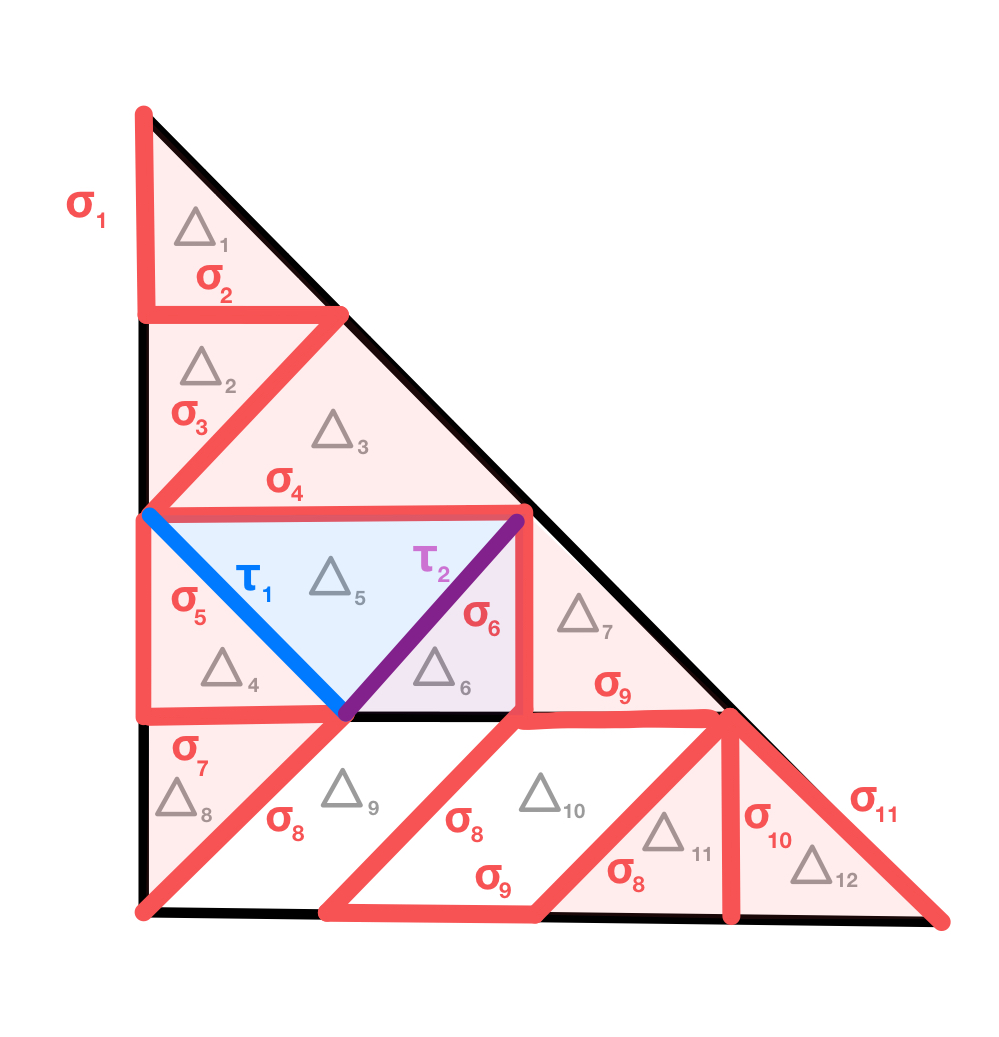}
    \end{tabular}
    \caption{Finding the coefficients $a_{ij}$ step by step.\label{fig:reconstruction}}
    \end{figure}
    In this example each of the factors in \eqref{eq:reconstructionexample1}, \eqref{eq:reconstructionexmple2} and \eqref{eq:reconstructionexample3} equals $1$ and thus the $k$-algebra $F$ from step 4 equals $F=k$.
    There is only one edge in $S$ of lattice length greater than $1$, namely $\sigma_4$ with $\vert\sigma_4\vert=2$.
    So the $F$-algebra $L$ in step 5 defined by the different choices for the deformation patterns is $L=\frac{F[u_4]}{(u_4^2-a_4c_4)}$.

    Finally, the $L$-algebra $\etale$ from step 6 equals $\etale=\frac{L[w_4]}{(w_4^2-d_4)}$ for some $d=d_4\neq 0$.
    Note that $\dim_k\etale=4=\operatorname{mult}_\C(A)$.

    \end{ex}

\subsection{Proof of the quadratically enriched correspondence theorem}
Let $C$ be the unique curve defined by $f\in \etale\Puiseux[x,y]$ we found in Subsection \ref{subsection:restoringthecurve}. Its quadratic weight $\Wel^{\A^1}_{\etale\Puiseux}(C)$ lives in $\GW(\etale\Puiseux)$. 
To prove our quadratically enriched correspondence theorem we need to compute 
\begin{equation}
\label{eq:TrEk}
    \Tr_{\etale\Puiseux/k\Puiseux}(\Wel^{\A^1}_{\etale\Puiseux}(C))\in \GW(k\Puiseux).
\end{equation}

We first compute $\Wel^{\A^1}_{\etale\Puiseux}(C)$. By the refined patchworking theorem, the nodes of~$C$ are in one to one correspondence with the nodes of the $C_l$ defined by $f_1,\ldots,f_N$ and the nodes of the deformation patterns and this correspondence preserves the type of node under the natural isomorphism $\GW(\etale\Puiseux)\cong \GW(\etale)$.
We have that
\begin{itemize}
    \item for a triangle $\Delta_l$ in $S$ with edges $\sigma_l$, $\sigma'_l$ and $\sigma''_l$,
        \begin{equation}
        \label{eq:Wel(C_l)}
          \Wel^{\A^1}_F(C_l)=\qinv{ (-1)^{\operatorname{Int}(\Delta_l)}\frac{\vert \Delta_l\vert^{\vert \Delta_l\vert}}{\vert\sigma_l\vert^{\vert\sigma_l\vert}\cdot\vert \sigma'_l\vert^{\vert \sigma'_l\vert} \cdot\vert\sigma''_l\vert^{\vert\sigma''_l\vert}}} \in \GW(F)
       \end{equation}
       where $\operatorname{Int}(\Delta_l)$ is the number of interior lattice points of $\Delta_l$ by \eqref{eq:quadraticweightofC'triangle}.
        In this case, all edges of $\Delta_l$ are odd and thus by Pick's theorem $\vert\Delta_l\vert$ is also odd, so  \eqref{eq:Wel(C_l)} equals
        \[\Wel^{\A^1}_F(C_l)=\qinv{(-1)^{\operatorname{Int}(\Delta_l)}\frac{\vert \Delta_l\vert}{\vert\sigma_l\vert\cdot\vert \sigma'_l\vert \cdot\vert\sigma''_l\vert}} \in \GW(F)\]
    \item for a parallelogram $\Delta_l$ in $S$, 
    \begin{equation}
        \Wel^{\A^1}_F(C_l)=\qinv{1} \in \GW(F)
    \end{equation}
    by Proposition \ref{thm:parallelogram}.
    \item Let $\sigma\in E^*(A)$ be an extended edge. Then a deformation pattern $f_{[\sigma]}$ defining the curve $C_{[\sigma]}$ has quadratic weight
    \begin{equation}
        \Wel^{\A^1}_{L_{[\sigma]}}(C_{[\sigma]})=\begin{cases}\qinv{ 1}\in \GW(L_{[\sigma]}) & \text{if $\vert \sigma\vert$ is odd}\\
        \qinv{ -u_{[\sigma]}}\in \GW(L_{[\sigma]}) & \text{if $\vert \sigma\vert$ is even}
        
        \end{cases} 
    \end{equation}
    by Proposition \ref{thm:quadratic 3.9} and where $L_{[\sigma]}$ is defined as in \eqref{eq:Lsigma}.
\end{itemize}
Multiplying the above quadratic weights and observing that every interior edge of $S$ appears twice in the product of the \eqref{eq:Wel(C_l)}'s over all triangles in $S$ and all boundary edges of~$S$ have lattice length $1$ yield
\begin{equation}
\label{eq:Wel(C)}
\Wel^{\A^1}_\etale(C)= \qinv{ \mkern-9mu \prod_{\quad\text{triangles $\Delta'$ in $S$}} \mkern-36mu (-1)^{\operatorname{Int}(\Delta')}\vert\Delta'\vert^{\vert\Delta'\vert} \cdot \mkern-9mu \prod_{[\sigma]\in E^*(A)\text{ even}}\mkern-18mu (-u_\sigma)}\in \GW(\etale),
\end{equation}
 where $\operatorname{Int}(\Delta')$ is the number of interior points of the triangle $\Delta'$ in $S$.
Now we can compute~\eqref{eq:TrEk} in $\W(k)$. 

Let us start with the case when all edges of $S$ are odd. Note that in this case the quadratic weight $\Wel^{\A^1}_\etale(C)=\qinv{ \prod_{\Delta'\text{ triangle in }S} (-1)^{\operatorname{Int}(\Delta')}\vert\Delta'\vert}$. Let us denote by $\operatorname{Int}(S)$ the total number of interior points of all the triangles in $S$.
Then,
\[\Wel^{\A^1}_{\etale\Puiseux}(C)=\qinv{(-1)^{\operatorname{Int}(S)} \prod_{\Delta'\text{ triangle in }S}\vert\Delta'\vert}\in \GW(\etale\Puiseux).\]
Denote by $\Wel_\etale^{\A^1}(C)= \qinv{(-1)^{\operatorname{Int}(S)} \prod_{\Delta'\text{ triangle in }S}\vert\Delta'\vert}\in \GW(\etale)$ the image of the quadratic weight $\Wel^{\A^1}_{\etale\Puiseux}(C)$ in $\GW(\etale\Puiseux)$ under the canonical isomorphism between $\GW(\etale\Puiseux)$ and $\GW(\etale)$ from Remark \ref{rmk:GW(K)}.
We compute \[\Tra{\Wel^{\A^1}_\etale(C)}=\Tr_{F/k}\circ\Tr_{L/F}\circ \Tra[\etale/L]{\Wel^{\A^1}_\etale(C)}.\]
In order to compute $\Tra[\etale/L]{\Wel^{\A^1}_\etale(C)}$,
recall from \eqref{eq:E} that $\etale=\frac{L[w_1,\ldots,w_n]}{\left(w_1^{\vert\sigma_1\vert}-d_1,\ldots,w_n^{\vert\sigma_n\vert}-d_n\right)}$. 
Set $L_0=L$ and $L_i=\frac{L_{i-1}[w_{i}]}{\left(w_i^{\vert\sigma_i\vert}-d_i\right)}$ for $i=1,\ldots,n$. Then we have that $\etale=L_n$, the dimension $\dim_{L_{i-1}}L_i=\vert\sigma_i\vert $ is odd, and  $\Tr_{\etale/L}=\Tr_{L_1/L_0}\circ\ldots\circ \Tr_{L_n/L_{n-1}}$.
We apply Proposition~\ref{prop:traces}(1) $n$ times in the odd case and obtain
\begin{align*}
    \Tra[\etale/L]{\Wel^{\A^1}_\etale(C)}
    =&\Tra[\etale/L]{\qinv{(-1)^{\operatorname{Int}(S)}\prod_{\text{triangles $\Delta'$ in $S$}}\vert\Delta'\vert}}\\
    =&\Tr_{L_1/L_0}\circ\ldots \circ \Tr_{L_n/L_{n-1}}\left({\qinv{(-1)^{\operatorname{Int}(S)}\prod_{\text{triangles $\Delta'$ in $S$}}\vert\Delta'\vert}}\right)\\
    \overset{\ref{prop:traces}}{=}&\left(\qinv{(-1)^{\operatorname{Int}(S)}\prod_{\text{triangles $\Delta'$ in $S$}}\vert\Delta'\vert\cdot \prod_{i=1}^n\vert\sigma_i\vert}\right)\\
    =&\left(\qinv{(-1)^{\operatorname{Int}(S)}\prod_{\text{triangles $\Delta'$ in $S$}}\vert\Delta'\vert\cdot \prod_{[\sigma]\in E^*_\text{marked}(A)}\vert\sigma\vert}\right).
\end{align*}
Next, we apply $\Tr_{L/F}$. Recall from \eqref{eq:Lsigma} and \eqref{eq:L} that 
\[L=\frac{F\left[\{u_{\sigma}:[\sigma]\in E^*(A)\}\right]}{\left(\{u_{\sigma}^{\vert\sigma\vert}-1:\text{$[\sigma]\in E^*(A)$}\}\right)}.\]
Again by applying Proposition \ref{prop:traces}(1) (odd case) successively, we get 
\begin{align*}
    &\Tr_{L/F}\circ \Tra[\etale/L]{\Wel^{\A^1}_\etale(C)}\\
    =& \Tr_{L/F}\left(\qinv{(-1)^{\operatorname{Int}(S)}\prod_{\text{triangles $\Delta'$ in $S$}}\vert\Delta'\vert\cdot \prod_{[\sigma]\in E^*_\text{marked}(A)}\vert\sigma\vert}\right)\\
    \overset{\ref{prop:traces}}{=}&\qinv{(-1)^{\operatorname{Int}(S)}\prod_{\text{triangles $\Delta'$ in $S$}}\vert\Delta'\vert\cdot \prod_{[\sigma]\in E^*_\text{marked}(A)}\vert\sigma\vert\cdot\prod_{[\sigma]\in E^*(A)}\vert\sigma\vert}.
\end{align*}
Finally, we have to apply $\Tr_{F/k}$. For this, recall that we found $F$ by successively adjoining~$\frac{\vert\Delta'\vert}{\vert\sigma\vert\vert\sigma'\vert}$ many roots for all triangles $\Delta'$ in $S$ (here $\sigma$ and $\sigma'$ are two sides of the triangle). Since we assume all edges are odd, by Pick's theorem \ref{thm:Pick} $\vert\Delta'\vert$ and thus also $\frac{ \vert\Delta'\vert}{\vert\sigma\vert\vert\sigma'\vert}$ are odd. So again by applying Proposition \ref{prop:traces} for each triangle, we get
\begin{align*}
    &\Tr_{F/k}\circ\Tr_{L/F}\circ \Tra[\etale/L]{\Wel^{\A^1}_\etale(C)}\\
    =&\Tr_{F/k}\left(\qinv{(-1)^{\operatorname{Int}(S)}\prod_{\text{triangles $\Delta'$ in $S$}}\vert\Delta'\vert\cdot \prod_{[\sigma]\in E^*_\text{marked}(A)}\vert\sigma\vert\cdot\prod_{[\sigma]\in E^*(A)}\vert\sigma\vert}\right)\\
    \overset{\ref{prop:traces}}{=}&\qinv{(-1)^{\operatorname{Int}(S)}\prod_{\text{triangles $\Delta'$ in $S$}}\vert\Delta'\vert\cdot \prod_{[\sigma]\in E^*_\text{marked}(A)}\vert\sigma\vert\cdot\prod_{[\sigma]\in E^*(A)}\vert\sigma\vert\cdot \prod_{\text{triangles $\Delta'$ in $S$}}\frac{\vert\Delta'\vert}{\vert\sigma\vert\vert\sigma'\vert}}\\
    =&\qinv{(-1)^{\operatorname{Int}(S)}\prod_{\text{triangles $\Delta'$ in $S$}}\vert\Delta'\vert^2}
    =\qinv{(-1)^{\operatorname{Int}(S)}}
\end{align*}
in $\W(k)$. In the last product in the second last row we divide by the lengths of the two edges where we know the coefficients (see Subsection \ref{subsection:lemma3.5} and step 3 and 4 in Subsection~\ref{subsection:restoringthecurve}). In this product, every marked extended edge not on the boundary appears exactly twice and any unmarked extended edge appears once. Thus, all the edge lengths cancel and the double areas come in squares (compare the second last row). The last row follows since  one has that $\qinv{ab^2}=\qinv{a}$ in $\GW(k)$.
Since the rank of $\Tr_{\etale/k}(\Wel_\etale^{\A^1}(C))$ equals 
\[\dim_k\etale=\prod_{\Delta'\text{ triangle in }S}\vert\Delta'\vert\]
and an element of $\GW(k)$ is completely determined by its rank and its image in $\W(k)$ we conclude that 
\[\Tr_{\etale/k}(\Wel^{\A^1}_\etale(C))=\mult_k(A)\]
in $\GW(k)$, in case all edges of $S$ are odd.

In case $S$ contains an even edge $\sigma$, then $\Tr_{L_\sigma/F}(\qinv{ u_\sigma})$ is a multiple of the hyperbolic form $h$  (see \eqref{eq:defpatterntraceeven}) and thus equals $0$ in $\W(k)$. It follows that the image of $\Tr_{\etale/k}(\Wel_\etale^{\A^1}(C))$ in $\W(k)$ is $0$. Consequently, $\Tr_{\etale/k}$ is a multiple of $h$ of rank $\dim_k\etale$, that is.
\[\Tr_{\etale/k}(\Wel^{\A^1}_\etale(C))=\frac{\dim_k\etale}{2}h=\mult_k(A)\]
in $\GW(k)$.
This completes the proof of the quadratically enriched correspondence theorem.

\subsection{Quadratically enriched vertex multiplicity}
Let $v$ be a vertex of $A$ with dual triangle $\Delta_v$ with edges $\sigma_v$, $\sigma_v'$ and $\sigma_v''$ in the dual subdivision of $A$.
We define the \emph{quadratically enriched vertex multiplicity} of $v$ to be
\begin{equation}
    m_v^{\A^1}\coloneqq \begin{cases}
        \qinv{ (-1)^{\operatorname{Int}(\Delta_v)}\vert\sigma_v\vert\vert\sigma_v'\vert\vert\sigma_v''\vert}+\frac{\vert\Delta_v\vert-1}{2}h &\text{ if $\sigma_v$, $\sigma_v'$ and $\sigma_v''$ are all odd}\\
        \frac{\vert \Delta_v\vert}{2}h& \text{else}
    \end{cases}
\end{equation}
where $\operatorname{Int}(\Delta_v)$ is the number of interior lattice points of $\Delta_v$.

For a simple tropical curve $A$ the quadratically enriched multiplicity $\operatorname{mult}_k^{\A^1}(A)$ of $A$ equals the product of its quadratically enriched vertex multiplicities:
\[\operatorname{mult}_k^{\A^1}(A)=\prod_{\text{$v$ a 3-valent vertex of $A$}}m_v^{\A^1}.\]
This is because the edges of the dual subdivision $S$ of $A$ which lie on the boundary $\partial \Delta$ all have length $1$. The edges in $S$ not on the boundary can be of length $\vert\sigma\vert>1$. However, each of these extended edges is an edge of exactly $2$ triangles in $S$. Thus in the case that all edges of $S$ are odd, when multiplying the vertex multiplicities $\vert\sigma\vert$ will appear twice in the non-hyperbolic term, that is as a square and in $\GW(k)$ it holds that $\qinv{ a}=\qinv{ ab^2}$.

In case, there is an even edge, the product of the vertex multiplicities will be hyperbolic, since multiplying with the hyperbolic always equals an integral multiple of the hyperbolic form.

The following is a motivation for the definition of the quadratically enriched vertex multiplicity. 
Assume $A$ is a tropical curve with Newton polygon a triangle and only one vertex $v$, i.e. $\Delta_v=\Delta$ is a triangle. 
Then $A$ has $3$ edges which are all unbounded. Assume two of them are marked by a point and assume that $\sigma_v$ and $\sigma_v'$ are the edges of $\Delta_v$ dual to these marked edges. Note that we do not assume that the edges of $\Delta_v$ have lattice length~$1$, that is, the curve $A$ is not necessarily a simple curve. If one follows the recipe to reconstruct all $f$ defining curves tropicalizing to $A$ described in Subsection \ref{subsection:restoringthecurve}, then the quadratically enriched count of curves tropicalizing to $A$ equals $m_v^{\A^1}$ in case all three edges are odd:
We again compute in $\W(k)$. 
\begin{align*}
    \Tr_{F/k}\circ\Tr_{\etale/F}\qinv{ (-1)^{\operatorname{Int}(\Delta_v)}\frac{\vert \Delta_v \vert}{\vert \sigma_v\vert\vert \sigma_v\vert\vert \sigma_v''\vert}}&=\Tr_{F/k}\qinv{ (-1)^{\operatorname{Int}(\Delta_v)}\frac{\vert \Delta_v\vert }{\vert \sigma_v''\vert}}\\
     &=\qinv{(-1)^{\operatorname{Int}(\Delta_v)} \frac{\vert\Delta_v\vert^2}{\vert \sigma_v\vert\vert \sigma_v\vert\vert \sigma_v''\vert}}\\
     &=\qinv{ (-1)^{\operatorname{Int}(\Delta_v)}\vert \sigma_v\vert\vert \sigma_v\vert\vert \sigma_v''\vert}
\end{align*}
where $F=\frac{k[w]}{\left(w^{\frac{\vert\Delta_v\vert}{\vert\sigma_v\vert\vert\sigma_v'\vert}}-D\right)}$
for some $D\in k^{\times}$ and $\etale=\frac{F[x,x']}{\left(x^{\vert\sigma_v\vert}-d, x'^{\vert\sigma_v'\vert}-d'\right)}$ for some $d,d'\in F^\times$.

When one of the edges of $\Delta_v$ is even, we define the quadratically enriched multiplicity to be a multiple of $h$. This ensures that the quadratic weight $\operatorname{mult}_k^{\A^1}$ of a simple tropical curve is a multiple of the hyperbolic form $h$ in case the tropical curve has an even edge.


\end{document}